\documentclass[review]{elsarticle}

\usepackage{lineno,hyperref}
\modulolinenumbers[5]

\journal{Journal of Computational Physics}

\usepackage{amsmath}
\usepackage{amsfonts}
\usepackage{amssymb}
\usepackage{amsthm}
\usepackage{mathtools}
\usepackage{algorithm}
\usepackage{algpseudocode}
\usepackage{graphicx}
\usepackage{caption}
\usepackage{subcaption}
\usepackage{multicol}
\usepackage{multirow}
\usepackage{hhline}
\usepackage{float}
\usepackage{color}
\usepackage{xcolor}
\usepackage{chngcntr}
\usepackage{tikz}
\usepackage{array}
\usepackage{comment}
\usepackage{geometry}
\usepackage{bm}
\usepackage{epsfig}

\usepackage{latexsym}
\usepackage{enumerate}
\usepackage{hyperref}
\usepackage{cleveref}
\usepackage{booktabs}

\graphicspath{{./figs/}}

\newcommand{\ds}{\displaystyle}

\newcommand{\bbeta}{{\boldsymbol\eta}}

\newcommand{\bPsi}{{\boldsymbol\Psi}}

\newcommand{\bzeta}{{\boldsymbol\zeta}}

\newcommand{\m}{{\mathbf{m}}}

\newcommand{\bx}{\mathbf{x}}
\newcommand{\bv}{\mathbf{v}}

\newcommand{\0}{{\mathbf{0}}}

\def\by{\mathbf{y}}

\newcommand{\bbR}{\mathbb{R}}
\newcommand{\bbZ}{\mathbb{Z}}

\newcommand{\cN}{\mathcal{N}}

\newcommand{\cM}{\mathcal{M}}
\newcommand{\cT}{\mathcal{T}}
\newcommand{\cO}{\mathcal{O}}
\def\wt{\widetilde}
\def\wh{\widehat}

\numberwithin{equation}{section}
\allowdisplaybreaks

\newtheorem{remark}{Remark}[section]

\usepackage[normalem]{ulem}

\begin{document}

\begin{frontmatter}

\title{Structurally informed data assimilation in two dimensions}

\author[UMBC]{Tongtong Li\fnref{myfootnote}\corref{mycorrespondingauthor}}
\cortext[mycorrespondingauthor]{Corresponding author}
\ead{tongtong.li@umbc.edu}

\author[dartmouth]{Anne Gelb\fnref{myfootnote}}
\ead{annegelb@math.dartmouth.edu}

\author[dartmouth]{Yoonsang Lee\fnref{myfootnote}}
\ead{yoonsang.lee@dartmouth.edu}

\address[UMBC]{Department of Mathematics and Statistics, University of Maryland, Baltimore County, Baltimore, MD 21250, USA}

\address[dartmouth]{Department of Mathematics, Dartmouth College, Hanover, NH 03755, USA}

\fntext[myfootnote]{This work is partially supported by the NSF grant DMS \#1912685, DOE ASCR \#DE-ACO5-000R22725, and DOD ONR MURI grant \#N00014-20-1-2595.}

\begin{abstract}
Accurate data assimilation (DA) for systems with piecewise-smooth or discontinuous state variables remains a significant challenge, as conventional covariance-based ensemble Kalman filter approaches often fail to effectively balance observations and model information near sharp features. In this paper we develop a structurally informed DA framework using ensemble transform Kalman filtering (ETKF). Our approach introduces gradient-based weighting matrices constructed from finite difference statistics of the forecast ensemble, thereby allowing the assimilation process to dynamically adjust the influence of observations and prior estimates according to local roughness. The design is intentionally flexible so that it can be suitably refined for sparse data environments. Numerical experiments demonstrate that our new structurally informed data assimilation framework consistently yields greater accuracy  when compared to  more conventional approaches.
\end{abstract}

\begin{keyword}
Data assimilation, ensemble transform Kalman filtering, structurally informed prior
\end{keyword}

\end{frontmatter}


\section{Introduction}
\label{sec:introduction}
Data assimilation (DA) provides a systematic approach for combining incomplete and inaccurate information from models and observations, enabling more accurate estimates of the system state \cite{Carrassi18}. It is widely used across various scientific domains, notably atmospheric sciences, geoscience, and oceanography, for applications such as climate modeling and weather forecasting. For instance, DA is used in operational weather forecast prediction centers, where the partial differential equations (PDEs) that describe atmospheric evolution are integrated to produce short-to-medium range forecasts \cite{Kalnay02}. In geoscience applications, DA helps address challenges induced by the chaotic nature of the atmosphere and oceans \cite{Carrassi08, NG11}. DA has also been widely used in modeling chemical constituents of the atmosphere \cite{Bocquet2015}, where the fundamental multivariate characteristic in chemistry and transport models is one of the key aspects to be accounted for. Despite its widespread use, traditional DA methods often overlook the structural characteristics of the underlying variables in these applications.

With the increasing availability of data from new measurement devices and advanced image processing techniques, DA has also gained  popularity in a broader range of applications, including geophysical fields such as sea ice modeling \cite{Mu18, Sakov12, Xie18, Sampson21, Cheng23}. In these applications, the sea ice state variable often contains sharp features in the sea ice cover, such as leads and ridges. It is therefore imperative that assimilation schemes used for such systems must not only preserve the physical principles of the underlying dynamics but also accurately recover solutions with discontinuous profiles, even when these features are not directly observed. 

Discrete-time filtering, from a probabilistic perspective, is designed to sequentially update the probability distribution of the state variable, conditioned on the accumulated data up to the current time. It is the most common and practical technique used in in geophysical applications. The state variable is treated as random, evolving from an assumed probabilistic distribution. When this distribution is Gaussian, and the system dynamics are linear with additive, state-independent Gaussian model and observation errors, \emph{Kalman filtering} \cite{Law15} provides an analytic solution, in the sense that the posterior distributions remain Gaussian, and their means and covariances can be computed in closed form. Many problems of interest generally depart from these assumptions, however, as they involve non-linearity and non-Gaussian behavior.

Such non-linearity and non-Gaussian behavior has motivated the development of ensemble-based Kalman filters, which include the ensemble Kalman filter (EnKF) \cite{Evensen94, Evensen03}, the ensemble adjustment Kalman filter (EAKF) \cite{Anderson01}, and the ensemble transform Kalman filter (ETKF) \cite{Bishop01, Tippett03}. These methods employ a two-step Monte Carlo approximation of the classical Kalman recursions on ensembles \cite{Evensen09}. In the \emph{forecast} step, a particle approximation of the filtering distribution is propagated through the transition dynamics to yield a ``forecast" ensemble for the next observation time. In the \emph{analysis} step, the forecast ensemble is updated through a linear transformation, producing an empirical estimate of the new filtering marginal. Once again, though, this linear transformation is estimated under Gaussian assumptions, and hence cannot consistently characterize the posterior distribution for non-Gaussian models \cite{Mandel11}. Furthermore, the update neglects structural features of the underlying state variables that could improve the characterization of errors, which often stem from model deficiencies and other sources. 

Specific prior knowledge about the underlying variable has been incorporated in variational DA approaches through the use of regularization terms. For instance, $\ell_1$ norm regularization was utilized in the analysis step of the method described in \cite{Freitag10, Budd11} to ensure sharp features are retained.  Similar $\ell_1$ norm regularization terms are used in image and signal recovery methods, as they are known to help preserve the sparsity in the corresponding sparse domain of the underlying signal or image (see e.g. \cite{CT}).  Analogous $l_1$ regularization methods have been adapted to promote sparse representation of hydro-meteorological states in the derivative or wavelet domains \cite{Foufoula-Georgiou14, Ebtehaj13}. Finally, we note that adding  an  $l_1$ norm regularization term  to the more standard $l_2$ objective function, that is, a mixed variational regularization approach, may lead to improved performance in certain applications \cite{Freitag13, Ebtehaj14}. Importantly, none of these methods incorporate  \emph{local} information regarding the sparse domain \cite{Adcock19,GelbScarnati2019}, which is critical for capturing fine-scale structures in complex geophysical systems.

By contrast, the structurally-informed data assimilation  method introduced in \cite{TLda}   explicitly incorporates local information, e.g.~the location of shock discontinuities,  into the assimilation process to improve performance when the underlying state variables exhibit such structures, and as such was partially inspired by the aforementioned variational DA approaches, specifically with respect to the equivalence between statistical and variational approaches. In particular the analysis step, where the prior distribution is updated to the posterior, can be understood as an optimization problem with the posterior mean being equivalent to the maximum a posterior (MAP) estimate for Gaussian distributions. From this perspective the analysis step solves an inverse problem with an objective function consisting of two components: one that quantifies the misfit between the observations and model state, or the {\em likelihood}, and the other that serves as an $\ell_2$ regularization term, which accounts for the difference between the model state and the prior ensemble mean from the forecast step, weighted by the covariance of the prior ensembles. 

As shown in \cite{TLda} for one-dimensional problems,  the prior sample covariance fails to capture local structure information, consequently making it less effective for state variables with discontinuous profiles. By contrast, the second moment of the gradient of the state variable is able to not only identify the discontinuous regions, but also distinguish shocks from steep gradients based on their relative magnitudes. Such structural information regarding underlying variables can then be used to better balance model and data information in the form of a weighting matrix ($W$ in \eqref{eq:ETRF min}). The approach can be enhanced by refining $W$ so that the solution at each domain grid point is selectively influenced by neighboring values.  Finally, the computational method used for advancing the PDE model was carefully chosen to ensure high-order convergence properties and resolved jump discontinuities.  

Extending this framework to two dimensions introduces new challenges.  In particular, discontinuities can form complex geometries so that directional information becomes critical, yet sparse observation networks make estimating reliable spatial correlations more difficult. This investigation seeks to develop a structurally-informed data assimilation framework for two-dimensional systems characterized by piecewise-smooth or discontinuous states. The inherent challenges are addressed by first formalizing the characterization of smooth and discontinuous regions.  Specifically,  we identify “smooth regions” as areas where the solution exhibits small numerical gradients, that is, reflecting locally continuous behavior, while “discontinuous regions” are conversely areas with large gradient, thereby indicating sharp transitions or jumps. Once this is successfully accomplished, we are then able to replace the conventional covariance-based weighting matrices in ensemble-based Kalman filters with matrices constructed from local gradient statistics of the forecast ensemble. As before, this enables the assimilation process to dynamically assign greater weight to observational data in discontinuous regions where model forecasts are unreliable, and rely more on the forecast in smooth regions. The approach is validated through numerical experiments on benchmark problems, including a linear advection equation and Burgers’ equation.  Our method improves assimilation accuracy and preserves key structural features more effectively than classical ensemble-based methods in both fully and sparsely observed settings.

The remainder of the paper is organized as follows.   Section \ref{sec:preliminaries} provides mathematical conventions used in this investigation and reviews the ETKF method. We present  our new structurally informed weighting matrix design for both fully and sparsely observed cases in Section \ref{sec:method}. Numerical experiments are conducted and analyzed in  Section \ref{sec:numerical}. Section \ref{sec:summary} provides some concluding remarks and ideas for future investigations.

\section{Preliminaries}
\label{sec:preliminaries}
We begin by establishing the mathematical foundations and notation used throughout the manuscript. To this end, in Section \ref{subsec:setup} we provide the general framework for sequential data assimilation.  We then review the ensemble transform Kalman filter (ETKF) in Section \ref{subsec: ETKF}, which serves as our prototype method. Finally we discuss the variance inflation and localization techniques and describe how they effectively mitigate potential challenges that arise when working with high-dimensional system finite ensembles  in Section \ref{subsec:inflate}. 

For purposes of consistency in presentation, we mainly follow the conventions used in \cite{Law15}. In this regard we denote $\bbZ^{+}$ as the set of non-negative integers, $C(X,Y)$ as the set of continuous functions mapping from space $X$ to space $Y$, and $\0$ as the null vector (or tensor as appropriate).  We also write the Hadamard (entrywise) product of two matrices $M:=(M_{i,j})_{i,j=1,\cdots,n}$ and $N=(N_{i,j})_{i,j=1,\cdots,n}$ as
$$M \odot N := (M_{i,j}N_{i,j})_{i,j=1,\cdots,n}.$$
Lastly, the squared weighted norm for a vector $\bv$ with respect to symmetric positive definite matrix $C$ is given by
\begin{equation}
\vert \bv \vert_C^2:= \vert C^{-\frac{1}{2}} \bv \vert^2 = \bv^T C^{-1} \bv.\label{eq:normsquared}
\end{equation}

\subsection{Problem set up}
\label{subsec:setup}
Let $v(\bx,t)$ denote the state variable of interest that describes the system on a spatial domain $\Omega \subset \bbR^2$ for a given spatial location $\bx$ and time $t \in [0,T]$. The spatial domain is discretized with spatial resolution $\Delta x$ and $\Delta y$ in the $x$ and $y$ directions, respectively. The temporal domain is discretized with the uniform time step $\Delta t = \frac{\#CFL}{\frac{1}{\Delta x} + \frac{1}{\Delta y}}$ to satisfy the Courant–Friedrichs–Lewy (CFL) stability condition \cite{Courant1928, Courant1967}. Denoting $\bv_l:=\bv(t_l)$ as the corresponding tensor at time instance $t_l = l\Delta t$, $l=0, \cdots ,L$, where $L = int(\frac{T}{\Delta t})  \in \bbZ^{+}$ is the final time step for simulation, we assume that the solution trajectory of $\bv$  follows the deterministic discrete model\footnote{For simplicity we redefine $T=L\Delta t$ to avoid fractional time steps.}
\begin{equation}\label{model: dynamic}
\bv_{l+1}=\bPsi(\bv_l), \quad l = 0,\dots, L-1.
\end{equation}
Here $\bPsi: \bbR^n \rightarrow \bbR^n$ is a deterministic discrete evolution operator and typically arises from a numerical discretization of an underlying continuous dynamical system governed by some differential equations. For ordinary differential equations (ODEs), $\bPsi$ is a temporal integration operator. For PDEs, temporal integration follows spatial domain discretization, for example using  finite difference or finite element methods. Importantly, $\bPsi$ is, by construction, a low rank approximation of the underlying dynamics. 

We further assume that the initial state $\bv_0$ follows a Gaussian distribution with mean $\m_0$ and covariance matrix $C_0$:
\begin{equation}\label{assump: initial Gaussian}
\bv_0 \sim \cN(\m_0, C_0). 
\end{equation}
Finally, we assume that the given observational data $\by \in \mathbb{R}^m$ are defined as 
\begin{equation}\label{model: observation}
\by_{l+1}=H\bv_{l+1}+\bbeta_{l+1}, \quad l = 0,\dots, L-1,
\end{equation}
where $H: \bbR^n \rightarrow \bbR^m$ is the linear observation operator, and $\bbeta_{l+1}$ represents the observation error, which is assumed to be an independent and identically distributed (i.i.d.) sequence, independent of $\bv_0$, with 
\begin{equation}\label{assump: noise Gaussian}
\bbeta_1 \sim \cN(\0,\Gamma).
\end{equation}
Here $\Gamma$ is assumed to be known and symmetric positive definite (SPD). In the sparse observation setting, $H$ is a sub-sampling operator selecting only a subset of the state grid points.

\subsection{Ensemble transform Kalman filter (ETKF)}
\label{subsec: ETKF}

When the system dynamics $\bPsi$ are linear and all associated distributions are Gaussian, the classical Kalman filtering method \cite{Kalman60} yields an optimal sequential update of the state probability distribution conditioned on data, as the posterior remains Gaussian and is entirely described by its mean and covariance \cite{Law15, Cohn97}. As already discussed in the introduction, the system dynamics are nonlinear in many real-world applications, however. Since the Gaussian assumptions no longer hold, prediction and uncertainty propagation pose a challenging computational problem. Various ensemble-based Kalman filter methods have been developed as practical alternatives, including but not limited to the Ensemble Kalman Filter (EnKF) \cite{Evensen94, Evensen03}, the Ensemble Adjustment Kalman Filter (EAKF) \cite{Anderson01}, and the Ensemble Transform Kalman Filter (ETKF) \cite{Bishop01, Tippett03}. All of these methods sequentially estimate the state distribution by propagating an ensemble of particles through the model and then update the ensembles using observational information. As it serves as the prototype for our investigation, below we briefly describe the ETKF method.

The ETKF method is realized in two steps: forecast and analysis. In the forecast step, the particle approximations $\{ \bv^{(k)}_l\}_{k=1}^{K}$ of the filtering distribution are propagated through the transition dynamic model in \eqref{model: dynamic} to yield the forecast ensembles $\{ \wh{\bv}^{(k)}_{l+1} \}_{k=1}^{K}$ at the next observation time $t_{l+1}$:
\begin{equation}\label{eq: EnKF predict 1}
\wh{\bv}^{(k)}_{l+1}=\bPsi(\bv^{(k)}_l).
\end{equation}
Once the ensemble has evolved, the forecast distribution is then approximately characterized by the sample prior mean and covariance matrix
\begin{subequations}
\begin{equation}
\ds \wh{\m}_{l+1}=\frac{1}{K}\sum_{k=1}^{K}\,\wh{\bv}^{(k)}_{l+1}, \label{eq: EnKF prior 1}
\end{equation}
\begin{equation}
\ds \wh{C}_{l+1} =\frac{1}{K-1} \sum_{k=1}^{K}\,(\wh{\bv}^{(k)}_{l+1}-\wh{\m}_{l+1})(\wh{\bv}^{(k)}_{l+1}-\wh{\m}_{l+1})^{T},  \label{eq: EnKF prior 2}   
\end{equation}
\end{subequations}
which serve as low-dimensional representations of the prior mean and covariance, respectively. The forecast ensembles $\{ \wh{\bv}^{(k)}_{l+1} \}_{k=1}^{K}$ are then updated to obtain posterior ensembles in the analysis step using new observation data $\by_{l+1}$ given by \eqref{model: observation}. This is accomplished by first obtaining the posterior mean and covariance through the standard Kalman update formula
\begin{subequations}
\begin{equation}
\m_{l+1}=(I-K_{l+1}H)\wh{\m}_{l+1}+ K_{l+1}\by_{l+1} \label{eq:mean_update}
\end{equation}
\begin{equation}
C_{l+1}=(I-K_{l+1}H) \wh{C}_{l+1}. \label{eq:cov_update}
\end{equation}
\end{subequations}
Here the Kalman gain $K_{l+1}$ and the innovation covariance $S_{l+1}$ are respectively given by
\begin{subequations}
\begin{equation}
K_{l+1}=\wh{C}_{l+1}H^{T}S_{l+1}^{-1} \label{eq:gain_update}
\end{equation}
\begin{equation}
S_{l+1} =H \wh{C}_{l+1}H^{T} + \Gamma, \label{eq:S_update}
\end{equation}
\end{subequations}
where $\Gamma$ is given in \eqref{assump: noise Gaussian}. The ETKF ensemble particles are then updated using a suitably chosen transformation operator. In this regard, following the implementation of ETKF in \cite{Law15}, we define the prior centered ensemble as
\begin{equation}
    \label{eq:priorcenteredensemble}
\wh{X}_{l+1} = \frac{1}{\sqrt{K-1}} \left[ \wh{\bv}^{(1)}_{l+1}-\wh{\m}_{l+1}, \cdots, \wh{\bv}^{(K)}_{l+1}-\wh{\m}_{l+1} \right],\end{equation}
yielding the transformation operator $T_{l+1}$ 
\begin{equation}\label{eq: ETKF T}
T_{l+1}=\left[ I+(H\wh{X}_{l+1})^T \,\Gamma^{-1}(H\wh{X}_{l+1}) \right]^{-1}.
\end{equation}

Finally, the posterior ensemble is obtained as 
\begin{equation}
\ds \bv^{(k)}_{l+1} = \m_{l+1}+\bzeta^{(k)}_{l+1},\quad k = 1,\dots,K,
\label{eq:zeta}
\end{equation}
where each $\bzeta^{(k)}_{l+1}$ is defined as the scaled $k$-th column of 
\begin{equation}\label{eq: ETKF zeta}
\wh{X}_{l+1} T^{\frac{1}{2}}_{l+1} = \frac{1}{\sqrt{K-1}} \left[ \bzeta^{(1)}_{l+1}, \cdots, \bzeta^{(K)}_{l+1} \right].
\end{equation}
Observe that \eqref{eq: ETKF zeta} is defined so that each $\bzeta_{l+1}^k \sim  \cN(\0,C_{l+1})$, $k = 1, \dots, K$. 

\begin{remark}\label{rem: min}
The posterior mean in the analysis step \eqref{eq:mean_update} can equivalently be obtained by solving
\begin{equation}
\m_{l+1} = \underset{\bv}{\arg\min} 
\left(\frac{1}{2}\vert \by_{l+1}-H \bv \vert_{\Gamma}^2 + \frac{1}{2} \vert \bv-\wh{\m}_{l+1}^{(k)} \vert_{W_{l+1}}^2\right),
\label{eq:ETRF min}
\end{equation}
where the weighting matrix is given by $W_{l+1}=\wh{C}_{l+1}$, the sample prior covariance defined in \eqref{eq: EnKF prior 2}. The posterior covariance matrix is then obtained through the standard Kalman update formula derived under the assumption of Gaussian distributions. That is, by Bayes' rule, the posterior density of the state ${\bv}_{l+1}$ can be written as
\begin{equation}\label{eq: KF post}
\exp \left( -\frac{1}{2} \vert \bv_{l+1} - \m_{l+1} \vert^2_{C_{l+1}} \right) \propto \exp \left(  -\frac{1}{2}\vert \by_{l+1}-H \bv_{l+1} \vert_{\Gamma}^2 -\frac{1}{2} \vert \bv_{l+1}-\wh{\m}_{l+1} \vert_{\wh{C}_{l+1}}^2 \right),
\end{equation}
where the right-hand side is comprised of the likelihood and prior densities. Since the distribution is Gaussian, the posterior mean coincides with the maximum a posteriori (MAP) estimate, which is equivalent to minimizing the negative log-posterior density, leading directly to \eqref{eq:ETRF min}.
\end{remark}

\subsection{Inflation and localization}
\label{subsec:inflate}

The ETKF approach approximates the prior covariance from the ensembles using sample prior covariance, $\wh{C}_{l+1}$ in \eqref{eq: EnKF prior 2}. This approximation tends to be poor, especially when the ensemble size $K$ is much smaller than the state dimension $n$, which is common for purposes of dimension reduction and computational feasibility \cite{Matthias16}. Indeed, several adverse effects have been noticed when working with finite ensembles in high-dimensional systems \cite{Furrer07}. Inflation and localization are two techniques developed to offset these effects in practical implementations of ensemble Kalman filters including ETKF.  We briefly review them below.

\subsubsection*{Variance inflation \cite{Anderson07, Law15}} Over time, standard ensemble Kalman filtering often reduces the ensemble spread, causing an ``overly confident'' prior distribution to become increasingly less responsive to new observations. Specifically, model biases and/or sampling error from small ensembles may underestimate the sample covariance resulting  in poor assimilation of new data \cite{Furrer07}. Inflation corrects this by artificially increasing the prior ensemble spread, effectively adding uncertainty to the prior estimate so that the analysis remains sensitive to incoming observations, and can be implemented either by adding noise to each ensemble member or by scaling deviations around the ensemble mean. Here we use multiplicative inflation, and scale the prior centered ensemble $\wh{X}_{l+1}$ defined in \eqref{eq:priorcenteredensemble} as
\begin{equation}
\wh{X}_{l+1} \leftarrow \alpha \wh{X}_{l+1}, \quad \alpha > 1,\label{eq:inflation}
\end{equation}
which effectively inflates the prior covariance $\wh{C}_{l+1}$ by a factor of $\alpha^2$.

\subsubsection*{Localization \cite{Evensen06}} In contrast to inflation, localization is used to mollify the spurious correlations that occur when small ensemble sizes overestimate correlations between spatially distant or uncorrelated state components \cite{Hamill01}. These long-range spurious correlations not only degrade assimilation quality, but may also introduce numerical instability, since the covariance matrix is nearly rank-deficient when the ensemble size is much smaller than the state dimension. Localization effectively enforces the physically reasonable assumption that correlations between state components decay with increasing spatial separation, and can be implemented by elementwise multiplication (Hadamard product) as
\begin{equation}
\wh{C}_{l+1} \leftarrow \wh{C}_{l+1} \odot \cT.\label{eq:localization}
\end{equation}
Here $\wh{C}_{l+1}$ is the  sample prior covariance in \eqref{eq: EnKF prior 2} and $\cT$ is defined as a sparse, symmetric positive definite (SPD) localization matrix encoding the desired correlation decay. Observe that localization acts as a pre-processor that effectively smooths long-range spurious correlations while preserving meaningful local relationships in the covariance structure. Common choices of $\cT$ include the compactly supported Gaspari–Cohn kernel \cite{Gaspari99}, which ensures smooth decay beyond a cutoff radius, and structured Toeplitz matrices with optional tapering (e.g., assigning $1$ on the diagonal, $0.5$ on first off-diagonals, and so on). Simpler alternatives include binary Toeplitz matrices with fixed bandwidth. In our numerical experiments, we adopt a simplified, bandwidth-based localization approach to balance computational efficiency with robustness, which will be described in more detail in Section \ref{subsub: sparse inflation}.

Finally, we note that inflation and localization are typically used together in ensemble data assimilation -- inflation restores the ensemble spread, while localization stabilizes updates by suppressing unrealistic correlations. In this way  data assimilation can be used effectively even in high-dimensional systems with limited ensemble sizes.

\section{Structurally informed prior design}
\label{sec:method}
Although inflation and localization have been effective in offsetting problems that arise from using a low dimensional approximation of the prior covariance, ensemble-based methods  still fail to yield accurate results when the underlying solution exhibits discontinuous profiles. As demonstrated in \cite{TLda} for the the one-dimensional setting, the covariance matrix does not capture any information regarding discontinuities, which are important features that by all accounts should be incorporated into the prior belief. 

This paper expands on the one-dimensional mathematical framework for data assimilation that was developed in \cite{TLda} to provide statistically more accurate approximations for the state variables with discontinuous profiles while maintaining the robustness and computational efficiency of the methods. The expansion to 2D systems is non-trivial since the discontinuities exhibit more complex geometries and require careful treatment of both directional gradients and spatial correlations. 

For simplicity and ease of notation, we focus on a single time instance and drop the time subscript $l$ in \eqref{eq:ETRF min}, noting that the procedure is applied independently at each time step. As discussed in Remark \ref{rem: min}, the posterior mean is estimated by solving 
\begin{equation}
\label{eq:obj_w}
\m = \underset{\bv}{\text{argmin}}\left(\frac{1}{2}\vert \by-H \bv \vert_{\Gamma}^2 + \frac{1}{2} \vert \bv-\wh{\m} \vert_{W}^2\right). 
\end{equation}
Since the formulation employs $l_2$ regularization, the minimization problem has a closed-form solution, which can be computed efficiently. Regarding the weighting matrix $W$, we summarize that while \eqref{eq: KF post} suggests that the weighting matrix $W$ is equivalent to the original prior covariance matrix $\wh{C}$ in \eqref{eq: EnKF prior 2}, Section \ref{subsec:inflate} demonstrates that inflation and localization are needed to enhance $W$ to avoid potential issues that may arise from the poor approximation of $\wh{C}$ or small ensemble size.  Finally, a  {\em structurally informed} weighting matrix $W$ that uses prior gradient information was designed in \cite{TLda} for 1D systems to replace $\wh{C}$. However, as also already noted, formulating $W$ for 2D systems brings additional challenges, which is the main focus of this investigation.

To this end, in what follows we first review the 2D construction of the (standard) weighting matrix, which uses the prior covariance matrix $\wh{C}$ along with appropriate inflation and localization.  We then describe how to effectively design a new 2D weighting matrix to incorporate structural information.  For clarity purposes, we write $W_C$ when the weighting matrix is constructed from the usual covariance matrix $\wh{C}$, and $W_S$ when it is determined from structural information. Finally we note that in each 2D formulation of $W$ we consider a rectangular domain $\Omega=[L_{x1}, L_{x2}]\times[L_{y1}, L_{y2}]$, discretized uniformly with mesh sizes $\Delta x$ and $\Delta y$ in the $x$ and $y$ directions, respectively. The discrete field is correspondingly represented by $\bv \in \bbR^{n_x \times n_y}$,  where $n_x = (L_{x2}-L_{x1})/\Delta x +1$ and $n_y = (L_{y2}-L_{y1})/ \Delta y + 1$. At each spatial grid point we denote the solution as $\bv_{i,j} = \bv(x_i,y_j)$, where $x_i = L_{x1} + (i-1)\Delta x$ and $y_j = L_{y1} + (j-1)\Delta x$, for $i=1, \dots, n_x$, $j=1, \dots, n_y$.

\subsection{Covariance-based weighting matrix}
\label{sub:cov}

The general formula for calculating the prior sample covariance matrix $\wh{C}$ is given in \eqref{eq: EnKF prior 2}. In higher dimensional problems, each ensemble member and the ensemble mean must first be vectorized, as they are originally defined on a physical grid. To better understand the internal structure of $\wh{C}$ and to motivate the construction of a structurally informed weighting matrix, we now explicitly derive how $\wh{C}$ is assembled from the spatial vector fields.

For 2D vector fields $\wh{\m} \in \bbR^{n_x \times n_y}$ defined in \eqref{eq: EnKF prior 1} and $\wh{\bv}^{(k)}  \in \bbR^{n_x \times n_y}$ defined in \eqref{eq: EnKF predict 1}, $k = 1,\dots, K$, the corresponding  ``flattened'' vectors, $\wh{\m}_{\text{flat}} \in \bbR^{n_x n_y}$ and $\wh{\bv}^{(k)}_{\text{flat}} \in \bbR^{n_x n_y}$, are obtained by stacking their columns in column-major order. That is, they are connected by the linear index 
\begin{equation}\label{eq: index mapping}
m=i+(j-1)n_x, \quad i = 1,\dots,n_x, \ j = 1,\dots,n_y,
\end{equation}
so that 
\begin{equation} \label{eq: flat}
\left(\wh{\m}_{\text{flat}}\right)_m = \wh{\m}_{i,j}, \quad \text{and} \quad  \left(\wh{\bv}^{(k)}_{\text{flat}}\right)_{m} = \wh{\bv}^{(k)}_{i,j}, \quad m = 1,\dots, n_xn_y.
\end{equation}
Under this mapping, the diagonal entry $\wh{C}_{mm}$ of the matrix $\wh{C} \in \bbR^{n_x n_y \times n_x n_y}$ corresponds to the variance $\wh{V}_{i,j}$ at the spatial location $(x_i,y_j)$, given by
\begin{equation}
\wh{V}_{i,j}=\frac{1}{K-1} \sum_{k=1}^{K} (\wh{\bv}_{i,j}^{(k)}-\wh{\m}_{i,j})^2, \quad i=1, \dots, n_x, \ j=1, \dots, n_y. 
\label{eq: variance}
\end{equation}

The full covariance matrix $\wh{C}$ in \eqref{eq: EnKF prior 2} can also be formed using its diagonal entries (i.e., the pointwise variances) $\wh{C}_{mm} = \wh{V}_{i,j}$ in \eqref{eq: variance} and the Pearson correlation matrix:
\begin{equation}
\label{eq: full covariance}
\wh{C}_{mm'}=\sqrt{\wh{C}_{mm}}\,\wh{R}_{mm'}\sqrt{\wh{C}_{m'm'}}, \quad m,m' = 1,\dots,n_x n_y,
\end{equation}
where $\wh{R}_{mm'}$ is the component of the sample Pearson correlation matrix $\wh{R}$, given by 
\begin{equation}\label{eq: correlation coefficient}
\wh{R}_{mm'} = \frac{\ds \frac{1}{K-1} \sum_{k=1}^{K} 
\left( \left(\wh{\bv}_{\text{flat}}\right)_{m}^{(k)}-\left(\wh{\m}_{\text{flat}}\right)_{m}\right)
\left( \left(\wh{\bv}_{\text{flat}}\right)_{m'}^{(k)}-\left(\wh{\m}_{\text{flat}}\right)_{m'}\right)}
{\ds 
\sqrt{\frac{1}{K-1} \sum_{k=1}^{K} \left( \left(\wh{\bv}_{\text{flat}}\right)_{m}^{(k)}-\left(\wh{\m}_{\text{flat}}\right)_{m}\right)^2} \sqrt{\frac{1}{K-1} \sum_{k=1}^{K} \left( \left(\wh{\bv}_{\text{flat}}\right)_{m'}^{(k)}-\left(\wh{\m}_{\text{flat}}\right)_{m'}\right)^2}}.
\end{equation}

The full matrix $\wh{C}$ in \eqref{eq: full covariance} captures important spatial correlation structures, allowing information to be propagated across neighboring regions. However, as discussed in Section \ref{subsec:inflate}, due to the limited ensemble size and possibly sparse observations, it is easily affected by sampling noise and long-range spurious correlations. 

One way to address these issues is to simultaneously employ two complementary strategies: (1) variance inflation, in which  a multiplicative inflation factor $\alpha^2>1$ is applied, as defined in \eqref{eq:inflation}, and (2) localization, which is obtained via a sparse symmetric positive definite localization operator $\cT \in \bbR^{n_x n_y \times n_x n_y}$, as defined in \eqref{eq:localization}. These combined strategies yield the {\em covariance-based weighting matrix}
\begin{equation}\label{eq: Wc}
W_C = \alpha^2 \wh{C} \odot \cT, \quad \alpha > 1,
\end{equation}
which can then be directly incorporated into \eqref{eq:obj_w}.
\begin{remark}\label{rem: full observe}
In an ideal setting with densely sampled direct observations, the observation operator $H$ in \eqref{eq:obj_w} is simply the identity matrix $I$. In this case we apply a simplified form of localization by adopting $\cT = I$ which retains only the diagonal entries of the prior covariance matrix, corresponding to an extreme or degenerate case of covariance localization.   Spatial correlations are entirely ignored and only pointwise variances are used, resulting in 
\begin{equation}\label{eq: Cd}
W_C^{D}= \alpha^2 \wh{C}^D = \alpha^2
\begin{pmatrix}
\wh{C}_{1,1} & 0 & \cdots &  0 \\
0 & \ddots & \ddots & \vdots\\
\vdots & \ddots & \ddots & 0 \\
0 & \cdots  &   0 & \wh{C}_{n_x n_y,n_x n_y} \\
\end{pmatrix}.
\end{equation}
\end{remark}

\subsection{Gradient-based weighting matrix}
\label{sec:gradbasedweighting}
In some sense, defining $W_C$ in \eqref{eq: Wc} through the prior covariance matrix $\wh{C}$ can be viewed as a way to incorporate spatial information into the posterior update, with $\wh{C}_{mm}=\wh{V}_{i,j}$, defined in \eqref{eq: variance}, being proportional to the amount of penalty in the posterior update at $(x_i,y_j)$. This weighting matrix is often insufficient for state variables with discontinuity profiles, however, and motivated the use of the second moment of the discrete spatial gradient as a viable alternative in \cite{TLda}, although the scope was limited to 1D.

We now aim to retain this intuition for the more challenging 2D environment. In this regard we note that the second moment statistic, originally motivated by its analogy to covariance,  plays an essential role in capturing structural information about the underlying state. It provides a spatially adaptive scaling effect that distinguishes smooth from non-sooth regions, thereby helping to balance the data  and model fidelity terms. As the geometry in higher dimensions becomes more complex, it is natural to consider a broader class of statistics. Moreover our approach bears similarities to $l_p$ regularization often used in inverse problems, for which  $p$ penalizes the solution's sensitivity to sharp features. 

It is also important in 2D to consider directional information, and gradient-based regularization with directional sensitivity has been extensively studied in the context of total variation (TV) regularization, particularly in imaging denoising problems, as first introduced in the seminal paper \cite{Rudin92}. TV regularization promotes solutions with sparse gradients \cite{Gonzalez17}, with its directional treatment primarily prescribed as either  (1) rotational invariant isotropic TV \cite{Chan2011} or (2) anisotropic TV \cite{Condat17}, which favors horizontal and vertical structures. 

Our design of weighting matrix $W$ is therefore motivated by the considering both the  local statistics of spatial gradients, along with the  flexible aggregation of the  directional contributions. This will be achieved by first defining a gradient-based statistical matrix $\wh{S}$ analogous to $\wh{C}$, which is then further refined via local structural variation information.

\subsubsection{The construction of $\wh{S}$}\label{subsec:gradientstatistics}

Analogous to the construction of $\wh{C}$ described in Section \ref{sub:cov}, we begin by defining the diagonal entries of $\wh{S}$ using gradient-based statistics. We then extend this to a full matrix via incorporating correlation information. This process proceeds through the following steps:

\begin{enumerate}
\item We first approximate the spatial gradients using second-order central differences in each coordinate direction. That is, for fixed $y = y_j$, and $x = x_i$ we define
\begin{equation}
\label{eq:gradensembles}
\partial_x\wh{\bv}^{(k)}_{i,j}=\frac{\wh{\bv}^{(k)}_{i+1,j}-\wh{\bv}^{(k)}_{i-1,j}}{2\Delta x}, \quad \partial_y\wh{\bv}^{(k)}_{i,j}=\frac{\wh{\bv}^{(k)}_{i,j+1}-\wh{\bv}^{(k)}_{i,j-1}}{2\Delta y},
\end{equation}
where $\Delta x$ and $\Delta y$ are respectively the uniform mesh size. For convenience we also impose periodic boundary, noting that neither periodicity, nor grid mesh uniformity,  nor the particular gradient approximation method chosen,  is intrinsic to our approach. 
\item We then define the directional gradient statistics as 
\begin{equation}
\label{eq: x gradient}
\wh{S}^{x}_{i,j}=\frac{1}{K} \sum_{k=1}^{K} \left| \partial_x \wh{v}^{(k)}_{i,j} \right|^{\vartheta}, \quad \wh{S}^{y}_{i,j}=\frac{1}{K} \sum_{k=1}^{K} \left| \partial_y \wh{v}^{(k)}_{i,j} \right|^{\vartheta},
\end{equation}
where we call $\vartheta>0$ the {\em moment} parameter. 
\item The aggregated gradient statistic is then formed by combining directional contributions as
\begin{equation}\label{eq: joint gradient}
\wh{S}^{D}_{i,j} = \left( \wh{S}^{x}_{i,j} \right)^{\varphi} + \left( \wh{S}^{y}_{i,j} \right)^{\varphi},
\end{equation}
where we have introduced the {\em aggregation} parameter $\varphi$. Similar to $\wh{V}_{i,j}$ in \eqref{eq: variance}, $\wh{S}^{D}_{i,j}$ will serve as the diagonal entries $\wh{S}_{mm}$ of the gradient-based statistical matrix using the linear index $m$ defined in \eqref{eq: index mapping}.
\item Finally, the full gradient-based statistical matrix $\wh{S}$ is given by 
\begin{equation}\label{eq: full gradient}
\wh{S}_{mm'}=\sqrt{\wh{S}_{mm}}\,\wh{R}_{mm'}\sqrt{\wh{S}_{m'm'}}, \quad m,m' = 1,\dots,n_x n_y,
\end{equation}
where the Pearson correlation matrix $\wh{R}_{mm'}$ is  defined in \eqref{eq: correlation coefficient}.
\end{enumerate}

\begin{remark}\label{rem:newparameters}
The construction of $\wh{S}$ introduces two parameters: (1) the moment parameter $\vartheta$, which controls sensitivity to  directional derivatives; and (2) the aggregation parameter $\varphi$, which governs how directional components are aggregated. The interaction  between $\vartheta$ and $\varphi$ affects how forcefully the weighting differentiates between smooth and nonsmooth regions. In this way we can explore a range of aggregation schemes, including those inspired by isotropic and anisotropic formulations.  For example, $(\vartheta, \varphi)$ = (2,1) corresponds to isotropic aggregation, as $\vartheta = 2$ emphasizes squared gradients analogous to energy or $\ell^2$ norms, and $\varphi = 1$ aggregates the directional components symmetrically via linear summation. In contrast, $(\vartheta, \varphi) = (1,2)$ yields anisotropic aggregation, since $\vartheta = 1$ reflects the $\ell^1$ norm, and the larger aggregation parameter $\varphi = 2$ amplifies directional discrepancies, thus treating the $x-$ and $y-$directions more distinctly.
\end{remark}

\subsubsection{Scaling and localization}
\label{subsub: sparse inflation}

Since $\wh{S}$ is obtained using ensembles, the issues arising from the finite ensemble size remain. Analogous to inflation in Section \ref{subsec:inflate}, here we introduce a global scaling parameter $\beta>0$ to adjust the overall magnitude of the gradient-based weighting matrix $\wh{S}$. The parameter $\beta$ plays a similar role as the inflation parameter $\alpha^2$ and ensures   the compatibility of $\beta \wh{S}$ with the inflated covariance matrix $\alpha^2 \wh{C}$ obtained from \eqref{eq:inflation}. To further suppress spurious correlations, we again apply a localization operator $\cT$ as is done in \eqref{eq:localization}, leading to a weighting matrix 
\begin{equation}\label{eq: Ws0}
\wh{W}_S=\beta \wh{S} \odot \cT.
\end{equation}

\subsubsection{Directional structure-based correlation refinement}
\label{subsec:correlationrefinement}
Although localization effectively limits spurious long-range correlations by enforcing distance-based decay, it does not account for abrupt changes, such as discontinuities or steep gradients, in the state variable itself. Such structural features can invalidate assumed correlations even between nearby points. The method described below  further refines  $\wh{W}_S$ in \eqref{eq: Ws0} by incorporating local structural variation information. 

Consider the points $\bx_m = (x_i,y_j)$ and $\bx_{m'}=(x_{i'},y_{j'})$, $m,m' = 1,\dots,n_xn_y$, with pairwise Euclidean distance computed as
\[\| \bx_{m} - \bx_{m'} \| = \sqrt{(x_i-x_{i'})^2 + (y_j-y_{j'})^2}. \]
We then define a directional discrepancy statistic to quantify sharp transitions between each pair of points as
\begin{equation}
\label{eq: directional discrepancy}
d_{m,m'} = \frac{\left|\left(\wh{\m}_{\text{flat}}\right)_{m} - \left(\wh{\m}_{\text{flat}}\right)_{m'}\right|}{\| \bx_{m} - \bx_{m'} \|}, 
\end{equation}
where $\wh{\m}_{\text{flat}}$ is the (flattened) sample prior mean \eqref{eq: flat}.
Large directional discrepancies indicate potential discontinuities or sharp transitions, motivating correlation refinement. We incorporate this information into the correlation structure using thresholding. Specifically, we define a so-called structurally-aware  Pearson correlation matrix $\wt{R}$ with components 
\begin{equation}
\label{eq: mod_coeffs}
\wt{R}_{mm'} =
\begin{cases}
0, & d_{m,m'} > d_{\text{thresh}}, \\
\wh{R}_{mm'}, & d_{m,m'} \le d_{\text{thresh}}.
\end{cases}
\end{equation}
Here $\wh{R}$ is the sample Pearson correlation matrix defined in \eqref{eq: correlation coefficient} and $d_{\text{thresh}}$ is a user-defined threshold. By replacing  \eqref{eq: correlation coefficient} in \eqref{eq: full gradient} with  $\wt{R}$, we prescribe a refinement strategy that effectively removes correlations across sharp gradients while retaining correlations in smooth regions. 

This refinement approach is similar to localization (Section \ref{subsub: sparse inflation}), in the sense that both aim to reduce unwanted correlations informed by physical properties. From this point of view, the directional refinement can be interpreted as applying an additional binary mask matrix $\cM \in \bbR^{n_x n_y \times n_x n_y}$, where each entry is defined by
\begin{equation}\label{eq:Mask}
\cM_{m,m’} =
\begin{cases}
0, & d_{m,m’} > d_\text{thresh}, \\
1, & d_{m,m’} \le d_\text{thresh},
\end{cases}
\end{equation}
leading to a final weighting matrix of the form
\begin{equation} \label{eq: W in mask}
W_S = \wh{W}_S \odot \cM = \beta \, \wh{S}\odot \cT \odot \cM.
\end{equation}
The practical construction of $\cT$ and $\cM$ specific to the observation sparsity pattern used in our numerical experiments is discussed in Section \ref{sec:numerical}. Several remarks are in order.

\begin{remark}
\label{rem:pairwiseapproach}
Unlike the aggregated gradient statistics described in Section \ref{subsec:gradientstatistics}, which measures local smoothness using central differences at each grid point, the refinement prescribed by \eqref{eq: mod_coeffs} focuses on pairwise directional discrepancies, which are essential for capturing structural features that affect correlations. This pairwise approach is critical because it directly assesses whether a discontinuity or sharp gradient lies along the line connecting two points, providing information that cannot be captured by pointwise gradient magnitude alone. As a result, the refinement can adaptively suppress correlations inconsistent with local structural variations, even when traditional distance-based localization would otherwise retain them.
\end{remark}

\begin{remark}\label{rem:drawbackpairwise}
A potential limitation in using \eqref{eq: mod_coeffs} is that it only considers discrepancies between line segments using their endpoints, that is, it does not explicitly account for structural features along the path connecting them. In theory, two points separated by a discontinuity barrier could have similar state values, leading to small directional discrepancies and retained correlations that should be suppressed. However, our refinement is built on top of a localization operator that restricts correlations to local neighborhoods. When the localization radius is moderate, the risk of such configurations within a small local region is minimal, and the refinement remains effective in suppressing correlations inconsistent with local structural variations.
\end{remark}

\begin{remark}
\label{rem:Mimple}
Since $\cM$ is defined based on pairwise directional discrepancies, its diagonal entries can be directly set to 1. Furthermore, because it is built on top of $\cT$ via elementwise multiplication, it is unnecessary to compute or store entries of $\cM$ where $\cT$ is 0, which can be trivially set to 0. As a result, the overall cost of constructing $\cM$ is significantly reduced due to the sparsity of $\cT$. 
\end{remark}

\begin{remark}
\label{rem:thresholding}
The directional discrepancy in \eqref{eq: directional discrepancy} is inspired by classical edge detection principles, particularly the use of undivided differences to identify jump discontinuities in piecewise smooth data \cite{gelb2006adaptive}. In smooth regions, such differences scale with the local mesh size $\cO(\Delta x) $, while across discontinuities they grow proportionally to the jump magnitude, which is independent of $\Delta x$. Accordingly, $d_{\text{thresh}}$ is selected to lie between the typical noise scale in smooth regions (i.e., $ \cO(\Delta x) $) and the larger divided differences observed at discontinuities. This enables reliable discrimination of structural transitions while avoiding false suppression of correlations in smooth areas.
\end{remark}

\section{Numerical examples}
\label{sec:numerical}
Our numerical experiments are designed to illustrate and evaluate the proposed structurally informed data assimilation framework in settings with sharp spatial transitions.  We begin with a linear advection problem, using an initial condition that produces discontinuous solutions. This allows us to assess the method's basic behavior and sensitivity to design parameters. We then focus on the Burgers’ equation, which introduces nonlinear dynamics while preserving the challenge of discontinuous profiles. As our goal is to address the structural challenges in DA, these problems serve as meaningful prototypes. In particular we note that unlike classical numerical methods for conservation laws, DA has yet to establish a common framework for PDEs that admit solutions with discontinuous profiles.  Finally, we address sparse observation scenarios by detailing the implementation of localization and directional refinement, and demonstrating performance under limited data. In all cases the single time instance posterior mean is estimated via \eqref{eq:obj_w}, and we compare results for which  the weighting matrix $W$ is chosen either as $W_C$, reflecting the traditional covariance weighting approach (Section \ref{sub:cov}), or as $W_S$, to reflect our new gradient-based structurally-informed approach (Section \ref{sec:gradbasedweighting}).

\begin{remark}\label{rem: lit review}
Most existing developments of DA methods focus on turbulent dynamics or low-dimensional chaotic systems, such as the Lorenz models \cite{Calvello2025meanfield, LeProvost2024Preserving}. Fewer studies directly target hyperbolic conservation laws, which exhibit sharp features such as shocks and contact discontinuities. Early efforts include nudging-based techniques derived from kinetic formulations, as in \cite{boulanger2015data}, which demonstrated applicability to 1D Burgers and Saint-Venant systems. A comparative study of ensemble-based and variational DA for 1D ideal magnetohydrodynamics (MHD) was presented in \cite{Arnal2024MHD}, while the method introduced in \cite{Hansen2024NormalScoreEnKF} applies a normal-score transformation within the EnKF framework to better handle non-Gaussian distributed states in 1D shock wave propagation. More recent work has begun to address the inherent difficulties of assimilating discontinuous states. For example, the ensemble transform particle filter in \cite{Subrahmanya2025FeatureDA} is extended via feature-aligned interpolation to preserve shocks and contacts in the compressible Euler equations. Such studies highlight the need for DA strategies that respect the structure of solutions governed by hyperbolic PDEs. Our framework contributes to this emerging direction by explicitly encoding structural information derived from gradient-based statistics into the assimilation step, and by introducing directional structural-based correlation refinement mechanisms tailored to discontinuous states.
\end{remark}

\subsection{Evaluation metrics}
\label{subsec:metrics}
We use several complementary error metrics to capture both pointwise accuracy and overall pattern agreement with the analytical solution to quantitatively evaluate the performance of different data assimilation methods.  These include the following:
\begin{enumerate}
\item \textbf{Pointwise absolute error}: At each grid point and observation time we calculate
\begin{equation} 
err(x_i,y_j,t) = \vert u_{i,j}(t) - u_{\text{true}}(x_i,y_j,t) \vert,\quad i = 1,\dots,n_x, \ j=1, \dots, n_y.
\label{eq:error}
\end{equation}
where $u_{\text{true}}(x_i,y_j, t)$ is the true solution at time $t$ for each $(x_i,y_j)$, and $u_{i,j}(t)$ is the corresponding posterior mean solution. As the observation data is available every $\Delta t_{obs}$ unit, we only compute the error at the time when the posterior mean is obtained by assimilating data, that is, $t=t_{q} = q \Delta t_{obs}$, where $q = 1, \dots, L_{obs}$, and $L_{obs} = \frac{T}{\Delta t_{obs}}$.
\item \textbf{Relative $l_1$ and $l_2$ errors}:
To summarize errors over the domain,  at each assimilation time $t_q$ we compute
\begin{equation}
\label{eq:l1err}
\ds err_{l1}(t_q) = \frac{\ds \sum_{i = 1}^{n_x} \sum_{j=1}^{n_y} |u_{i,j}(t_q)- u_{\text{true}}(x_i,y_j,t_q)| / (n_x n_y)}{\ds \sum_{i = 1}^{n_x} \sum_{j=1}^{n_y} |u_{\text{true}}(x_i,y_j,t_q)|/ (n_x n_y)} 
\end{equation}
and 
\begin{equation}
\label{eq:l2err}
err_{\ell_2}(t_q) = \frac{\ds \sqrt{\sum_{i = 1}^{n_x} \sum_{j=1}^{n_y} (u_{i,j}(t_q)- u_{\text{true}}(x_i,y_j,t_q))^2/ (n_x n_y)}}{\ds\sqrt{\sum_{i = 1}^{n_x} \sum_{j=1}^{n_y} (u_{\text{true}}(x_i,y_j,t_q))^2/ (n_x n_y)}}.
\end{equation}

\item \textbf{Pattern correlation}:  We  employ  a statistical measure that quantifies the similarity between two spatial patterns, namely, the Pearson correlation coefficient computed between the true solution and the posterior mean solution and given by
\begin{equation}
\label{eq:pcorr}
Pcorr (t_q) =  \frac{\ds \sum_{i = 1}^{n_x} \sum_{j=1}^{n_y} \left(u_{i,j}(t_q)-\bar{u}(t_q)\right)\left(u_{\text{true}}(x_i,y_j,t_q)-\bar{u}_{\text{true}}(t_q)\right)}{\ds \sqrt{\sum_{i = 1}^{n_x} \sum_{j=1}^{n_y} \left(u_{i,j}(t_q)-\bar{u}(t_q)\right)^2} \sqrt{\sum_{i = 1}^{n_x} \sum_{j=1}^{n_y} \left(u_{\text{true}}(x_i,y_j,t_q)-\bar{u}_{\text{true}}(t_q)\right)^2}}.
\end{equation}
Here $\ds \bar{u}(t_q)=\frac{1}{n_x n_y} \sum_{i = 1}^{n_x} \sum_{j=1}^{n_y} u_{i,j}(t_q)$ and $\ds \bar{u}_{\text{true}}(t_q) = \frac{1}{n_x n_y} \sum_{i = 1}^{n_x} \sum_{j=1}^{n_y} u_{\text{true}}(x_i,y_j,t_q)$ are respectively the spatial means of the posterior and true fields.

\item \textbf{Long-term performance metrics}: We calculate the mean relative $l_1$ and $l_2$ errors
\begin{equation}
\label{eq:meanl1err}
\ds e_{l*} = \frac{1}{L_{obs}-q_0+1}\sum_{q=q_0}^{L_{obs}} err_{l*}(t_q)
\end{equation}
with $*=\{1,2\}$ and mean pattern correlation
\begin{equation}
Pc = \frac{1}{L_{obs}-q_0+1}\sum_{q=q_0}^{L_{obs}} Pcorr(t_q),
\end{equation}
where the calculation starts at $q_0$ to discard the burn in period. We choose $q_0 = L_{obs}/2$ in our implementations so that metrics are averaged over the second half of the assimilation window.
\end{enumerate}

\subsection{Linear advection equation with dense observations}
\label{subset:linear}
We first test our method on a linear dynamics system with full observations to both illustrate the idea and systematically investigate performance for a range of gradient weighting parameters.

\subsubsection*{Problem setup} We consider the 2D  linear advection equation within a square domain $[0,1]\times[0,1]$ for $t \in [0,2]$:
\begin{equation}
\begin{aligned}
& u_t + 0.5 u_x -  u_y = 0, \\
& u(t=0) = u_0(x,y) = \begin{cases}
 2y+0.4 & \text{ if } (x,y) \in [0.4, 0.6]\times [0.3,0.4], \\
 1.2 & \text{ if } (x,y) \in [0.4, 0.6]\times [0.4,0.6], \\
 -y+1.8 & \text{ if } (x,y) \in [0.4, 0.6]\times [0.6,0.8], \\
 1 & \text{ otherwise }.
\end{cases}
\end{aligned}
\label{eq:2Dlinearadvection}
\end{equation}
The discontinuous profile in both the $x$ and $y$ directions is chosen with intent, so that the impact of using structural information via weighting matrix $W_S^D$ in \eqref{eq: Sd} can be adequately assessed. For simplicity we impose periodic boundary conditions. The analytical solution is then given by
\begin{equation}
\label{eq:u_analytic}
    u_{\text{true}}(x,y,t) = u_0 \left( (x-0.5t) \text{ mod } 1, (y+t) \text{ mod } 1 \right).
\end{equation}

\subsubsection*{Numerical solver} The spatial domain is  uniformly discretized with $\Delta x = \Delta y = 10^{-2}$ so that $n_x=n_y=101$. Numerical experiments conducted in \cite{TLda} demonstrated that accurate data assimilation results critically depend on using an accurate numerical PDE solver for the forecasting step. We therefore employ the fifth order finite difference weighted essentially non-oscillatory (WENO) method \cite{Jiang96}, which yields high-order accuracy in smooth regions while still resolving discontinuities.  In this regard WENO is particularly effective when paired with our new weighting matrix, which puts more weight on the prior in smooth regions. To maintain high order accuracy,  we use the total variation diminishing third-order Runge-Kutta (TVDRK3) for time integration  \cite{Liu94}, where we fix time step $\Delta t = 5 \times 10^{-3}$ to ensure the CFL condition is satisfied with $\#CFL = 1$.

\subsection*{Observation and ensemble setup} We generate  observation data in \eqref{model: observation} with $\Gamma = \gamma^2 I$, $\gamma = 0.01$, using the analytical solution \eqref{eq:u_analytic}.  We assume that the data are {\em densely} observed at regular time intervals, that is, at every point in the spatial domain  with $\Delta t_{obs} = 5 \Delta t$ = $2.5 \times 10^{-2}$. Finally, we consider  $K=100$ ensembles, where each is initialized using $u_0(x,y)$ in \eqref{eq:2Dlinearadvection}, perturbed by additive i.i.d. noise with mean $0$ and standard deviation $0.1$.

\subsection*{Gradient statistics and structural features}  

As discussed in Remark \ref{rem: full observe}, a diagonal weighting matrix is appropriate when observations are dense, with the gradient-based weighting matrix  given by
\begin{equation}\label{eq: Sd}
W_S^{D}= \beta \wh{S}^D = \beta
\begin{pmatrix}
\wh{S}_{1,1} & 0 & \cdots &  0 \\
0 & \ddots & \ddots & \vdots\\
\vdots & \ddots & \ddots & 0 \\
0 & \cdots  &   0 & \wh{S}_{n_x n_y,n_x n_y} \\
\end{pmatrix}.
\end{equation}
Following the discussion surrounding  \eqref{eq: Ws0}, we  note that $\beta$ should generally be chosen so that $W_S^D$ has comparable magnitude to $W_C^D$ in \eqref{eq: Cd} used in the classical ETKF method. This can be accomplished by specifying 
\begin{equation}\label{eq:betat}
\beta:=\wt{\beta} / \Vert \wh{S}^D \Vert_{\max},  \end{equation} 
where $\Vert \cdot \Vert_{\max}$ is the elementwise infinity norm. This normalization ensures that the gradient-based weighting matrix $W_S^D$ operates on a similar scale as the  covariance-based weighting matrix $W_C^D$.

\begin{figure}[h]
\centering
\includegraphics[width=0.4\textwidth]{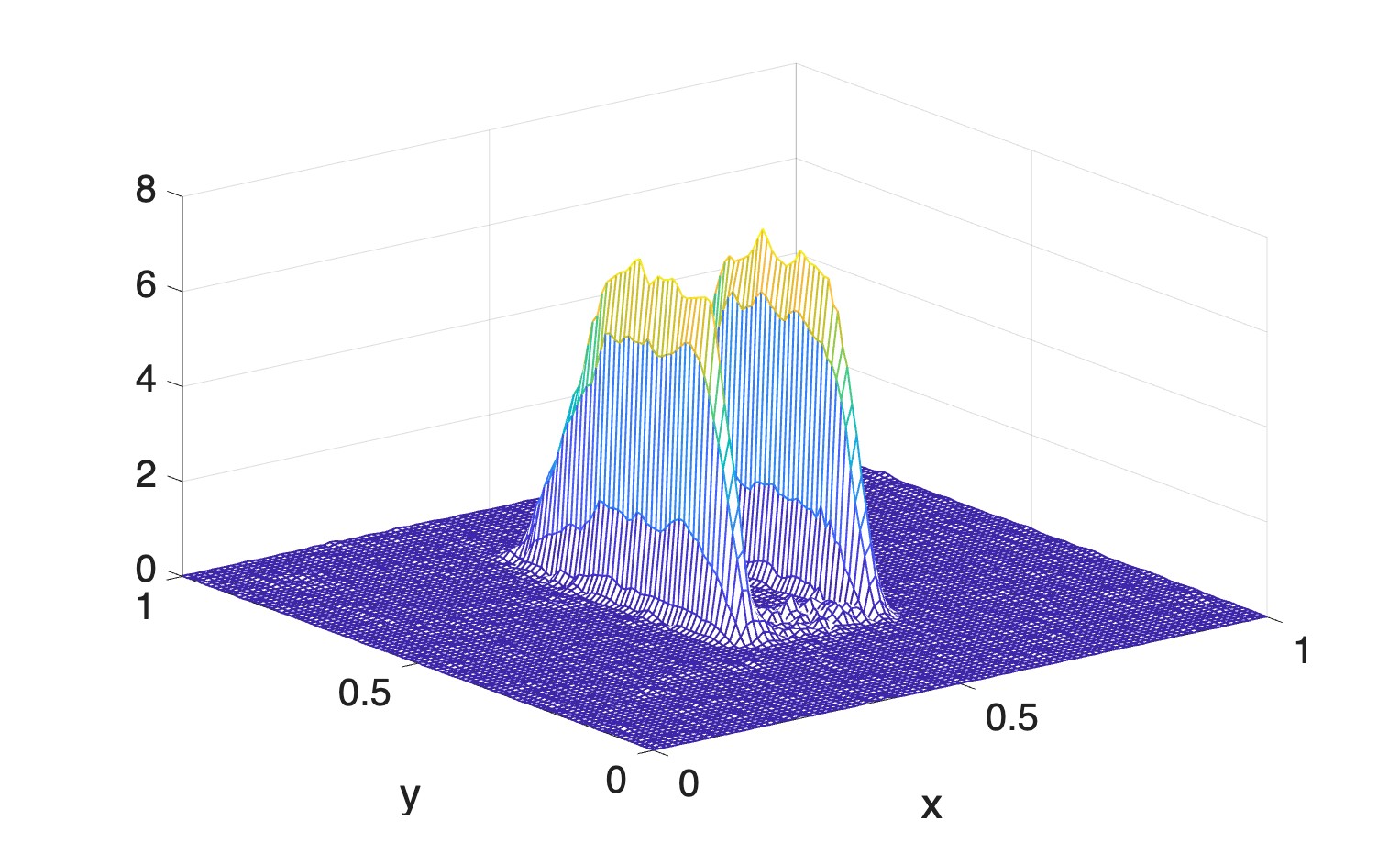}
\includegraphics[width=0.4\textwidth]{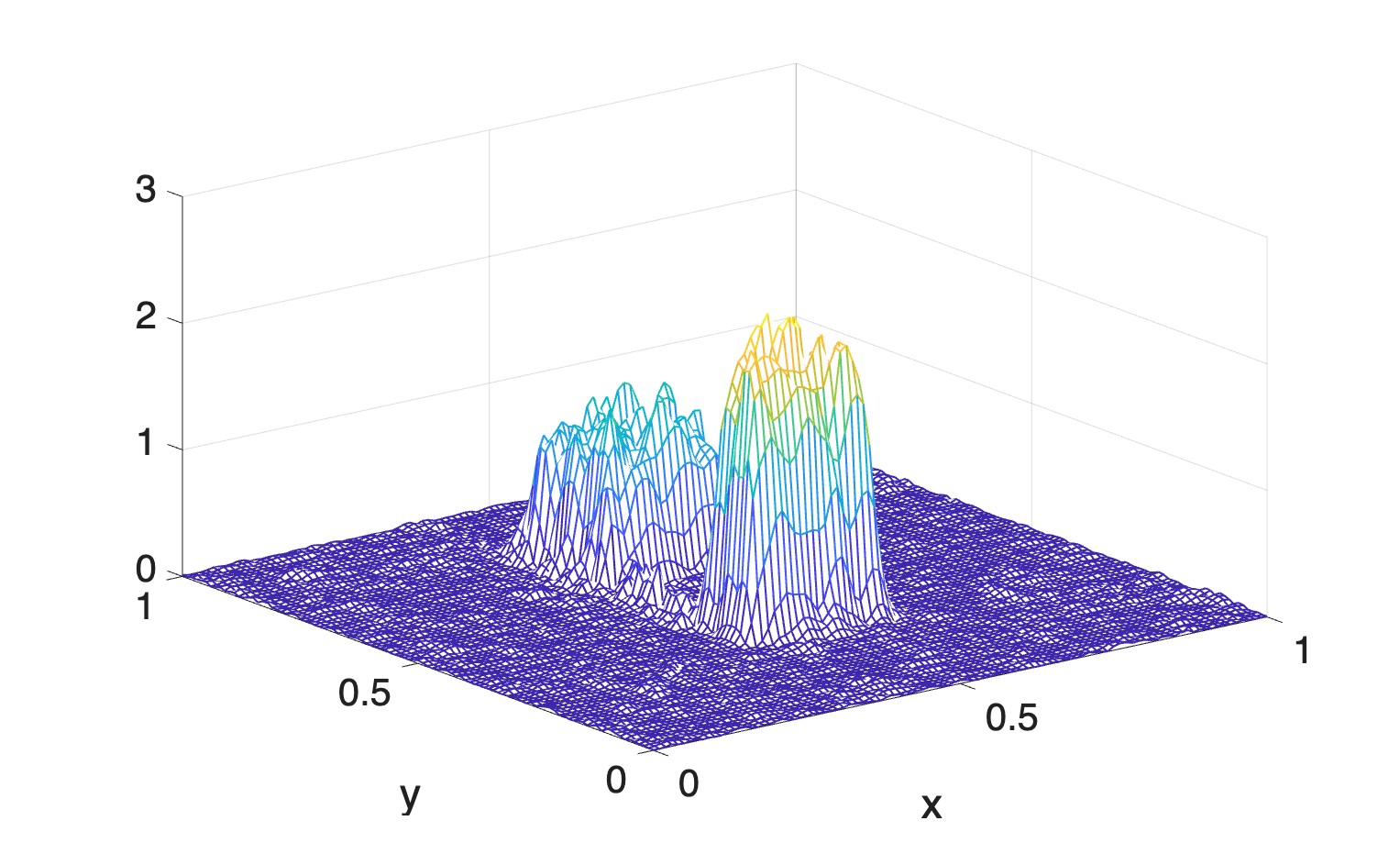}\\
\includegraphics[width=0.4\textwidth]{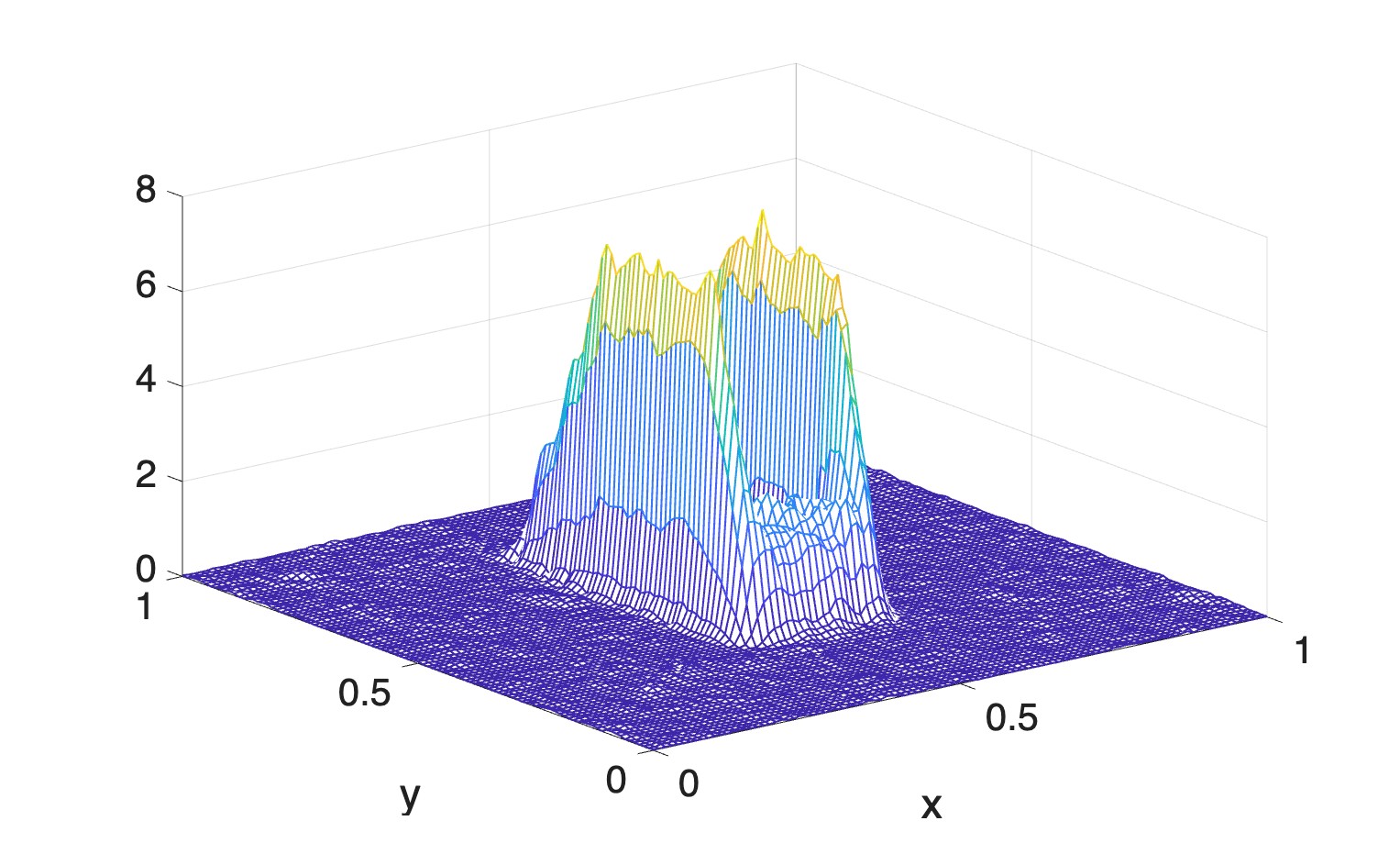}
\includegraphics[width=0.4\textwidth]{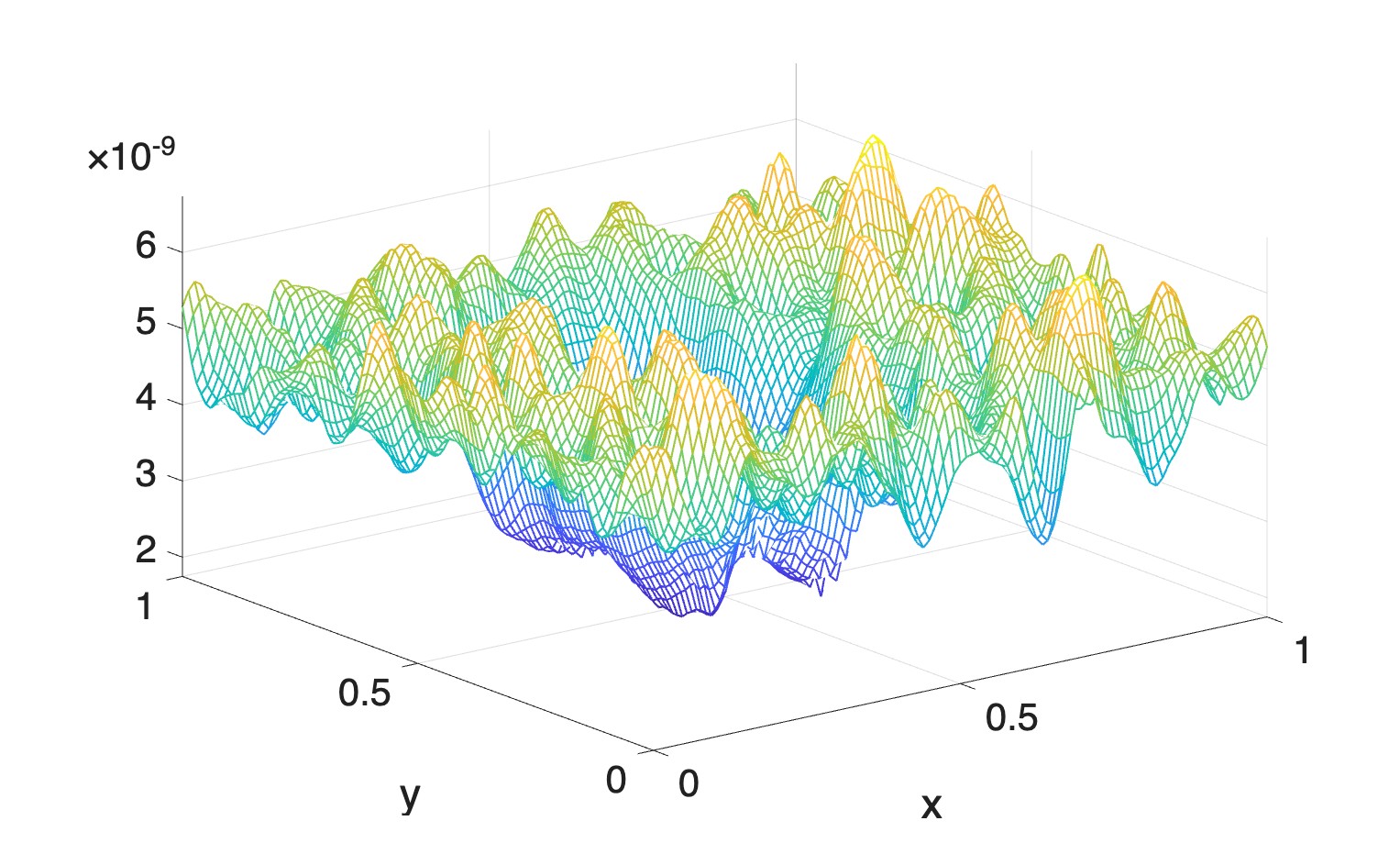}
\caption{Gradient statistics \eqref{eq: x gradient} for (top-left)  $x$  and (top-right) $y$; (bottom-left) aggregate gradient statistic \eqref{eq: joint gradient}; and (bottom-right) variance $\wh{V}_{i,j}$ \eqref{eq: variance} at time $t =2$. Here $(\vartheta,\varphi)=(1,1)$ and $\wt\beta = 10^{-3}$ in \eqref{eq:betat} is used to determine $\beta$ in \eqref{eq: Sd}.}
\label{fig: csm}
\end{figure}

The benefits of employing $W = {W}_S^D$ in \eqref{eq: Sd} can be observed in Figure \ref{fig: csm}, which displays the $x$ and $y$ gradient statistics, $\wh{S}_{i,j}^{x}$ and $\wh{S}_{i,j}^{y}$ in \eqref{eq: x gradient},  the joint gradient $\wh{S}^D_{i,j}$ in \eqref{eq: joint gradient}), along with the variance $\wh{V}_{i,j}$ \eqref{eq: variance} calculated at time $t = 2$.  The corresponding  parameters are $(\vartheta,\varphi)=(1,1)$ and $\wt{\beta}=10^{-3}$ as defined in \eqref{eq:betat}.  As parameter optimization is not the focus of our investigation, there was no further parameter tuning. It is evident that in addition to correctly identifying discontinuity locations, the gradient steepness is also captured.  In particular, the steepest gradient statistics occur along $x=0.4$, $y \in [0.3,0.8]$ and $x=0.6$, $y \in [0.3,0.8]$ respectively, and is reflected by the gradient statistic magnitudes. Figure \ref{fig: csm}(bottom-right) demonstrates that the variance, by contrast, is not able to delineate any local structural information.

\subsection*{Method performance comparison} Table \ref{tab:cov and gradient all} compares the evaluation metrics in Section \ref{subsec:metrics} for weighting matrices $W_C^D$ in \eqref{eq: Cd} and $W_S^D$ in \eqref{eq: Sd} with prespeficied pairs of $(\vartheta,\varphi)$. For each case, we test on some preselected parameter values, $\alpha$ for $W_C^D$ and $\wt{\beta}$ for $W_S^D$ respectively. Since for the traditional covariance weighting approach using $W_C^D$ performance metrics outcomes vary somewhat depending on parameter choice $\alpha$, we present the results for both $\alpha=4$ and $\alpha=6$. Observe that the former choice of $\alpha$ gives a smaller relative $l_1$ error $e_{l_1}$, while the latter perform better according to relative $l_2$ error $e_{l_2}$ and pattern correlation $Pc$. As the structurally data assimilation framework yields more consistent results across the metrics, we show the metrics as functions of $(\vartheta,\varphi)$ with the selected $\wt{\beta}$ providing best performance for all the metrics.

\begin{table}[h!]
\centering
\begin{tabular}{c c c c c c c c c }
\toprule
parameter & \multicolumn{2}{c}{$W_C^D$ \eqref{eq: Cd}} & \multicolumn{6}{c}{$W_S^D$ \eqref{eq: Sd}} \\ 
\midrule
$(\vartheta,\varphi)$ & & & (1/2,\ 1) & (1/2,\ 2) & (1,1) & (1,2) & (2,1) & (2,2) \\ 
\midrule
 $\alpha$/$\wt{\beta}$  & 4.0e-00  & 6.0e-00  & 1.0e-04  & 1.0e-03  & 1.0e-03 & 1.0e-01  & 1.0e-01 & 1.0e+02 \\
\midrule
$e_{l_1}$  \eqref{eq:l1err} & 2.2e-03 & 2.7e-03 & 1.6e-03 & \textbf{1.1e-03} & \textbf{1.1e-03} & 1.2e-03 & 1.2e-03 & 1.4e-03\\
$e_{l_2} \eqref{eq:l2err}$  & 5.9e-03 & 5.2e-03 & 4.7e-03 & \textbf{2.5e-03} & \textbf{2.5e-03} & 3.1e-03 & 3.1e-03 & 3.2e-03\\
$Pc$ \eqref{eq:pcorr} & 99.00\% & 99.37\% & 99.49\% &  \textbf{99.86\%} & \textbf{99.86\%} & 99.78\% & 99.78\% & 99.76\% \\
\bottomrule
\end{tabular}
\caption{The mean relative $l_1$, $l_2$ errors and the mean pattern correlation for weighting matrices $W_C^D$ in \eqref{eq: Cd}  and $W_S^D$ in \eqref{eq: Sd} for different parameters $\alpha$, $(\vartheta,\varphi)$ and $\wt{\beta}$.}
\label{tab:cov and gradient all}
\end{table}

In comparing the performance metrics in Table \ref{tab:cov and gradient all}, we see that employing the structurally-informed weighting matrix $W_S^D$ for all combinations of $(\vartheta,\varphi)$ and $\wt{\beta}$ consistently yields better outcomes than employing the standard covariance matrix $W_C^D$. Among all the cases, the structurally-informed weighting matrix $W_S^D$ with $(\vartheta,\varphi)$ being $(1/2,2)$ and $(1,1)$ provides the best results. Indeed, we observe  that the performance metrics are identical for $(\vartheta_1,\varphi_1)$ and  $(\vartheta_2,\varphi_2)$ whenever $\vartheta_1\varphi_1 = \vartheta_2\varphi_2$,   possibly due to the symmetric structure of \eqref{eq:u_analytic}, although a general understanding of the moment and aggregate parameter sensitivity with respect to their isotropic and anisotropic formulations requires further study.

\subsection*{Posterior solution cross-section}
\begin{figure}[h!]
\centering
\includegraphics[width=0.24\textwidth]{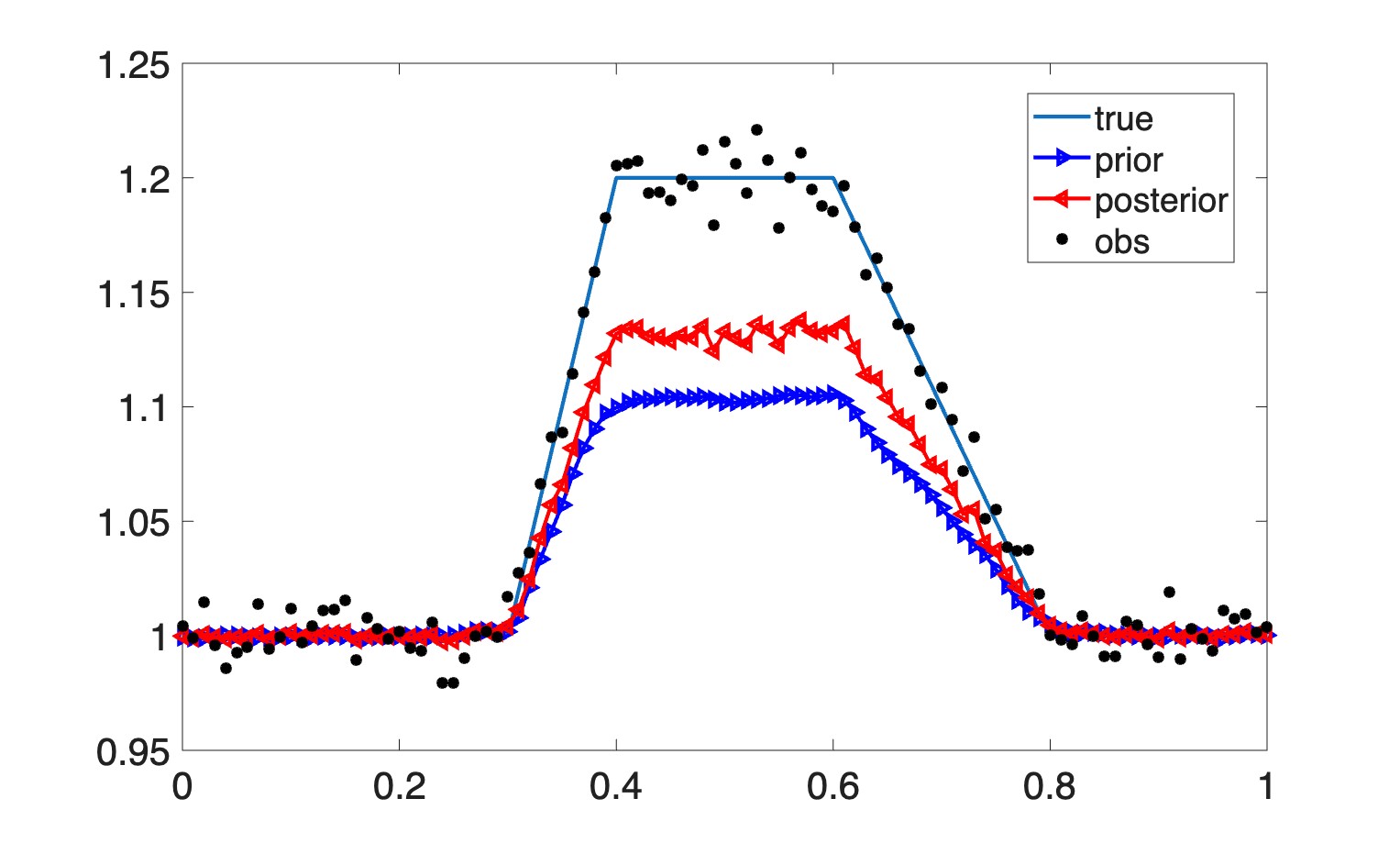}
\includegraphics[width=0.24\textwidth]{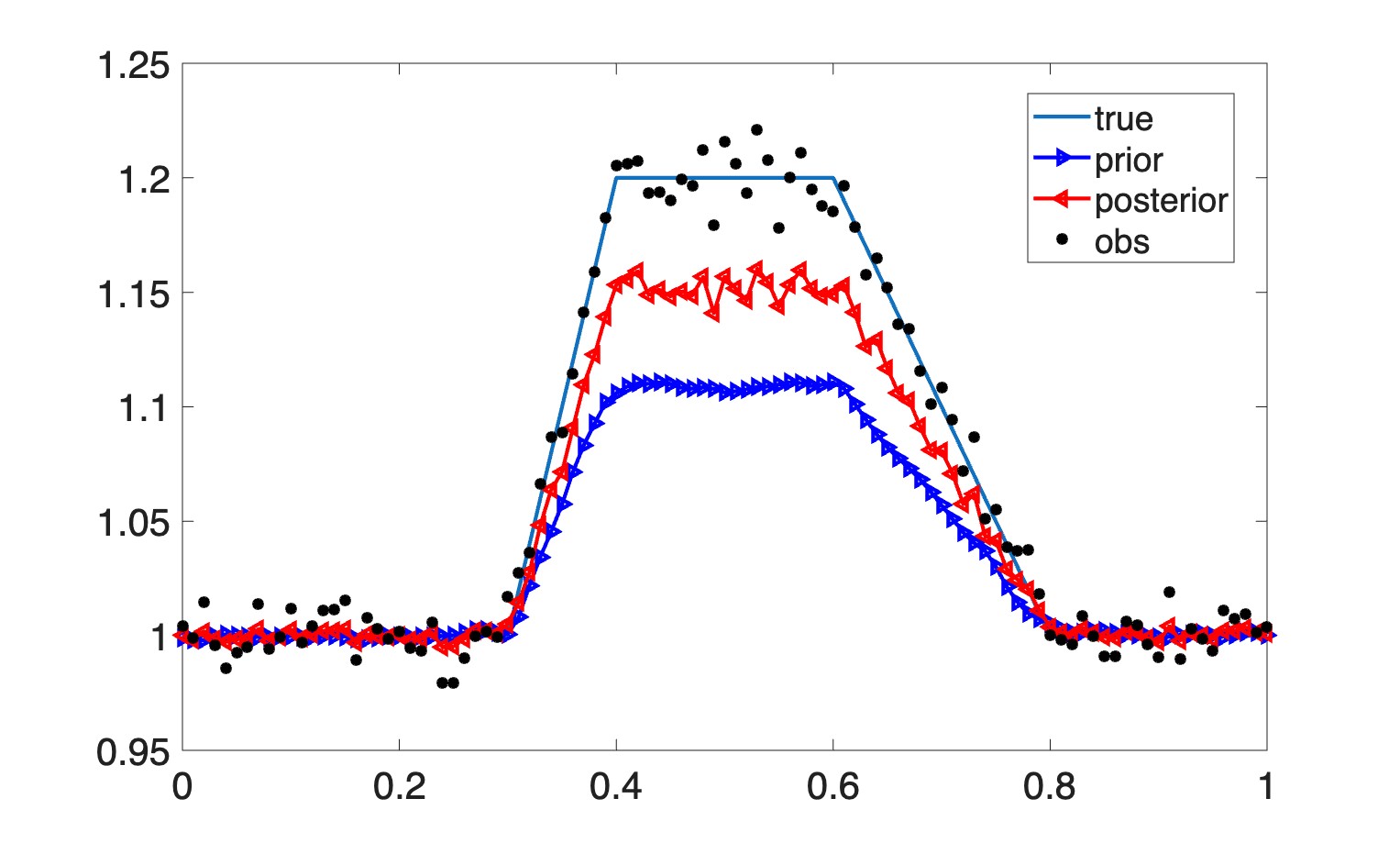}
\includegraphics[width=0.24\textwidth]{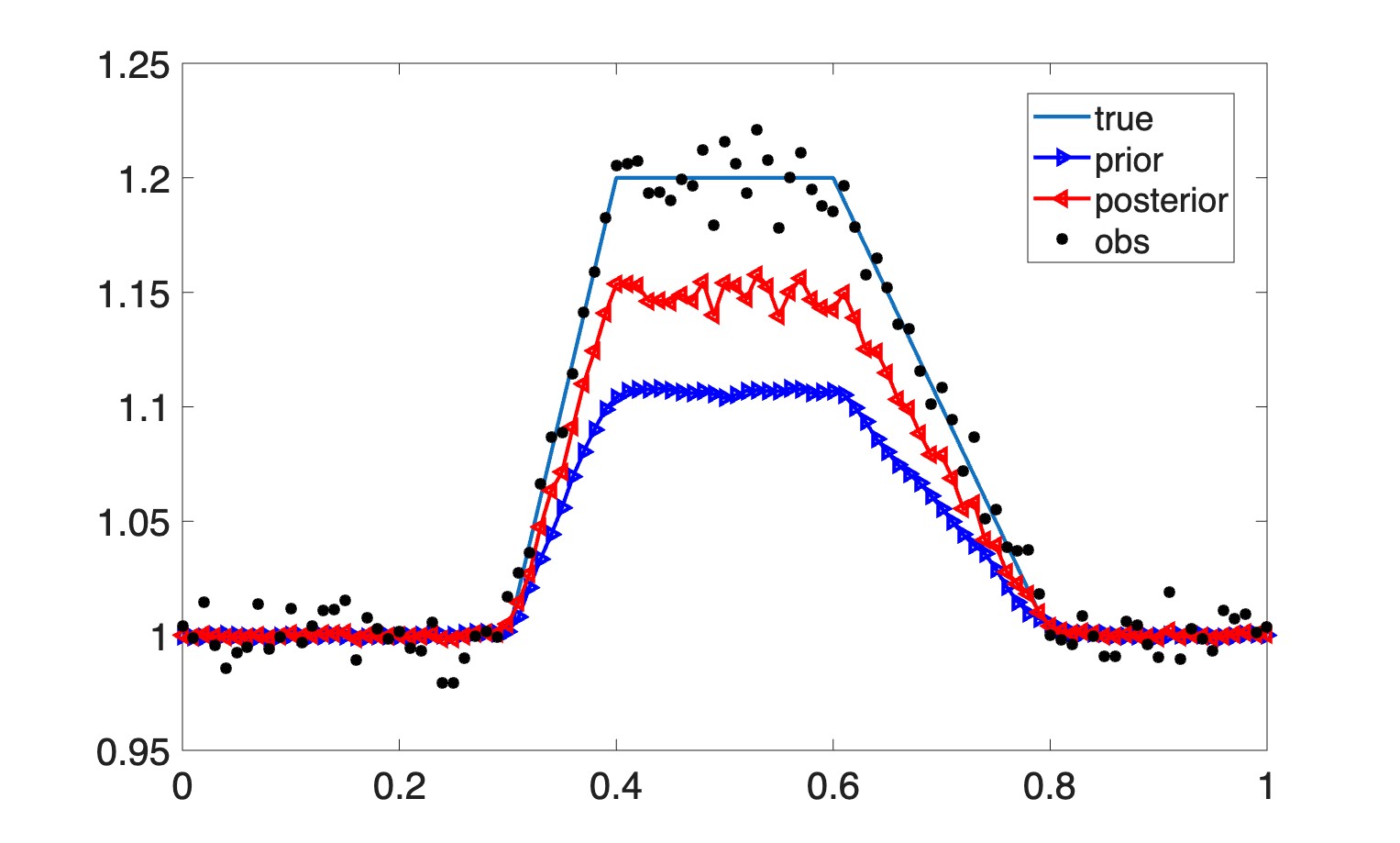}
\includegraphics[width=0.24\textwidth]{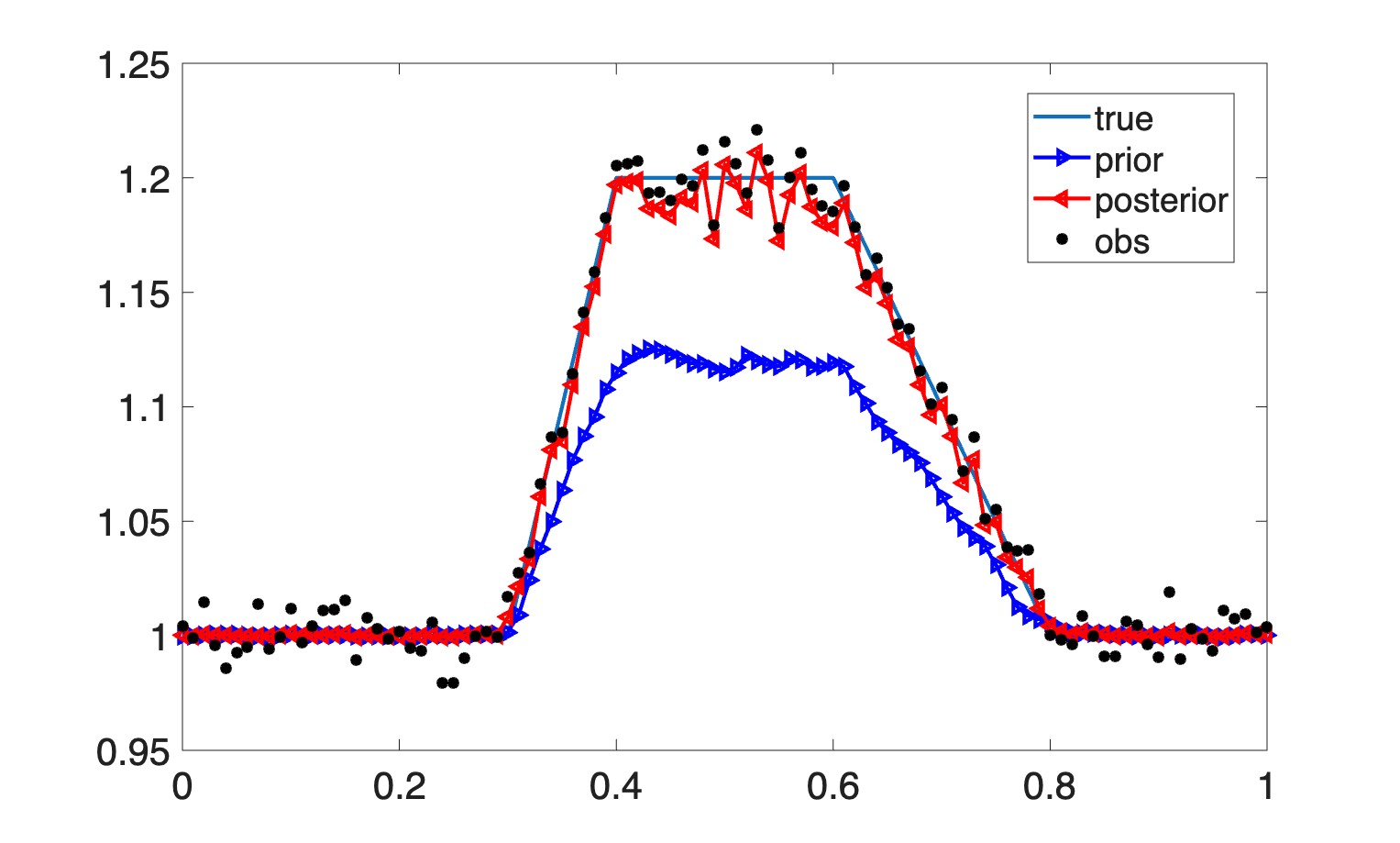} \\
\includegraphics[width=0.24\textwidth]{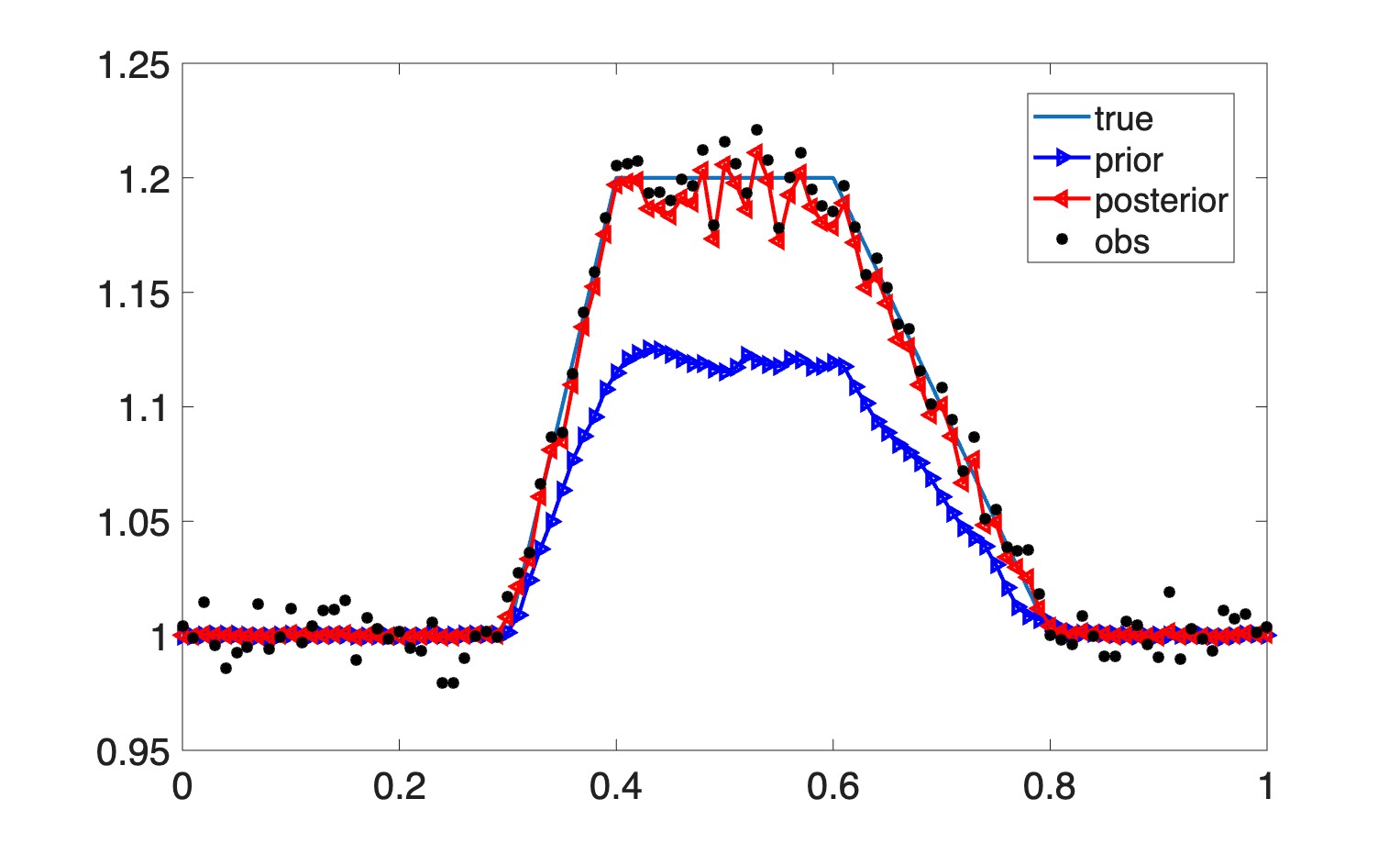}
\includegraphics[width=0.24\textwidth]{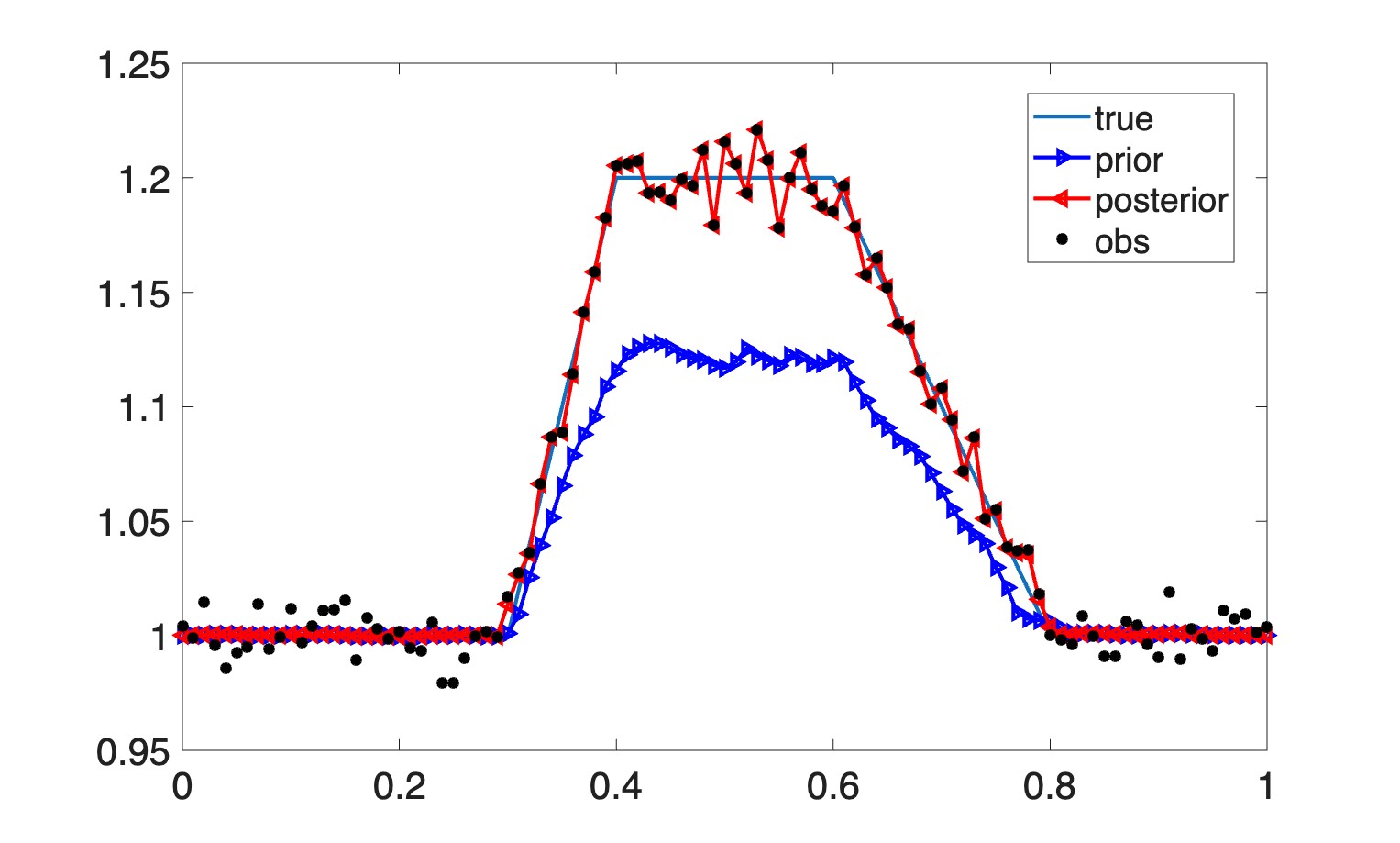}
\includegraphics[width=0.24\textwidth]{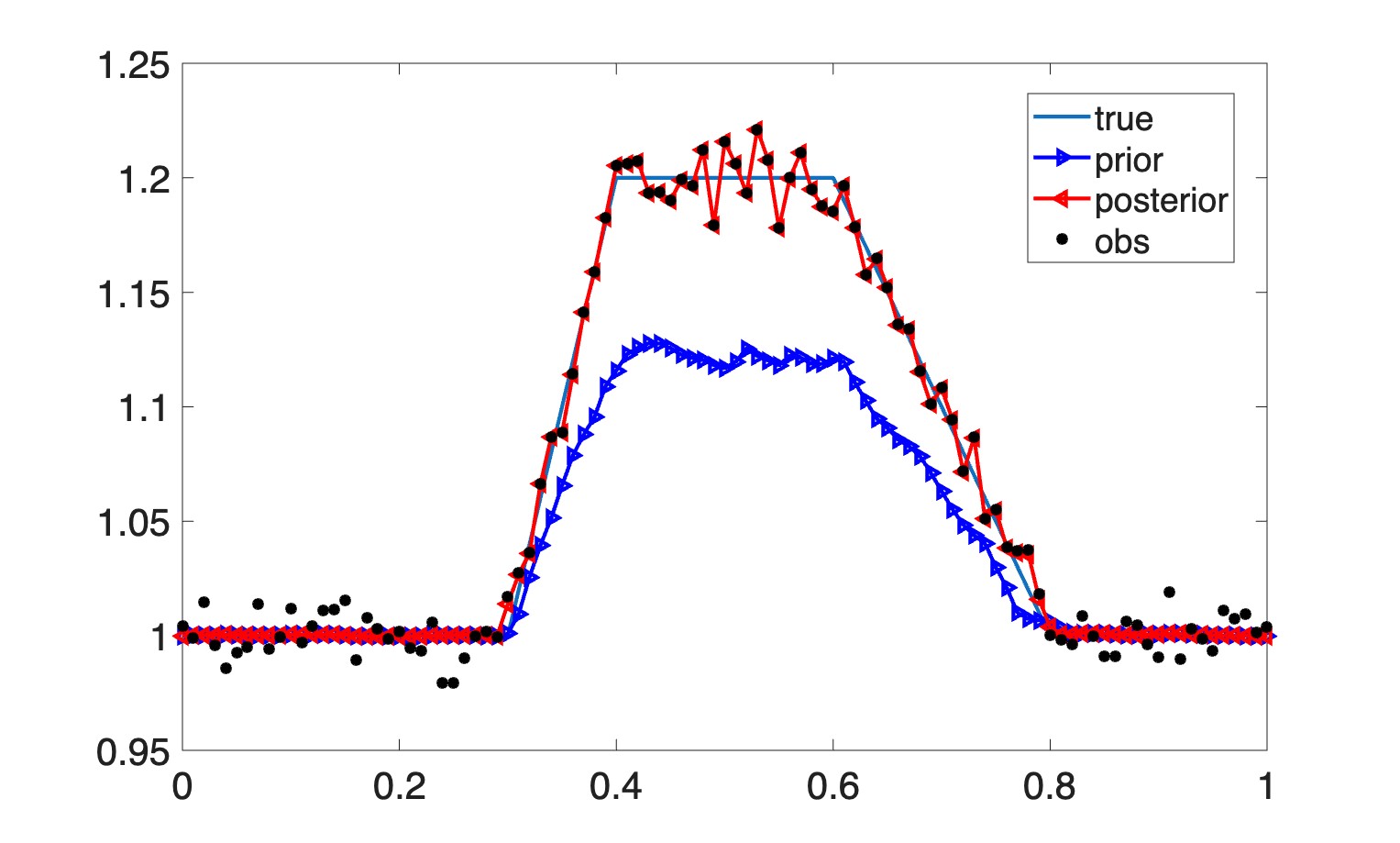}
\includegraphics[width=0.24\textwidth]{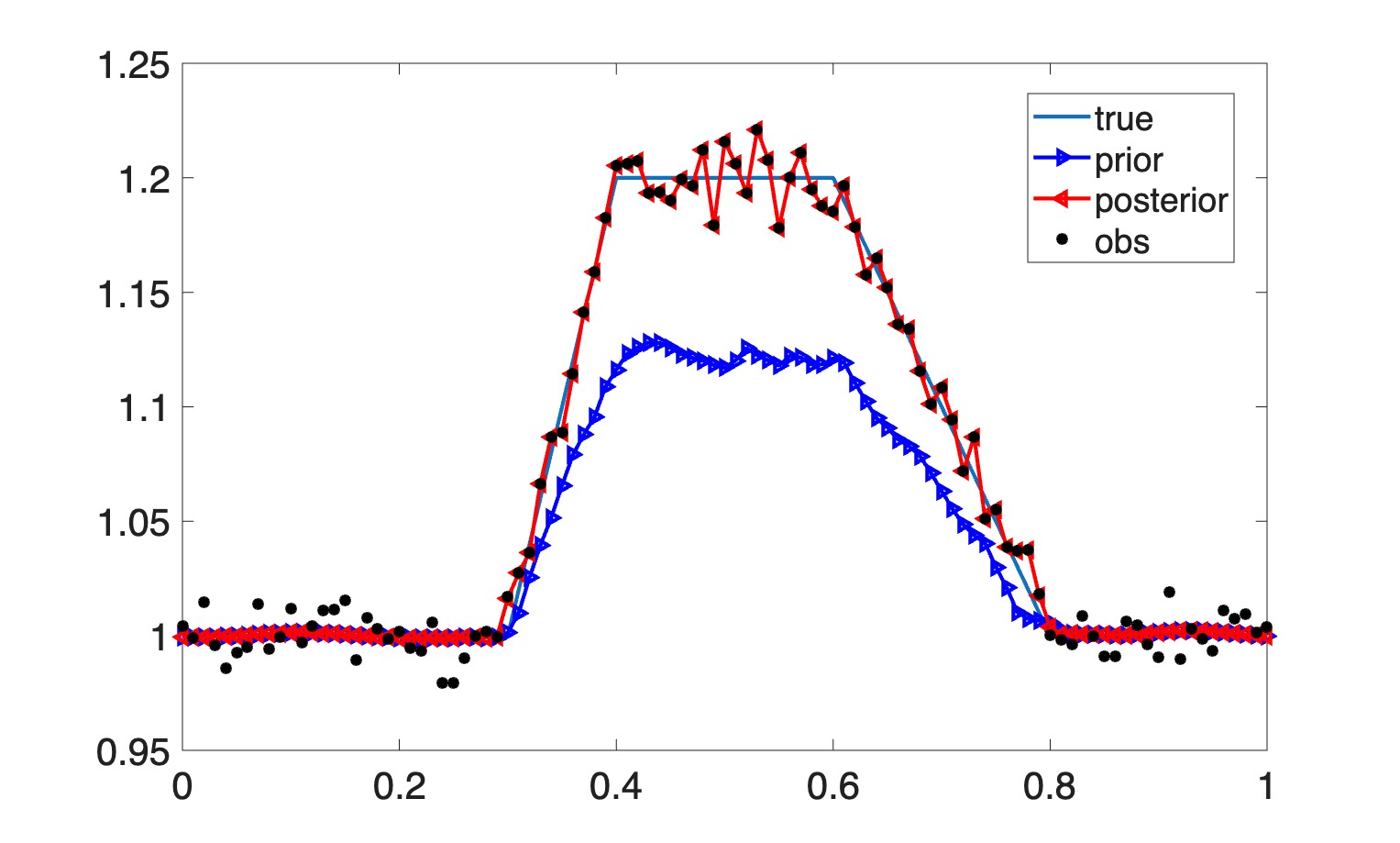} 
\caption{Numerical solutions for the posterior mean $\m$ for $W_C^D$ with $\alpha = 4, 6$ and $W_S^D$ with $(\vartheta,\varphi) =  (1/2,1), (1/2,2), (1,1), (1,2), (2,1)$, and $(2,2)$  respectively (left to right). 
}
\label{fig: grad plot}
\end{figure}

Finally, Figure \ref{fig: grad plot} displays the results at the 1D cross section of $u(0.4,y,2)$ for each of the eight parameters choices displayed in Table \ref{tab:cov and gradient all}. This particular cross section is chosen to show the impact of the weighting matrix in situations where the prior solution is very inaccurate, in this case due to diffusivity of the numerical PDE solver at $x = 0.4$. 

A desirable outcome would be one in which the  (inaccurate) prior does not dominate the posterior solution.  It should neither be entirely discarded.
Here we observe that choosing $W_C^D$ or $W_S^D$ with $(\vartheta,\varphi) = (1/2,1)$ yields  posterior solutions which are heavily influenced by the inaccurate prior.  By contrast, choosing $W_S^D$ with $(\vartheta,\varphi) = (1,2), (2,1)$ or $(2,2)$ produces posterior solutions that completely resemble the observations. A more balanced result is obtained when using $W_S^D$ with $(\vartheta,\varphi) = (1/2,2)$ or $(1,1)$, which is consistent with what is reported in Table \ref{tab:cov and gradient all}.

\subsection{Burgers' equation with dense observations}
\label{subsec:burgers-dense}

We use the 2D Burgers’ equation with dense observations to test the performance of the structurally-informed DA method on nonlinear advection problems with shock formation. 

\subsubsection*{Problem setup and numerical configuration} The 2D Burgers' equation within a square domain $[0,1]\times[0,1]$ for $t \in [0,2]$ is given by
\begin{equation}
\label{eq:2Dburgers}
 u_t + \left(\frac{u^2}{2}\right)_x +  \left(\frac{u^2}{2}\right)_y = 0. 
\end{equation}
Here we consider initial conditions 
\[ u(t=0) = u_0(x,y) = \begin{cases}
 1.2 & \text{ if } (x,y) \in [0.4, 0.6]\times [0.4,0.6], \\
 1 & \text{ otherwise},
\end{cases}
\]
with periodic boundary conditions  imposed in each direction.  As a closed form analytical solution is generally not available, we approximate the true solution using a more highly resolved fifth order WENO scheme and TVDRK3 time integration scheme with $\Delta x = \Delta y = 10^{-3}$ and $\Delta t = 2 \times 10^{-5}.$

The numerical solution is simulated using $\Delta x = \Delta y = 10^{-2}$ and $\Delta t = 2 \times 10^{-3}.$ The data assimilation setup is identical to that used for the linear advection equation, except that observations are now available at intervals of $\Delta t_{obs}=5\Delta t = 10^{-2}$. We compare the posterior solutions obtained using  $W_C^D$ in \eqref{eq: Cd} with inflation parameter $\alpha = 4$ with those obtained using $W^D_S$ in \eqref{eq: Sd} for fixed aggregate parameter $\varphi = 1$ and moment parameters $\vartheta=1/2, 1$, and $2$. The tuning parameter $\wt{\beta}$ is chosen as in the linear advection equation (see Table \ref{tab:cov and gradient all}).

\subsubsection*{Performance Comparison} 
\begin{figure}[h!]
\centering
\includegraphics[width=0.24\textwidth]{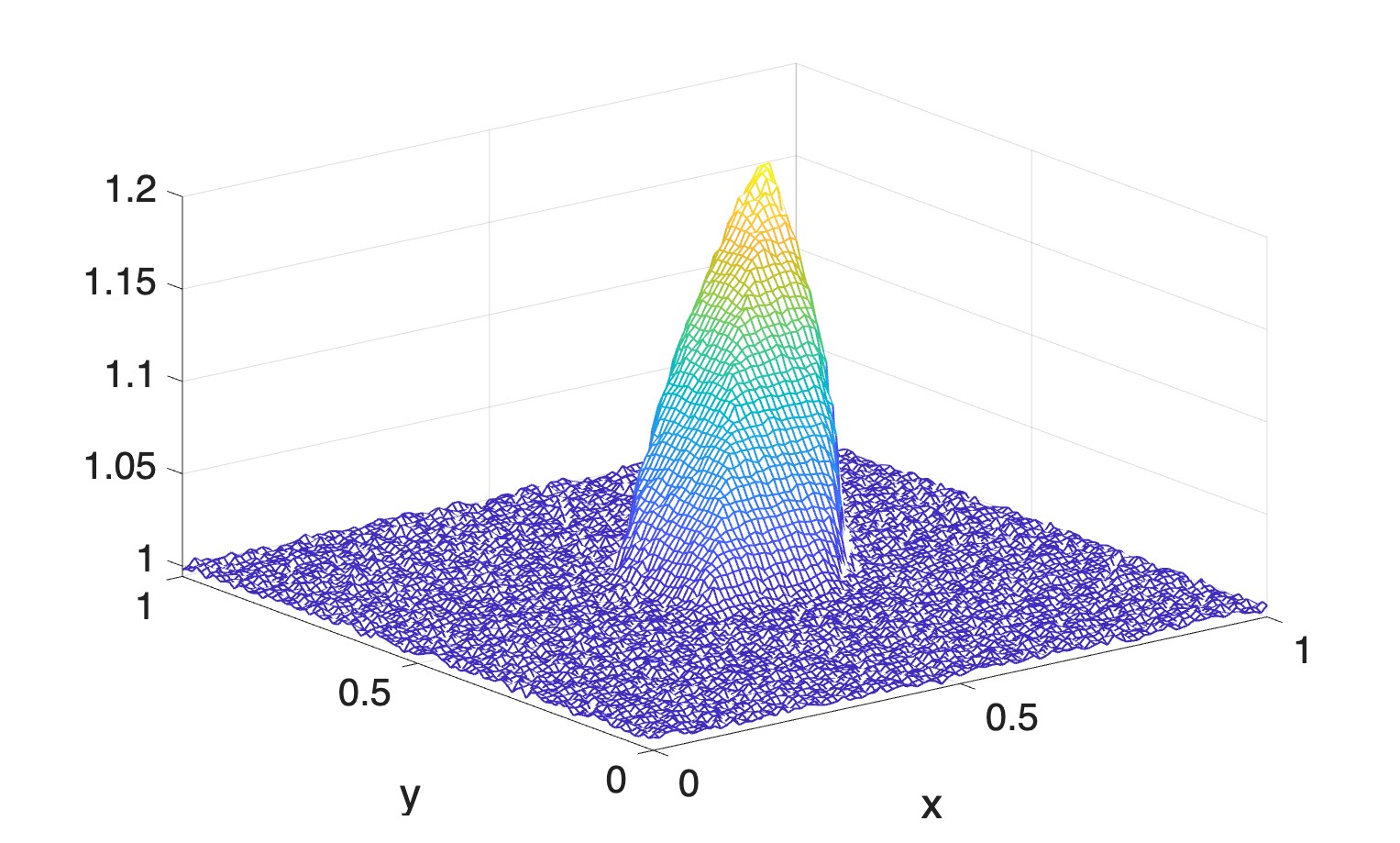}
\includegraphics[width=0.24\textwidth]{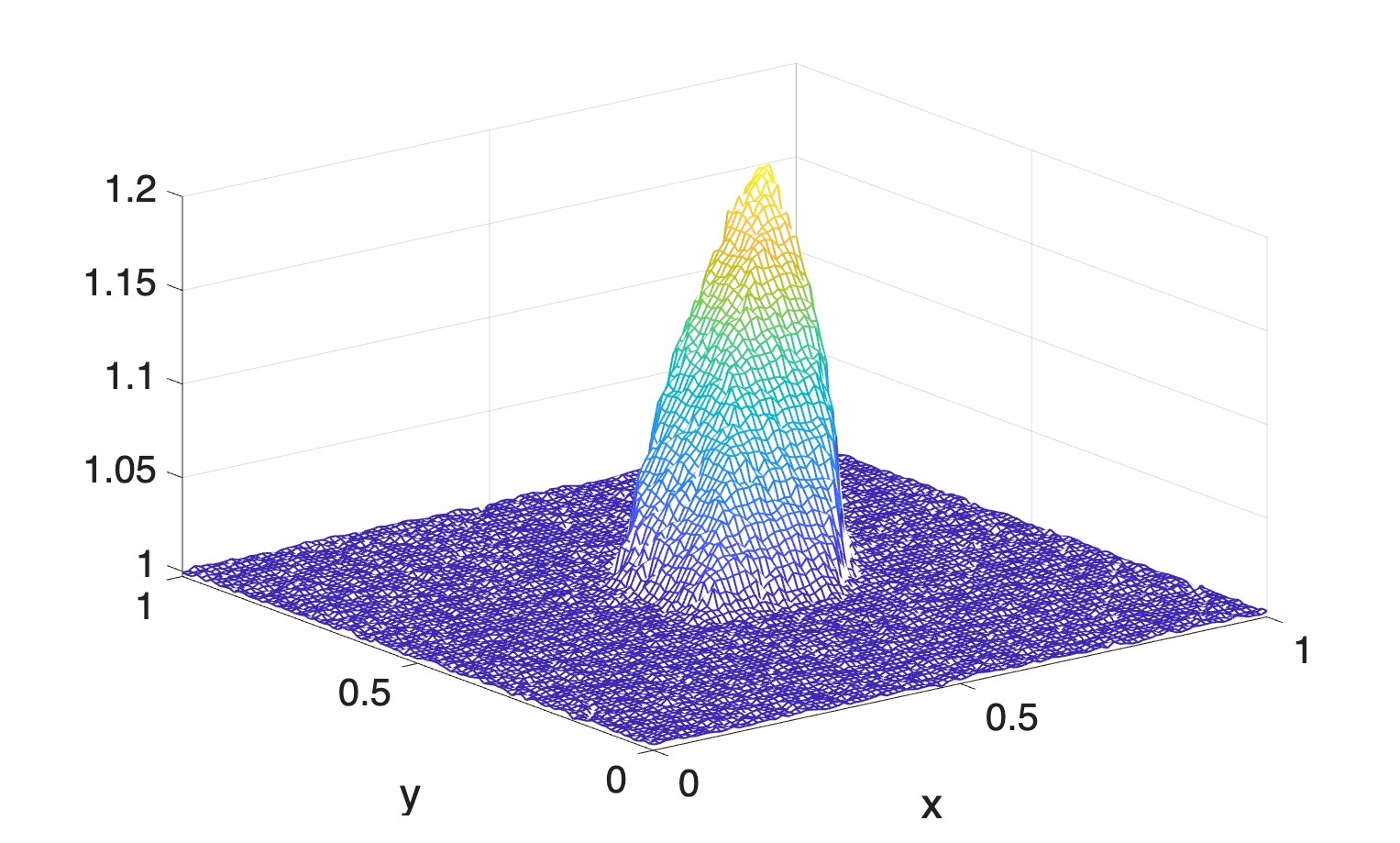}
\includegraphics[width=0.24\textwidth]{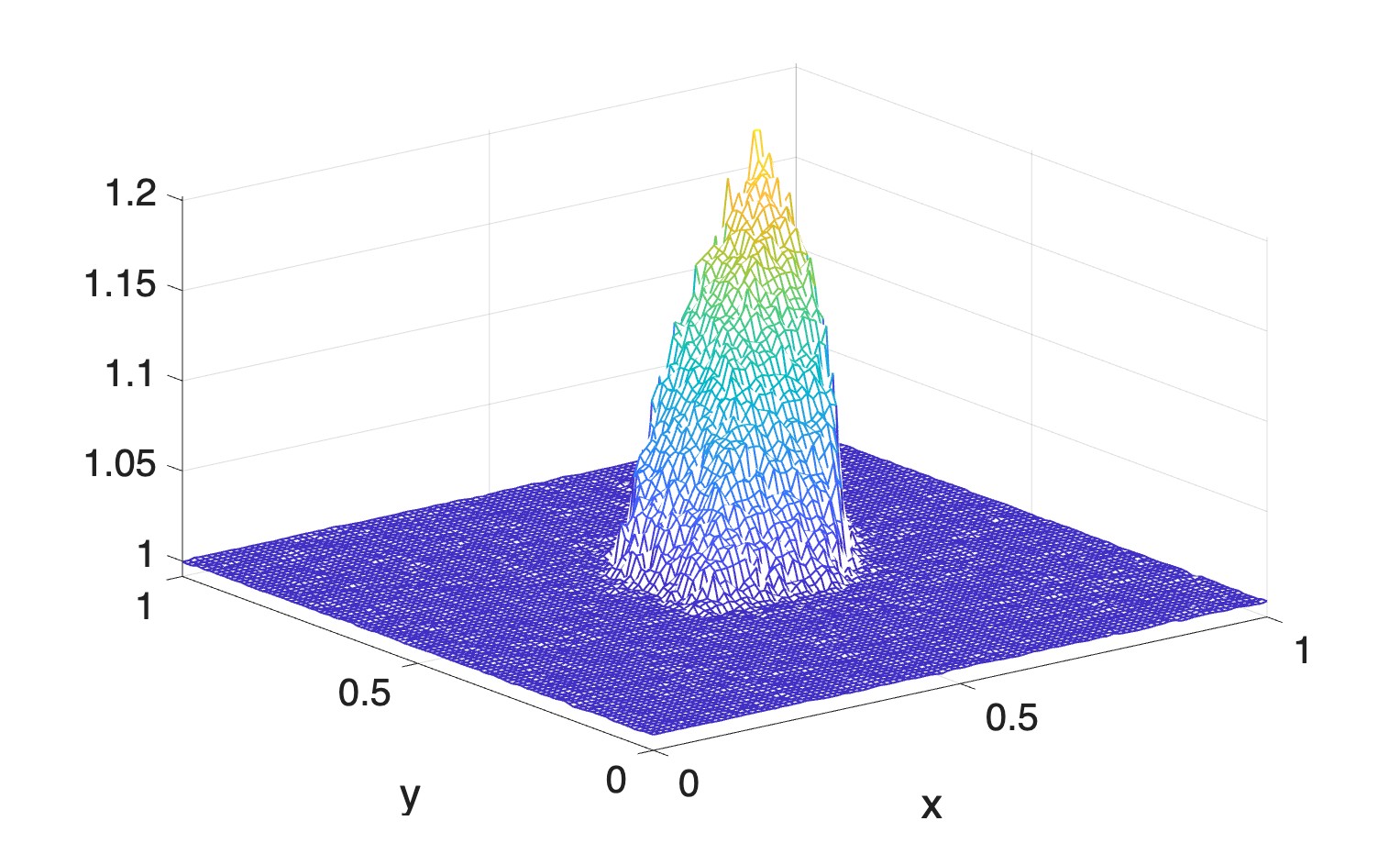}
\includegraphics[width=0.24\textwidth]{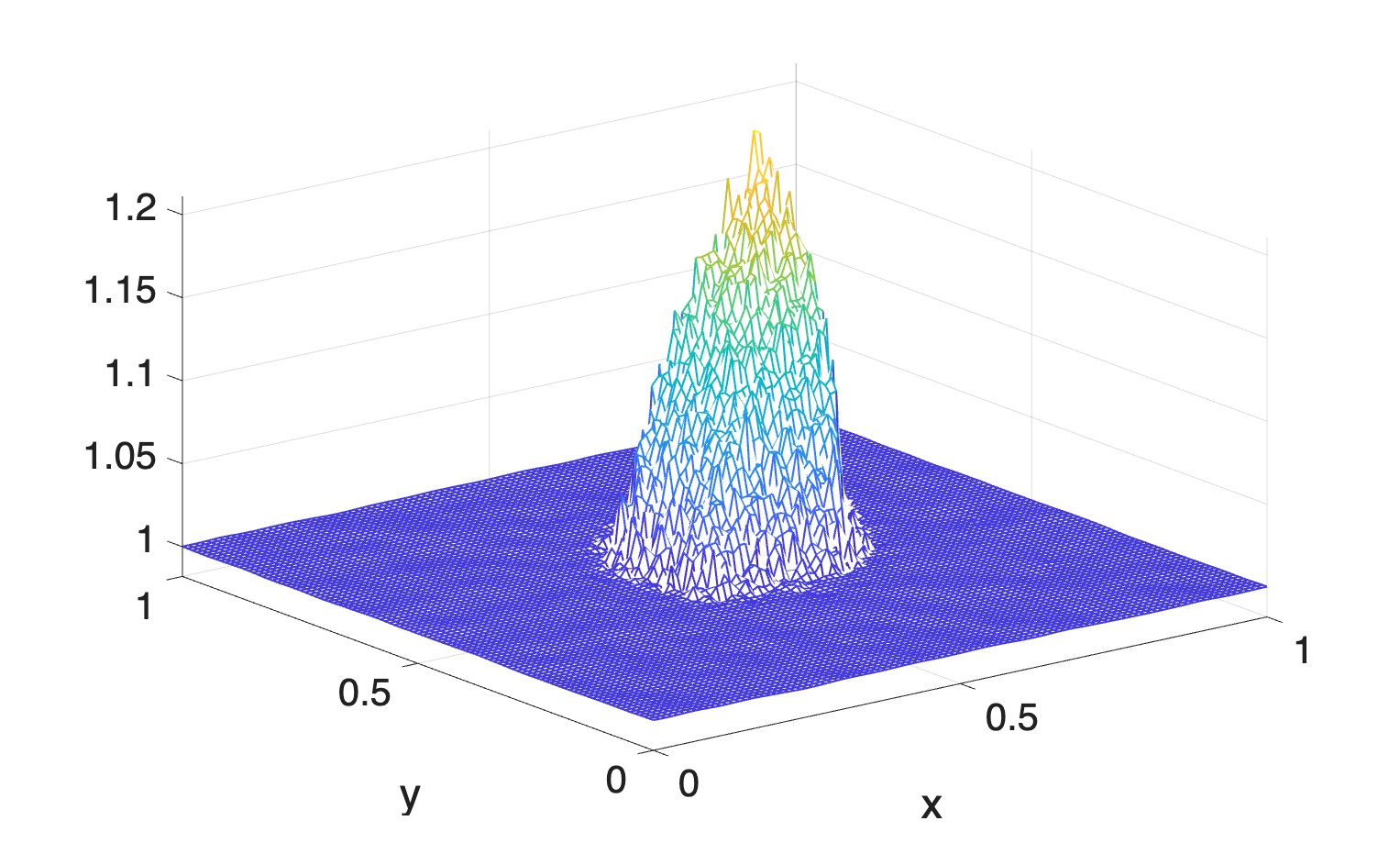}
\caption{Posterior mean solution $\m$ \eqref{eq:mean_update} at time $t=2$ obtained by $W_C^D$ in \eqref{eq: Cd} (left), and $W_S^D$ in \eqref{eq: Sd} with $(\vartheta,\varphi)=$ (1/2, 1) (middle left), (1,1) (middle right), and (2,1) (right).}
\label{fig: burgers mesh}
\end{figure}

\begin{figure}[h!]
\centering
\includegraphics[width=0.24\textwidth]{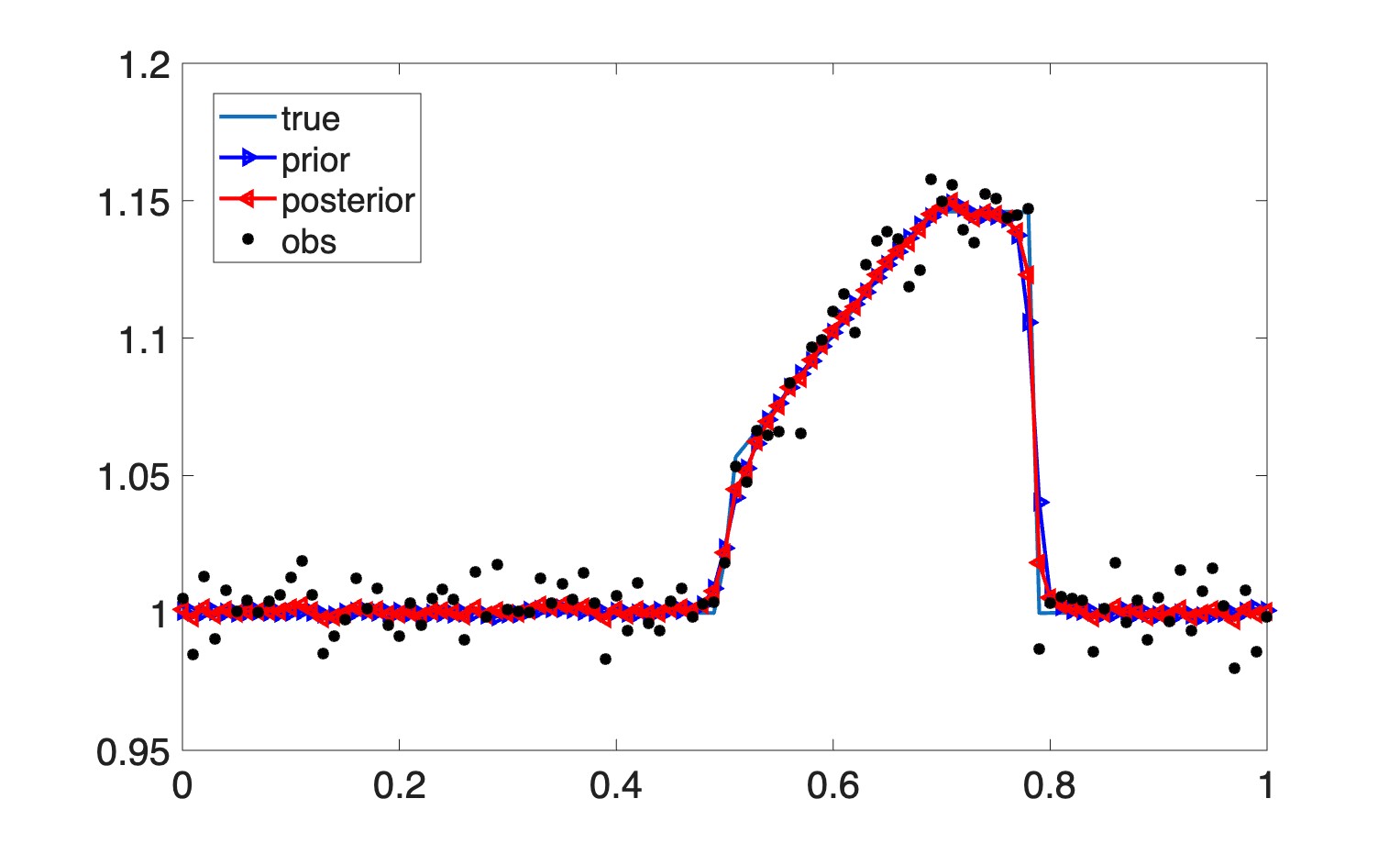}
\includegraphics[width=0.24\textwidth]{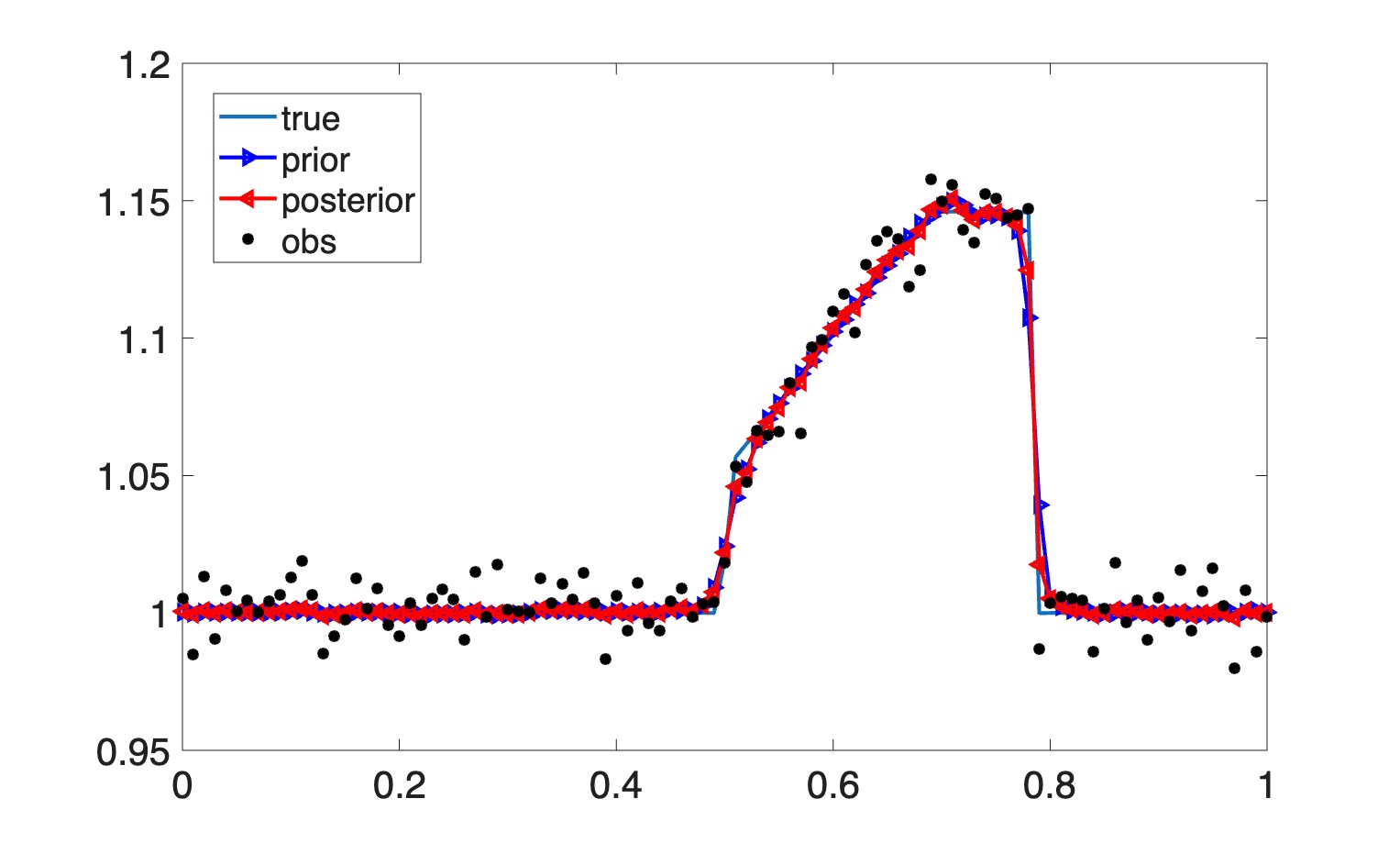}
\includegraphics[width=0.24\textwidth]{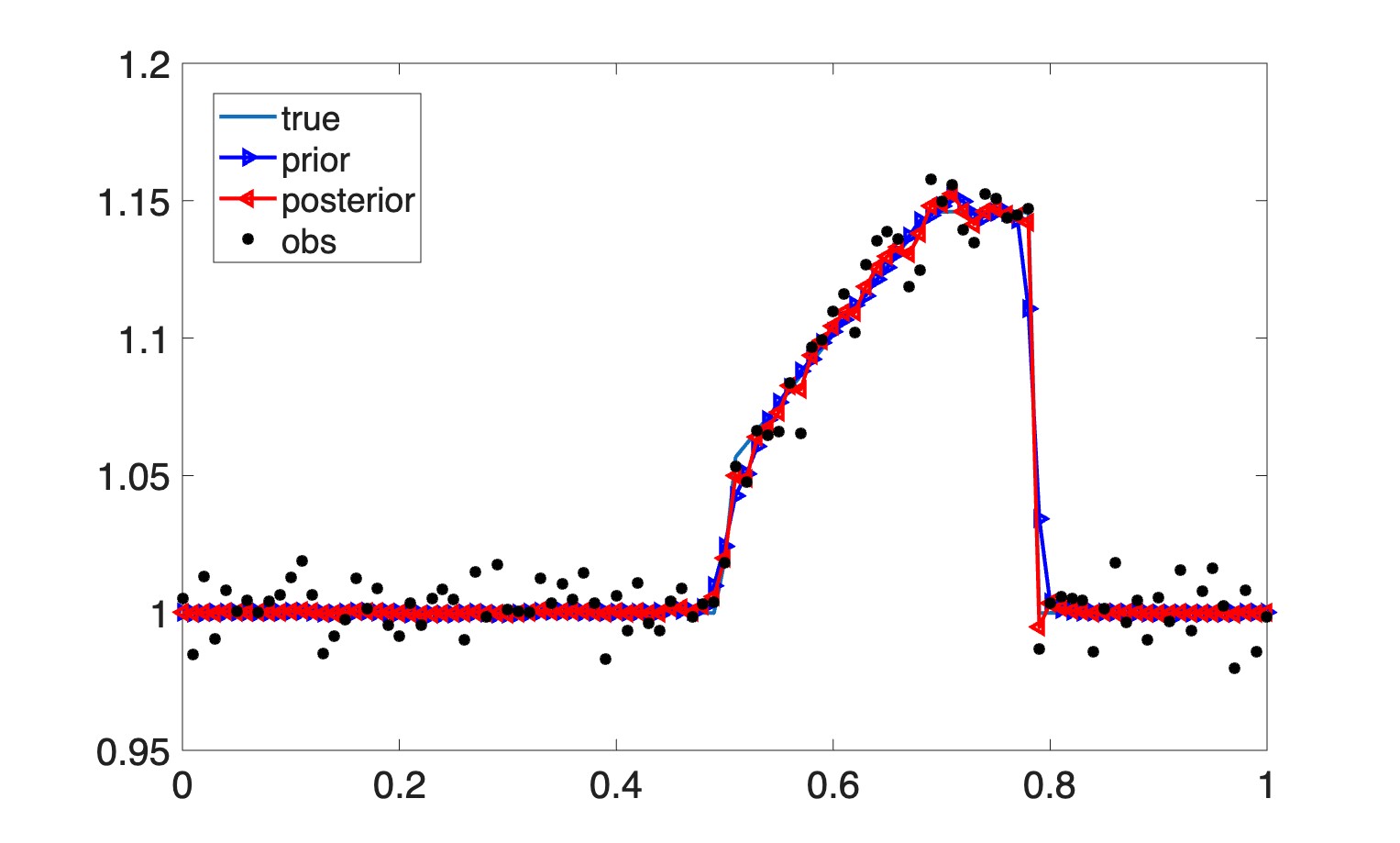}
\includegraphics[width=0.24\textwidth]{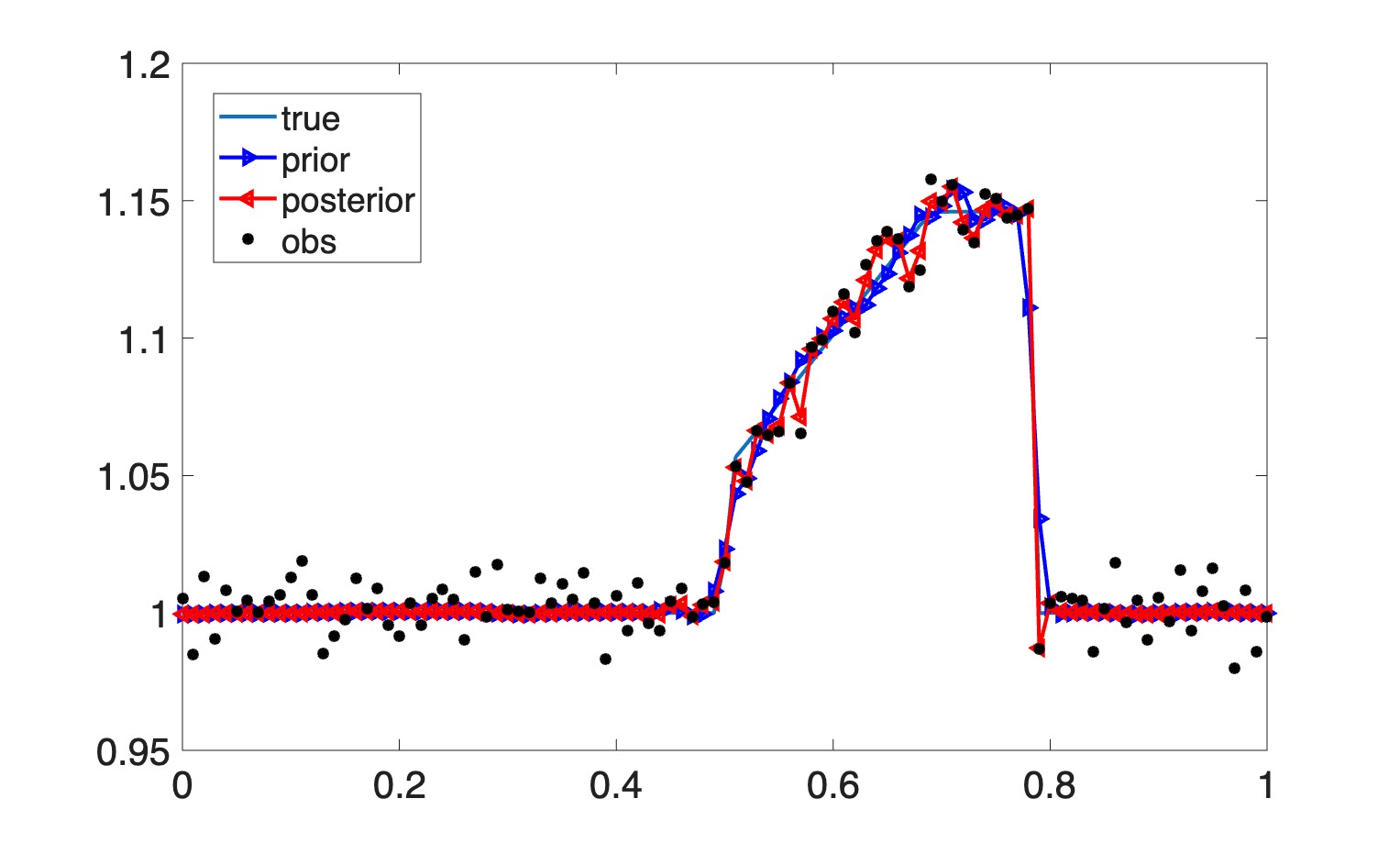}
\\
\includegraphics[width=0.24\textwidth]{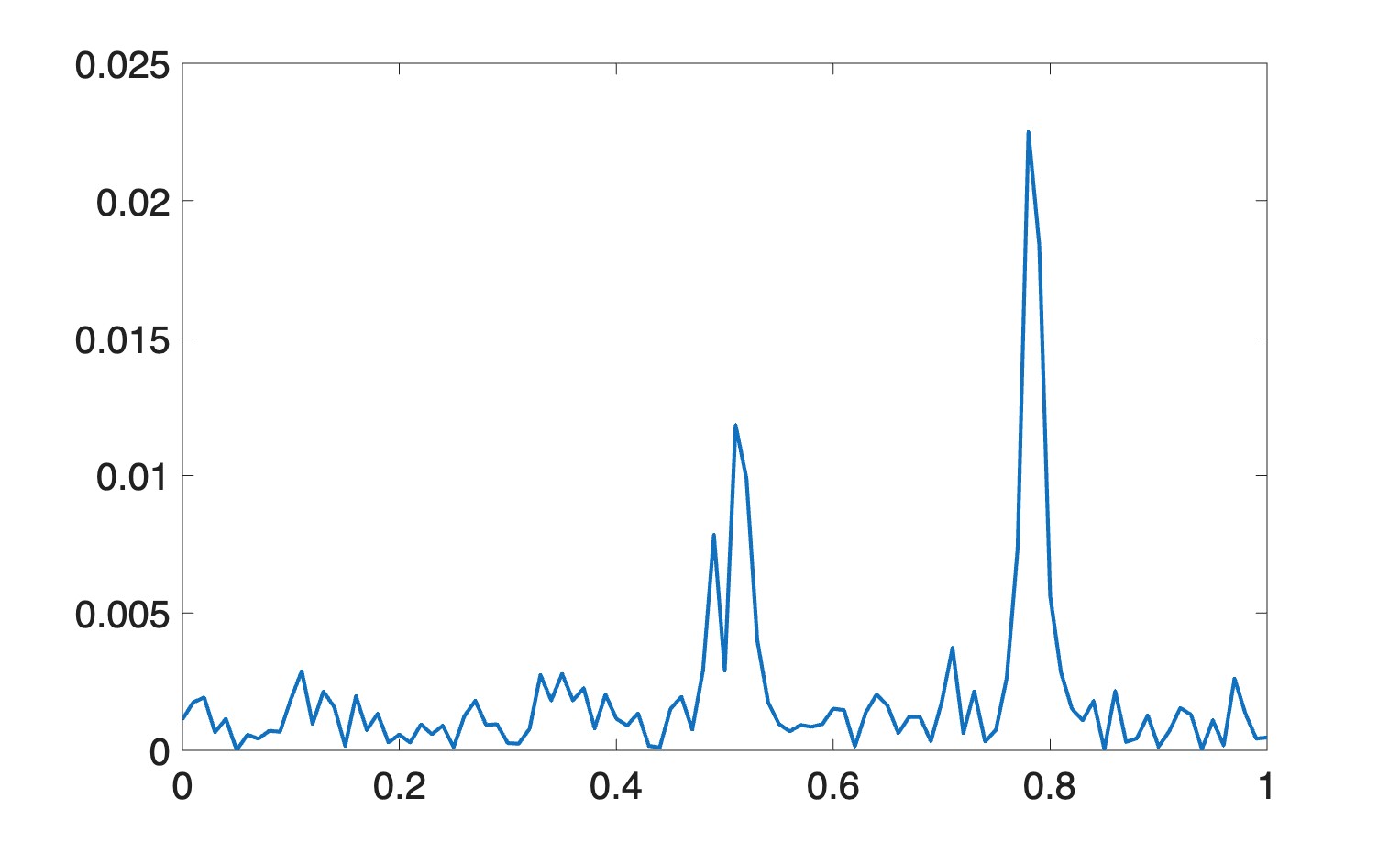}
\includegraphics[width=0.24\textwidth]{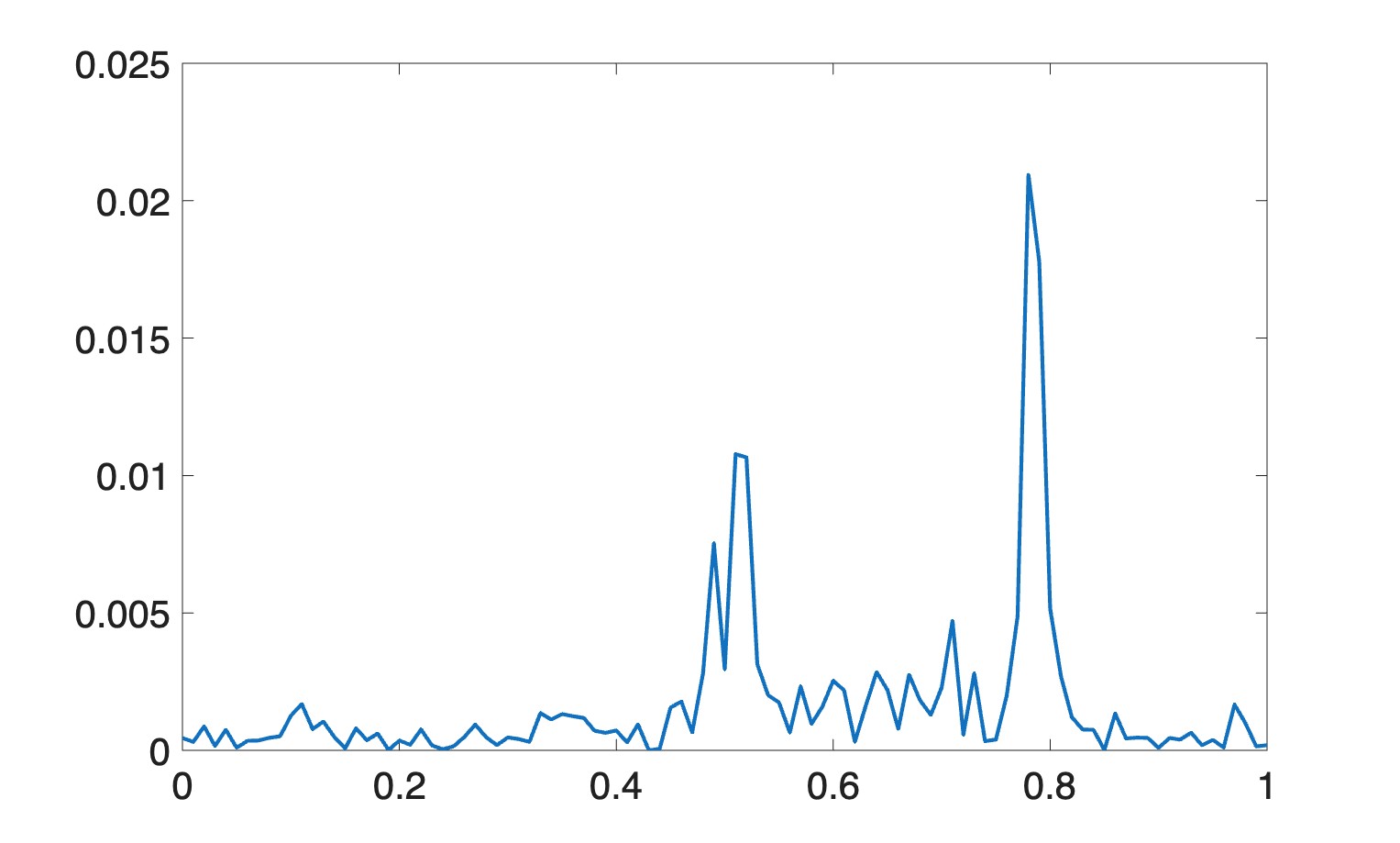}
\includegraphics[width=0.24\textwidth]{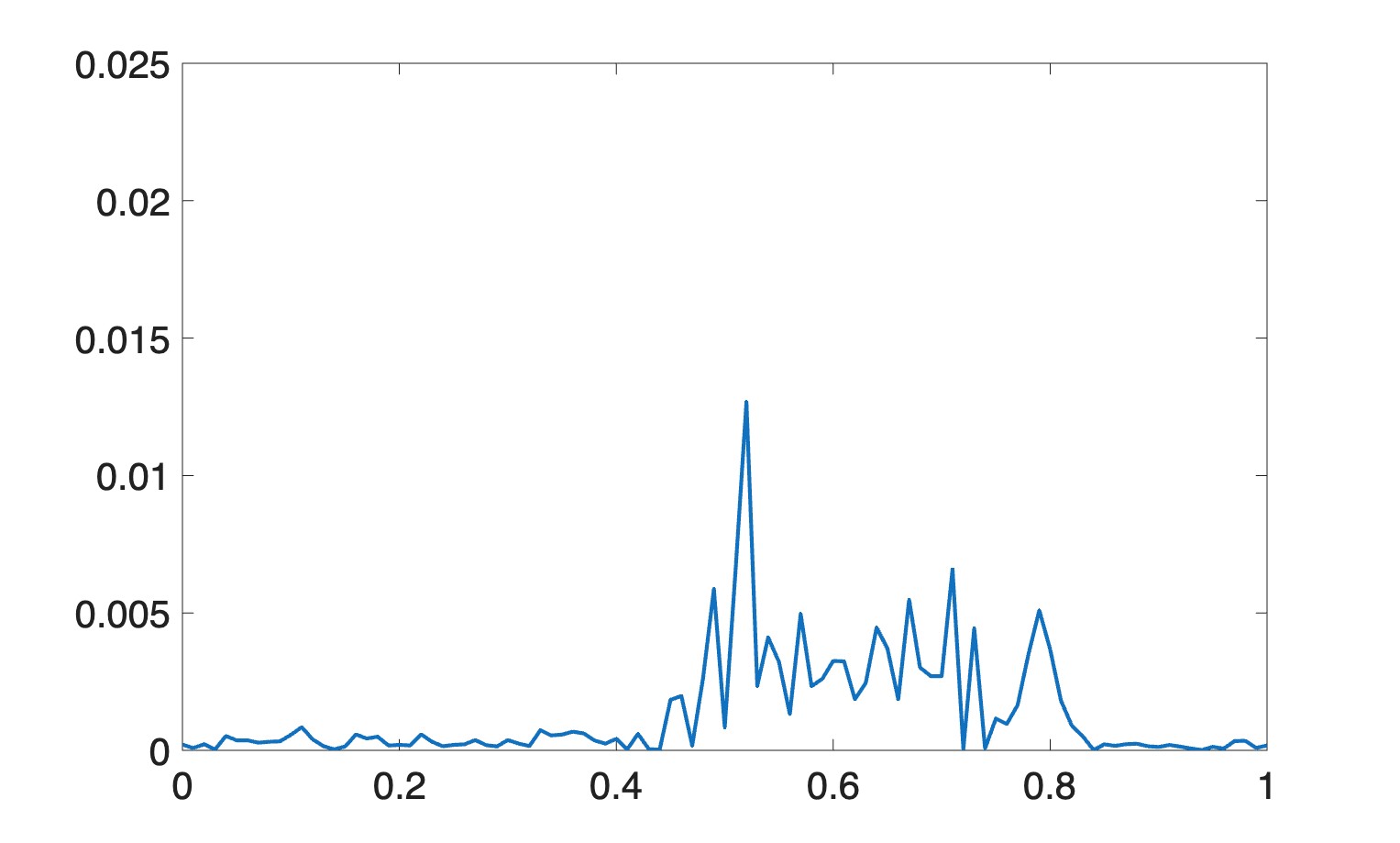}
\includegraphics[width=0.24\textwidth]{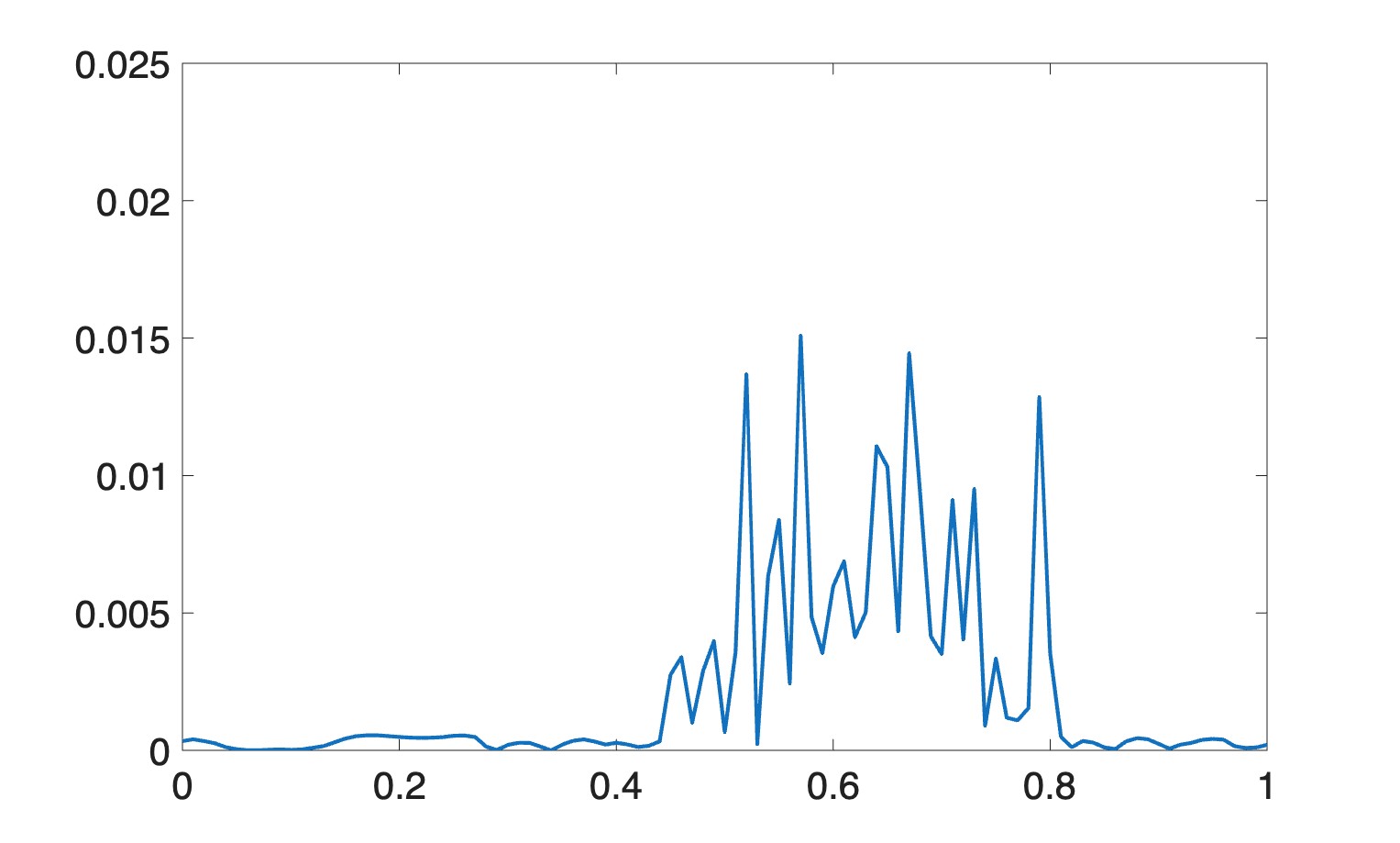}
\caption{Numerical solutions for the posterior mean $\m$ (top) and the corresponding pointwise error $err$ (bottom) for $W_C^D$ and $W_S^D$ with $(\vartheta, \varphi)=$ (1/2, 1), (1,1) , and (2,1) (from left to right) respectively.}
\label{fig: burgers oned}
\end{figure}

Figure \ref{fig: burgers mesh} compares the DA mean solution $\m$ in \eqref{eq:mean_update} at time $t = 2$ for $W_C^D$ in \eqref{eq: Cd}  and  $W_S^D$ in \eqref{eq: Sd}, while Figure \ref{fig: burgers oned} compares the corresponding cross section solution at $y = 0.7$.  Observe that $W_S^D$ yields overall better accuracy which is especially noticeable in smooth regions.  It is also able to distinguish  smooth from discontinuous regions, especially for increasing values of $\vartheta$.  In particular, the prior has greater influence in smooth regions, while the observational data has more impact in  discontinuous regions.  When $\vartheta = 2$, the solution relies too heavily on the observations. Hence we see that the moment parameter must be chosen carefully to balance the prior and observations contributions effectively across all regions. In this case $\vartheta=1$ effectively achieves this balance. 

\subsubsection*{Time-series error metrics} 
\begin{figure}[h!]
\centering
\includegraphics[width=0.3\textwidth]{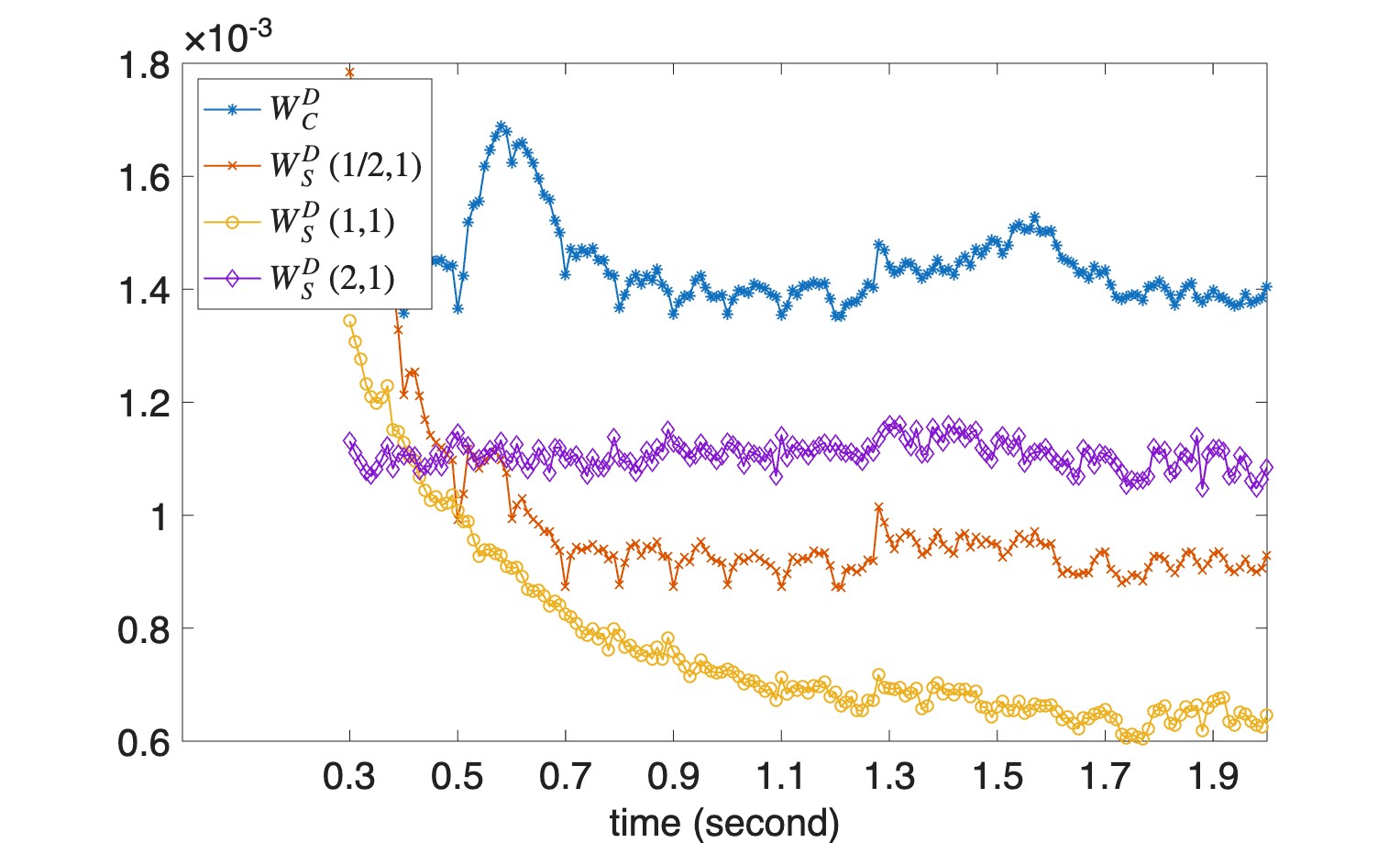}
\includegraphics[width=0.3\textwidth]{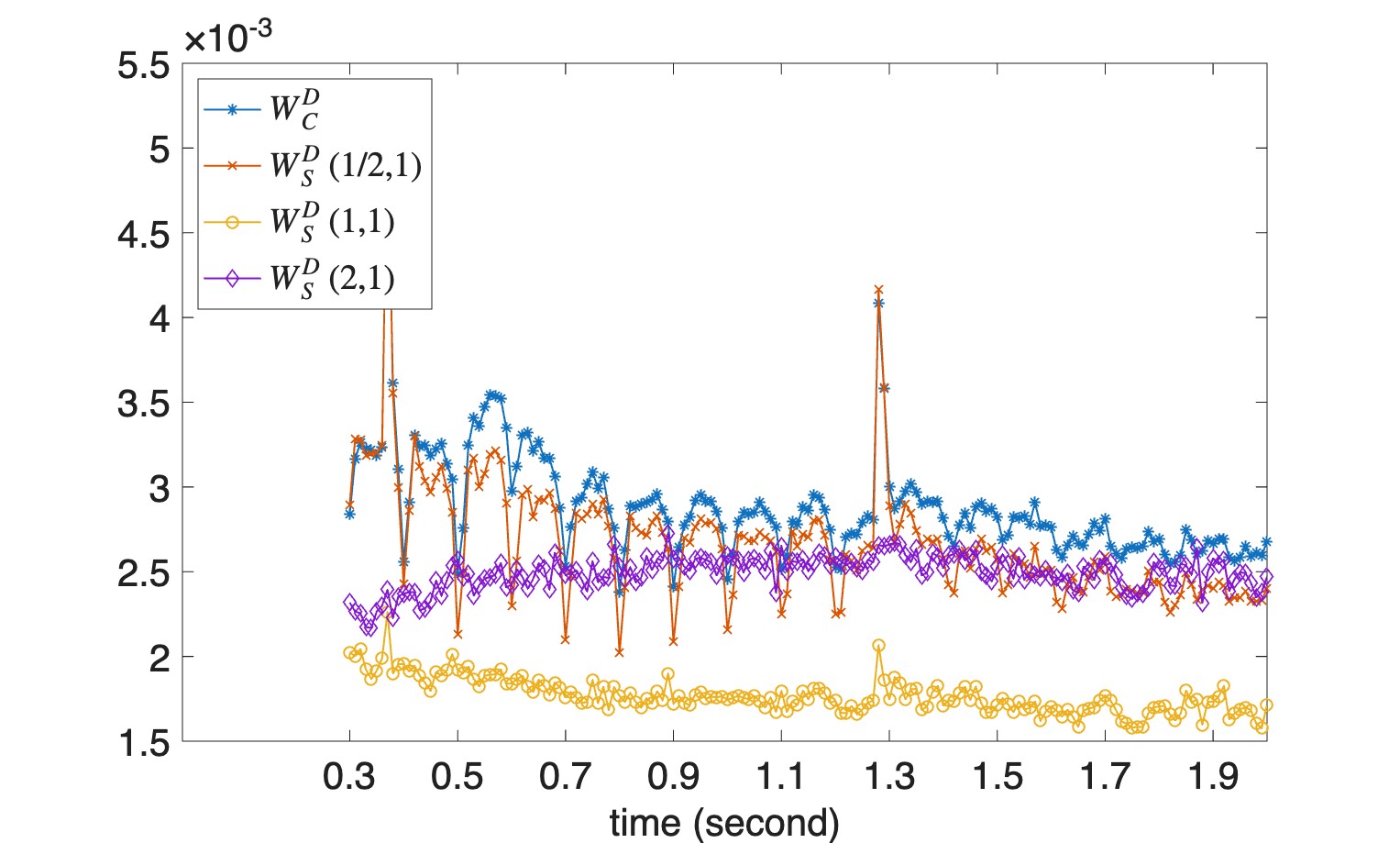}
\includegraphics[width=0.3\textwidth]{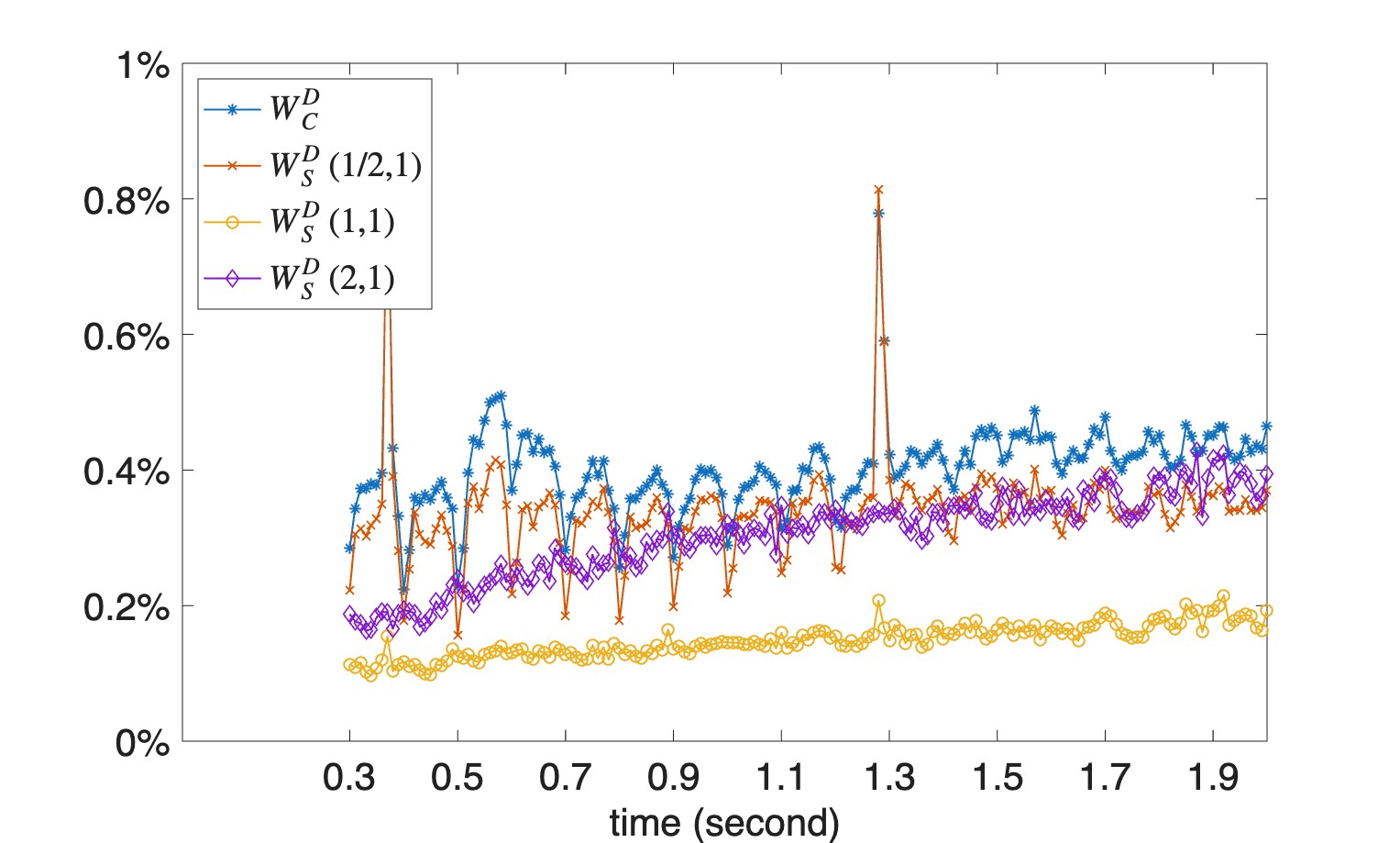}
\caption{Relative errors $err_{l_1}$ \eqref{eq:l1err} (left) and  $err_{l_2}$ \eqref{eq:l2err} (middle), along with the complementary pattern correlation $1-Pcorr$ \eqref{eq:pcorr} (right), obtained for the posterior mean  using $W_C^D$ and $W_S^D$ with momentum parameters $\vartheta = 1/2, 1, 2$.}
\label{fig: burgers err}
\end{figure}

Figure \ref{fig: burgers err} displays the relative errors $err_{l_1}$, $err_{l_2}$, and the complementary pattern correlation $1-Pcorr$, obtained using \eqref{eq:l1err}, \eqref{eq:l2err}, and \eqref{eq:pcorr}, respectively, for the four choices of parameters used in Figure \ref{fig: burgers mesh} and Figure \ref{fig: burgers oned}, over the time domain $[0.3,2]$. The initial time interval $[0,0.3]$ is omitted to avoid bias from initialization. As with the results reported in Figure \ref{fig: burgers mesh} and Figure \ref{fig: burgers oned}, $W_S^D$ consistently outperforms $W_C^D$, with $(\vartheta,\varphi)=(1,1)$ achieving the smallest relative errors and complementary pattern correlation.

\subsection{Sparse observation pattern and implementation of localization and directional refinement}
\label{subsec:sparse-impl}

\begin{figure}[h!]
\centering
\includegraphics[width=0.6\linewidth]{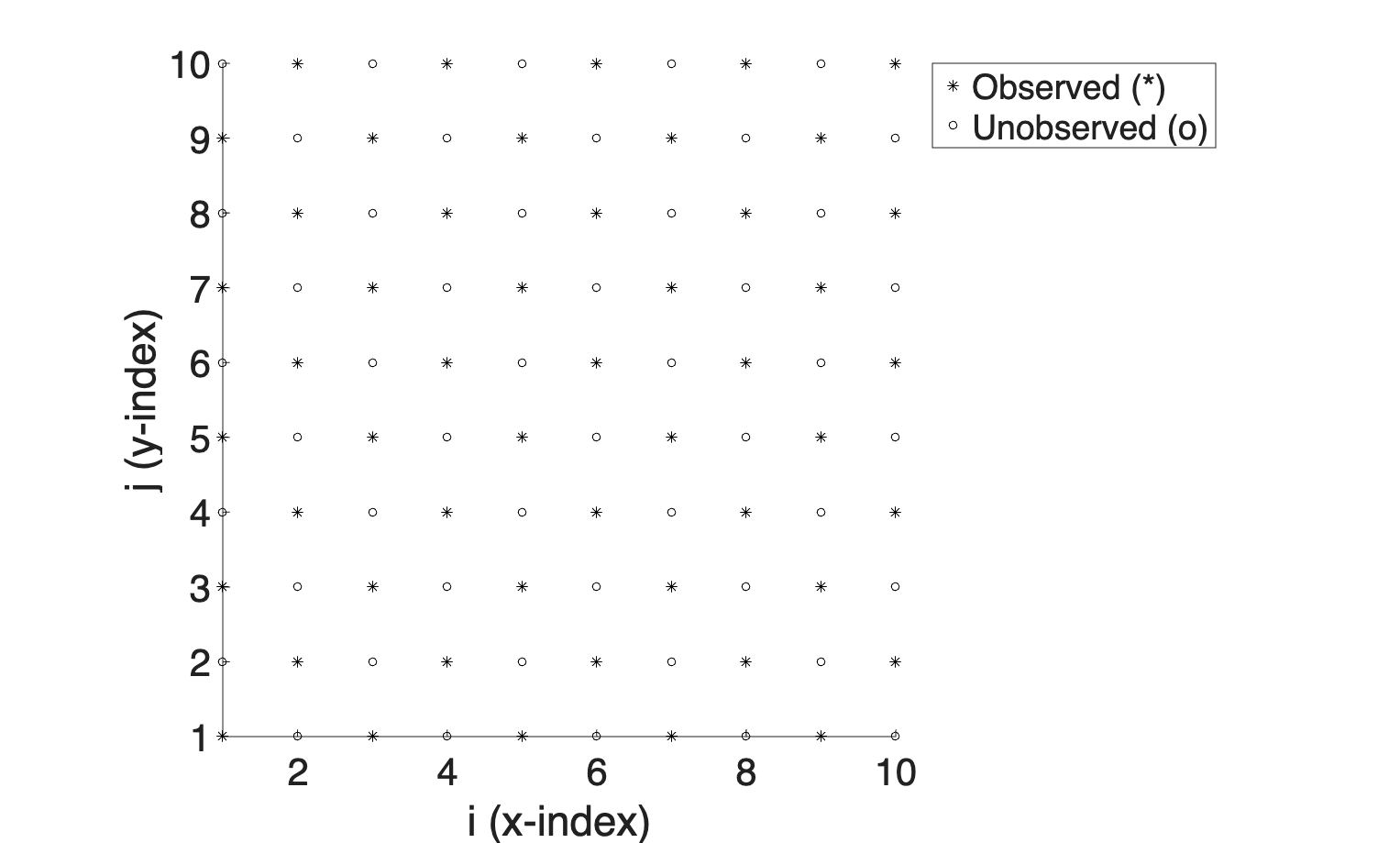}
\caption{Checkerboard sparsity pattern illustration.}
\label{fig: checkerboar}
\end{figure}

We now focus on the more challenging case of sparse observations.  To this end we first detail the design of the sparse observation pattern and describe how the localization and directional refinement strategies can be efficiently implemented. We adopt a classic checkerboard sparsity pattern on the 2D grid: if $i$ and $j$ denote the horizontal and vertical grid indices, observations are placed at points satisfying
$$\text{Observation at }(i,j) \Leftrightarrow (i+j) \text{ mod } 2 = 0,$$
resulting in an alternating pattern both horizontally and vertically across the domain, as illustrated in Figure \ref{fig: checkerboar}.

Following \eqref{eq:localization}, we now seek a localization matrix $\cT$ that is designed to accommodate this checkerboard arrangement. In this regard, we note that a grid point either contains an observation or its  four immediate neighbors (up, down, left, right) do.  Therefore, to efficiently obtain localization that allows correlations between neighbors, we design $\cT$ to be a sparse five-banded symmetric matrix.  Specifically, in the flattened vector representation of the domain (using column-major stacking according to \eqref{eq: index mapping}), the main diagonal of $\cT \in \bbR^{n_x n_y \times n_x n_y}$ represents self-correlations. The first upper and lower off-diagonals correspond to neighbors one index apart in the flattened vector, which for points in the same column represent vertical neighbors.\footnote{Since the flattened index does not distinguish column boundaries, modulus checks on the index are applied to exclude spurious connections across columns.} Finally, as their relative index shift directly matches one full column's stride, the $n_x$-th upper and lower off-diagonals correspond to left and right horizontal neighbors. In terms of the assigned values, we set 1 on the main diagonal (self-correlations) and 0.5 on the off-diagonal entries corresponding to valid immediate neighbors. This tapering enforces distance-based decay, acting as a simplified alternative to the Gaspari-Cohn kernel. We are now able to construct the localization matrix with entries given by
\begin{equation} \label{eq: cT}
\cT_{m,m'} = 
\begin{cases}
1, & m = m', \\
0.5, & m'-m=1, \text{and } m \text{ mod } n_x \neq 0, \\
0.5, & m-m'=1, \text{and } m' \text{ mod } n_x \neq 0, \\
0.5, & |m-m'|=n_x, \\
0, & \text{otherwise}.
\end{cases}
\end{equation}

We next describe the construction of the binary directional refinement mask $\cM$ introduced in \eqref{eq:Mask}, which complements the localization $\cT$ in the weighting matrix design.\footnote{We note that utilizing $\cM$ may also be appropriate in the dense observation case for more complex dynamical systems.} Observe first that the first upper and lower off-diagonals correspond to immediate neighbors in the $x$-direction, and the $n_x$-th upper and lower off-diagonals are related to immediate neighbors in the $y-$direction. This enables an efficient directional refinement by computing pairwise discrepancies for all relevant pairs in a single step per direction. Second,  recall from Remark \ref{rem:Mimple} that $\cM$ is sparse, with non-zero values only on the four off-diagonals corresponding to $\cT$ in \eqref{eq: cT}. We now proceed as follows: For the $x$-direction, we compute a matrix of one-sided discrete directional derivatives on the prior mean:
\begin{equation}\label{eq: direction derivative x}
D^x_{i,j} = \frac{\left|\wh{\m}_{i+1,j} - \wh{\m}_{i,j} \right|}{\Delta x}, \quad i=1, ..., n_x-1, \ j=1, ..., n_y. 
\end{equation}
which yields a derivative field of size $(n_x-1) \times n_y$, indexed by $i$ for row and $j$ for column. This corresponds to discrepancies between left-right neighbor pairs in the physical space (i.e. \eqref{eq: directional discrepancy} with $m$ and $m'$ corresponding to $i+1,j$ and $i,j$), and provides information to refine the first upper and lower off-diagonals (of size $n_x n_y -1$) in the weighting matrix. To ensure alignment with the vectorized structure, we first pad $D^x$ by appending an additional row of zeros at the bottom to account for artificial connections introduced by column-stacking. We then remove the last entry to eliminate extraneous connections at the last column boundary. This resulting modified matrix is then flattened column-wise to yield a vector $D^x_\text{flat}$, which is compared against the threshold $d_\text{thresh}$ to determine binary refinement $\cM$. We denote the corresponding vector for the first upper and lower off-diagonals by $M^{x}_{\text{offd}} \in \{0,1\}^{n_x n_y -1}$, with entries given by
\begin{equation}\label{eq: Mxoff}
\left( M^{x}_{\text{offd}}\right)_k =
\begin{cases}
0, & \left( D^x_\text{flat} \right)_k > d_\text{thresh}, \\
1, & \left( D^x_\text{flat} \right)_k \le d_\text{thresh},
\end{cases}
\quad k = 1, \cdots, n_x n_y-1,
\end{equation}
reflecting whether directional discrepancies exceed the user-defined threshold. These values are then inserted into $\cM$ at the appropriate positions along the first upper and lower off-diagonals, controlling the influence between horizontally adjacent grid points in physical space base on the detected gradient strength.

The $y-$direction directional derivatives matrix is similarly obtained as
\begin{equation}\label{eq: direction derivative y}
D^y_{i,j} = \frac{\left|\wh{\m}_{i,j+1} - \wh{\m}_{i,j}\right|}{\Delta y}, \quad i=1, \dots, n_x, \ j=1, \dots, n_y-1,
\end{equation}
yielding a derivative field of size $n_x \times (n_y-1)$. Flattening this matrix column-wise produces a vector $D^y_\text{flat}$, whose entries are then compared with the threshold $d_\text{thresh}$ to construct the $n_x$-th upper and lower diagonals of the reduced binary mask $\cM$ in the same way as above. Finally, we write $\cM$ as a sparse matrix with five-banded symmetric structure:
\begin{equation}\label{eq: cM eff}
\cM=
\begin{pmatrix} 
1 & \left(M^{x}_{\text{offd}}\right)_1 & \cdots & \left(M^{y}_{\text{offd}}\right)_1 & \cdots & 0\\
\left(M^{x}_{\text{offd}}\right)_1 & \ddots & \ddots & \ddots & \ddots & \vdots \\
\vdots & \ddots & \ddots & \ddots & \ddots & \left(M^{y}_{\text{offd}}\right)_{n_x\times (n_y-1)}  \\
\left(M^{y}_{\text{offd}}\right)_1 & \ddots & \ddots & \ddots & \ddots & \vdots \\
\vdots & \ddots & \ddots &  \ddots & \ddots &\left(M^{x}_{\text{offd}}\right)_{n_xn_y-1} \\
0 & \cdots  &   \left(M^{y}_{\text{offd}}\right)_{n_x\times (n_y-1)} & \cdots & \left(M^{x}_{\text{offdiag}}\right)_{n_xn_y-1} & 1
\end{pmatrix}.    
\end{equation}
\subsection{Burgers' equation with sparse observations}

Our final numerical experiment investigates the performance of the structurally informed data assimilation framework in  a sparse observation environment (Figure \ref{fig: checkerboar}) for 2D Burgers’ equation in \eqref{eq:2Dburgers}.  We use the same initial and boundary conditions, along with $\gamma = 0.005$ for $\Gamma = \gamma^2 I$ in \eqref{model: observation}.   The observation data are again generated using a highly resolved fifth order  WENO  scheme and TVDRK3 for time integration. 

Our experiments are conducted using the same numerical solver and ensemble settings as in the the dense observation case, and we compare results using (1) the covariance-based weighting $W_C$ in \eqref{eq: Wc} with $\alpha = 4$; (2) the gradient-based weighting $\wh{W}_S$ in \eqref{eq: Ws0} without correlation refinement; and (3)  the gradient based weighting with directional structure-based correlation refinement ${W}_S$ in \eqref{eq: W in mask}. We use localization matrix $\mathcal{T}$ in \eqref{eq: cT} for all three  weighting matrix constructions and  choose $(\vartheta,\varphi)=(1,1)$ for both versions of $\wh{W}_S$.\footnote{Test results (not shown here) confirm that using other parameter choices yield results analogous to the dense observation case, with $(\vartheta,\varphi) = (1,1)$ and $(\vartheta,\varphi) = (\frac{1}{2},2)$ being optimal.} Finally, we set $\wt{\beta}=1.0e-04$ in \eqref{eq:betat}, and $d_\text{thresh}=4$ in \eqref{eq: mod_coeffs}.  

\begin{figure}[h]
\centering
\includegraphics[width=0.24\textwidth]{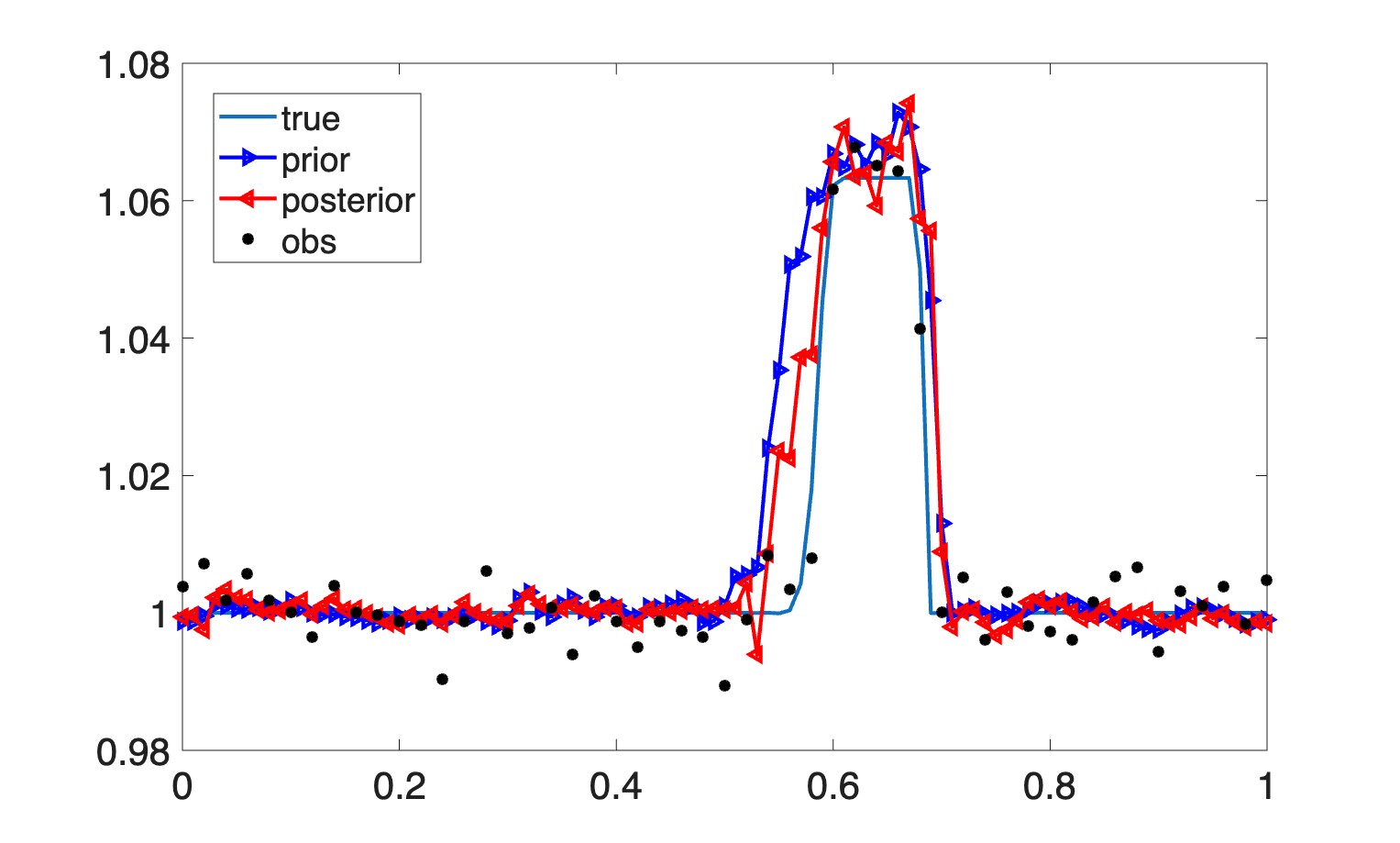}
\includegraphics[width=0.24\textwidth]{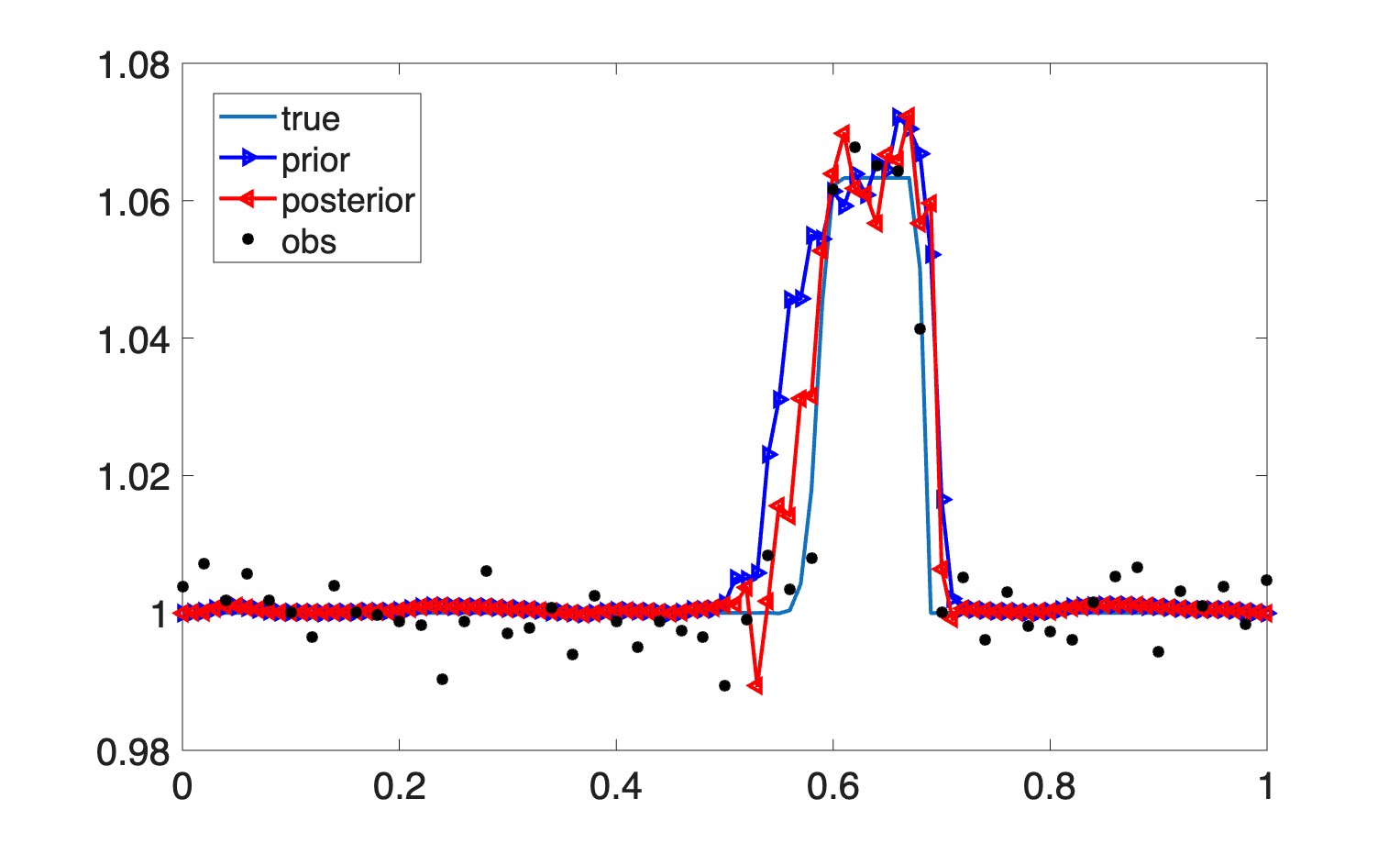}
\includegraphics[width=0.24\textwidth]{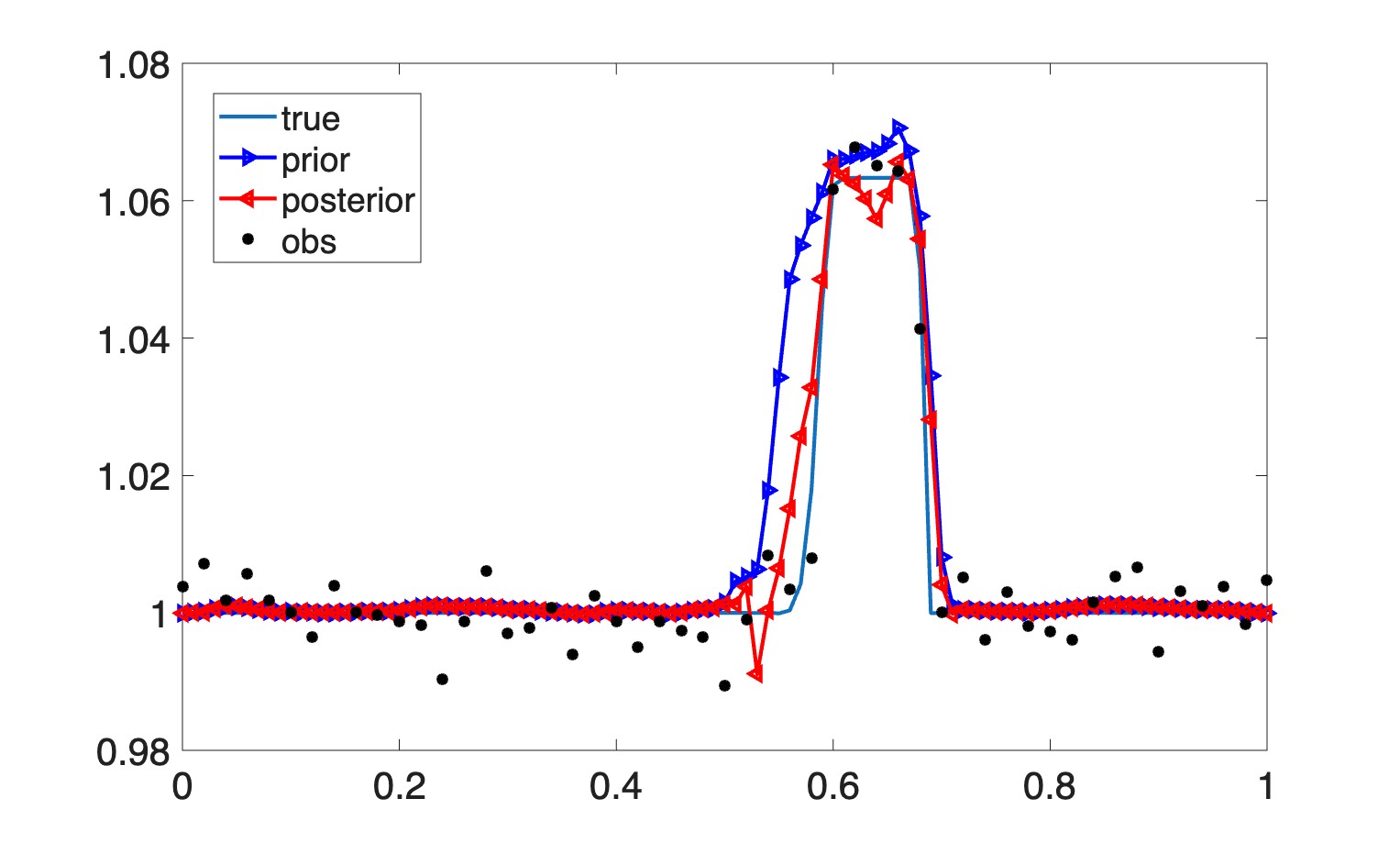}
\\
\includegraphics[width=0.24\textwidth]{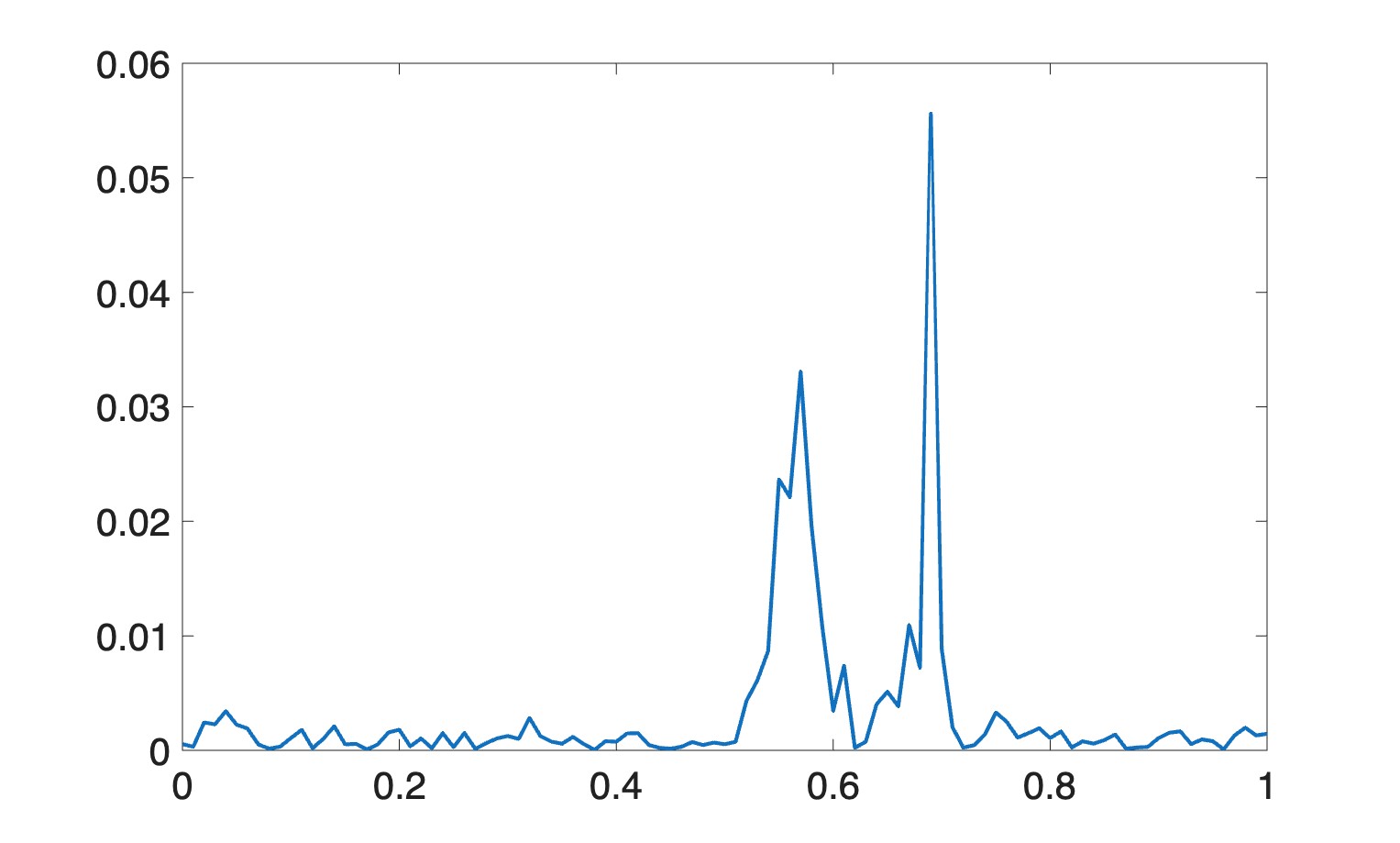}
\includegraphics[width=0.24\textwidth]{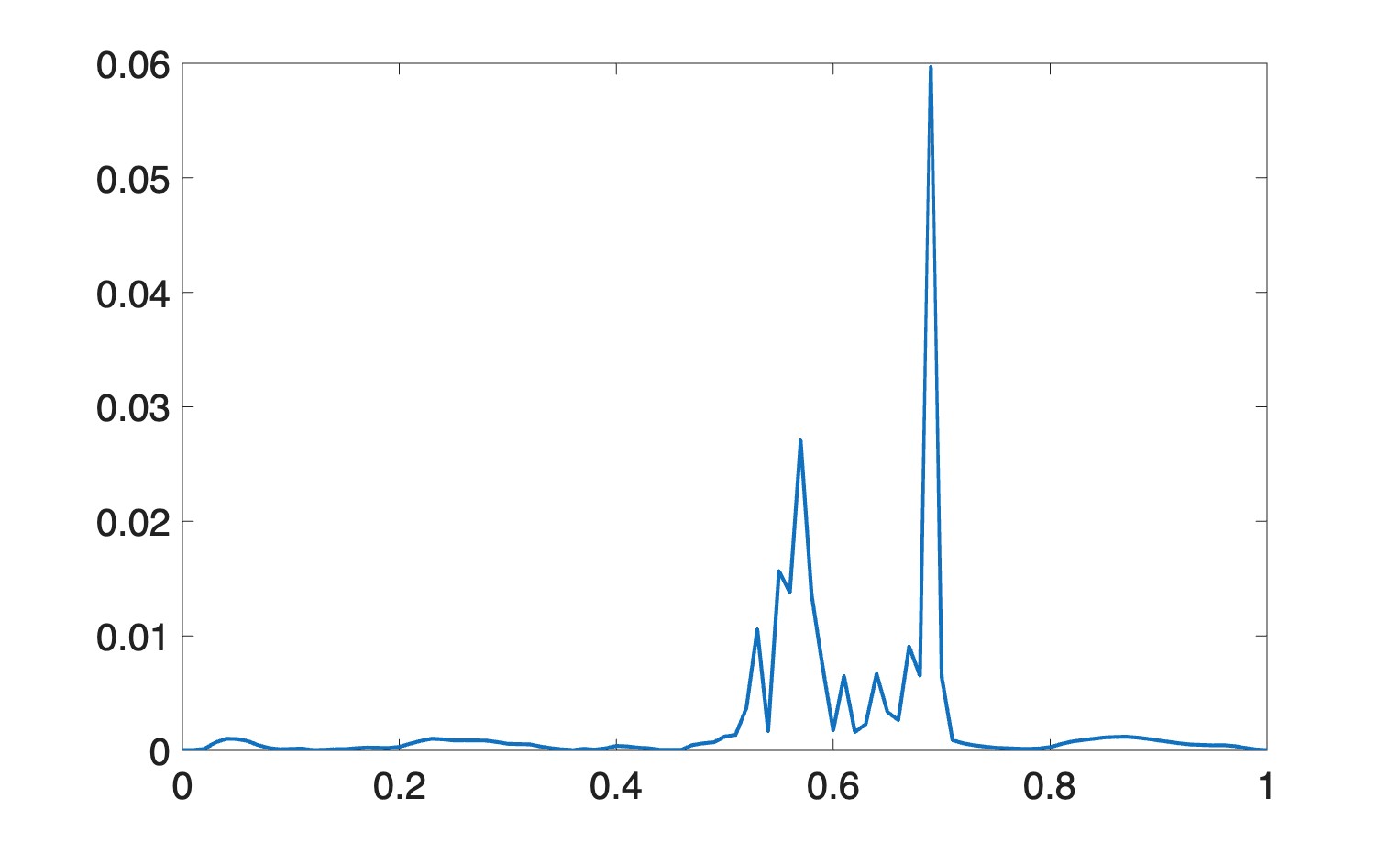}
\includegraphics[width=0.24\textwidth]{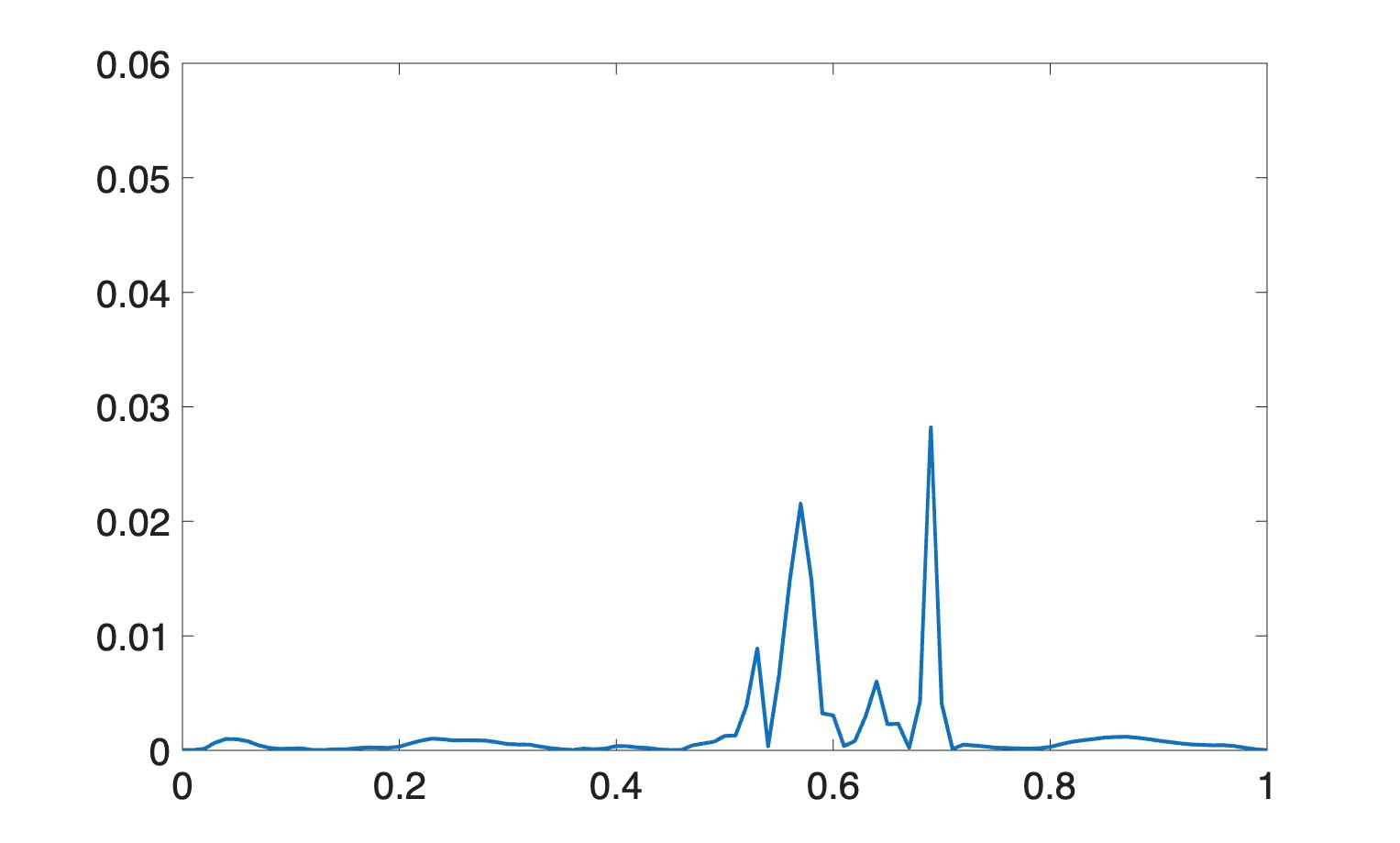}
\caption{Numerical solutions for the posterior mean solutions $\m$ (top) and the corresponding pointwise errors $err$ (bottom) for $W_C$, $\wh{W}_S$, and $W_S$ (left to right).}
\label{fig: burgers sparse oned}
\end{figure}

Figure \ref{fig: burgers sparse oned} shows the numerical solutions and the corresponding pointwise posterior errors along the one-dimensional cross section at $y=0.7$ at time $t=1$, obtained using $W_C$, $\wh{W}_S$ and $W_S$, respectively. It is evident that the results appear to be less oscillatory in smooth regions when using gradient-based weighting $\wh{W}_S$ and $W_S$ compared to covariance-based weighting $W_C$. Furthermore, incorporating correlation refinement helps to mitigate large errors near the discontinuities.

\begin{figure}[h]
\centering
\includegraphics[width=0.3\textwidth]{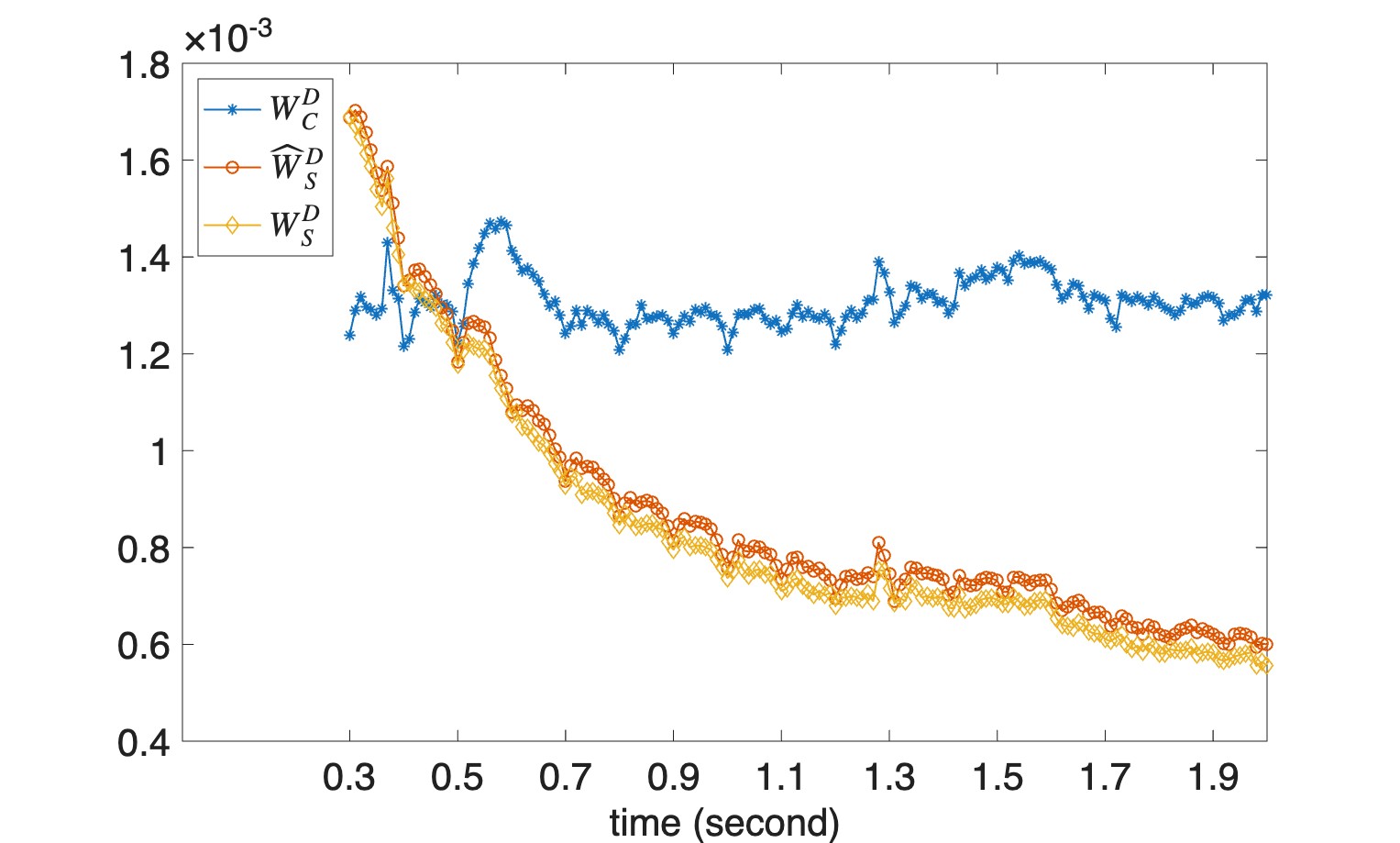}
\includegraphics[width=0.3\textwidth]{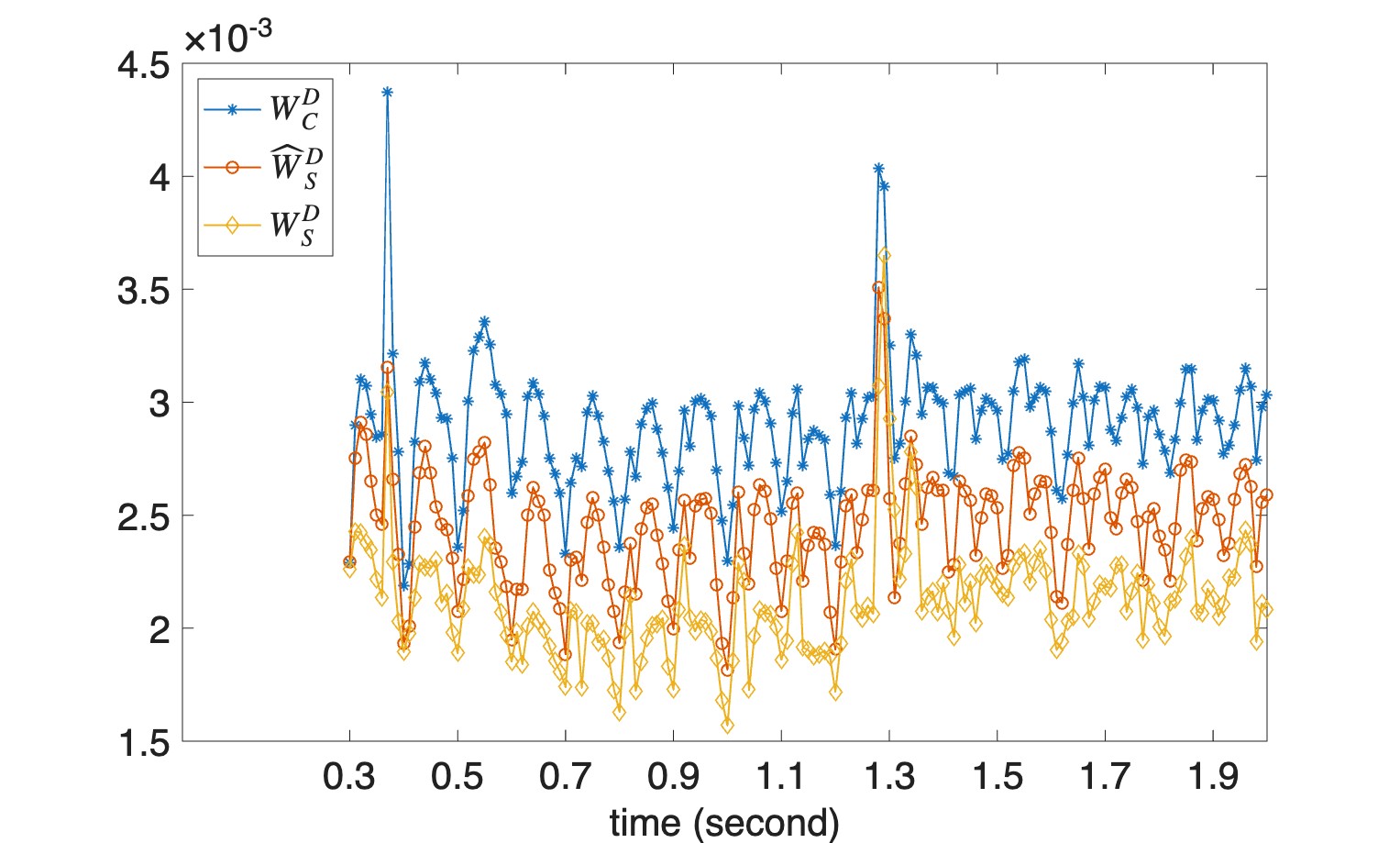}
\includegraphics[width=0.3\textwidth]{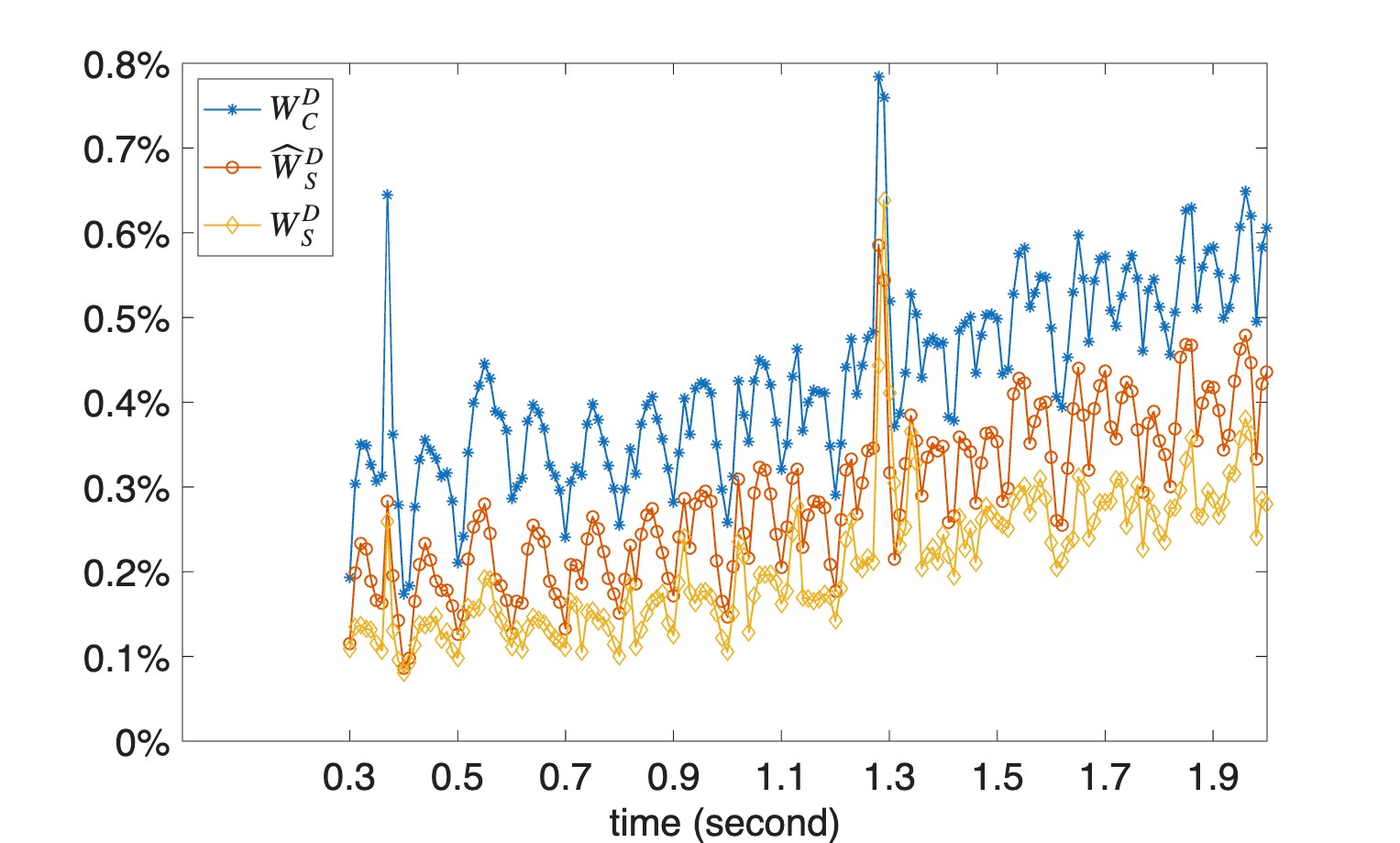}
\caption{The plots of the relative errors $err_{l_1}$  (left), $err_{l_2}$ (middle), and the complementary pattern correlation $1-Pcorr$ for $W_C$, $\wh{W}_S$, and $W_S$.}
\label{fig: burgers sparse err}
\end{figure}

Figure \ref{fig: burgers sparse err} displays the relative errors $err_{l_1}$, $err_{l_2}$, and the complementary pattern correlation $1-Pcorr$, obtained using \eqref{eq:l1err}, \eqref{eq:l2err}, and \eqref{eq:pcorr}, respectively, over the time domain $[0.3,2]$. The plots further verify  that the gradient-based weighting provides more accurate results, and that the correlation refinement further helps to reduce errors in the discontinuity regions, achieving better accuracy and higher pattern correlation. These results confirm that incorporating structural information through gradient-based weighting and directional refinement is essential for reliable assimilation in sparsely observed, discontinuous systems.

\section{Concluding Remarks}
\label{sec:summary}

We have developed a structurally informed data assimilation framework for two-dimensional systems characterized by piecewise-smooth or discontinuous states by designing a flexible gradient-based weighting matrix that is designed to accommodate both full and sparse data environments.  A novel directional refinement strategy is introduced to ensure that spurious correlations are suppressed across sharp transitions, especially in sparse data environments. Our numerical experiments demonstrate that using the structurally informed data assimilation framework consistently improves  
accuracy and preserves key structural features compared to standard ensemble-based techniques in both dense and sparse data environments.

While our current framework  uses  high-order finite difference schemes (WENO) for the prediction step, other numerical methods, such as discontinuous Galerkin or finite element methods may be more practical for problems in higher dimensions or on non-rectangular grids. Moreover, as our approach builds on a Bayesian interpretation of data assimilation,  it provides the foundation to rigorously study uncertainty quantification which will allow us to assess confidence in the assimilated estimates. Both ideas will be explored in future investigations, along will testing and adapting our approach to problems with more complex dynamics.

\bibliography{mybibfile}

\end{document}